\documentclass[10pt]{amsart}
\usepackage{amsthm,amsfonts,amsmath,amssymb,mathrsfs,verbatim,cite,graphicx,hyperref, eucal}

\newtheorem{Thm}{Theorem}[section]
\newtheorem{Def}[Thm]{Definition}
\newtheorem{Lm}[Thm]{Lemma}
\newtheorem{Prop}[Thm]{Proposition}
\newtheorem{Cor}[Thm]{Corollary}

\newtheorem{state}{Theorem}

\theoremstyle{definition}
\newtheorem*{ack}{Acknowledgements}

\theoremstyle{remark}
\newtheorem*{convention}{Convention}

\numberwithin{equation}{section}

\def\cal{\mathcal}
\def\Bbb{\mathbb}
\def\mf{\mathfrak}

\def\<{\langle}
\def\>{\rangle}
\def\dprime{{\prime\prime}}

\DeclareMathOperator{\hght}{ht}
\DeclareMathOperator{\wt}{wt}

\DeclareMathOperator{\sgn}{sgn}
\DeclareMathOperator{\ind}{ind}

\def\ag{{\bf a}}
\def\de{{\bf d}}

\def\a{\alpha}
\def\b{\beta}

\def\th{\theta}
\def\l{\lambda}

\def\e{\varepsilon}
\def\teps{\tilde{\e}}

\def\boldeta{\boldsymbol{\eta}}
\def\boldaleph{\boldsymbol \aleph}
\def\boldsigma{\boldsymbol{\sigma}}
\def\dotbvs{\dot{\boldsymbol{\varsigma}}}
\def\ddotbvs{\ddot{\boldsymbol{\varsigma}}}
\def\dotY{\dot{Y}}
\def\dotcalY{\dot{\mathcal{Y}}}
\def\ddotY{\ddot{Y}}
\def\ddotcalY{\ddot{\mathcal{Y}}}
\def\dotcalZ{\dot{\mathcal{Z}}}
\def\ddotcalZ{\ddot{\mathcal{Z}}}

\def\Re{\Bbb R}

\def\Z{\Bbb Z}

\def\A{\cal A}
\def\E{\cal E}
\def\Do{\cal D}

\def\h{\mf h}

\newcommand{\fcoef}[3]{\begin{matrix}\underline{#1}\!\!\!\!\! &&\!\!\!\!\! \underline{#2}\\  & #3 & \end{matrix}}
\newcommand{\scoef}[5]{\begin{matrix}#1\!\!\!\!\!&&\!\!\!\!\! \underline{#2}\!\!\!\!\! &&\!\!\!\!\! #3\\  &#4&& #5& \end{matrix}}
\newcommand{\pcoef}[2]{\left(#1\; ;#2\right)}
\newcommand{\ppcoef}[3]{\left(#1\; ;#2\; ;#3\right)}

\begin{document}
\title[]
{Generalized exponents of small representations. II.}
\author{Bogdan Ion}
\thanks{E-mail address: {\tt bion@pitt.edu}}
\date{April 15, 2009}
\address{Department of Mathematics, University of Pittsburgh, Pittsburgh, PA 15260}
\address{University of Bucharest, Faculty of Mathematics and Computer Science, Algebra and Number Theory research center, 14 Academiei St., Bucharest, Romania}


\maketitle

\section*{Introduction} 

This work is part of a series that aims to give manifestly non-negative formulas for generalized exponents, and in fact all $t$-weight multiplicities, corresponding to small weights.
It is a direct continuation of \cite{ion-geI} where the overall structure of the argument was illustrated for root systems of type $A$. The notation and conventions on finite root systems used here are those set up in \cite[$\S$ 1]{ion-geI}. 

The results contained in this paper correspond to those in \cite[$\S$ 3]{ion-geI} which form the backbone of the entire argument. More precisely, we give formulas for Fourier coefficients (or partition function coefficients) of the degenerate Cherednik kernel
$$
C(t)=\prod_{\a\in R^+} \frac{1-e^\a}{1-te^\a}
$$
Unlike the usual partition function, the answer in this case  encodes only the combinatorics of minimal expressions of a weight as a sum of roots and, in consequence, is amenable to a full description. Sections \ref{II-parametrization} and \ref{invariants} are devoted to developing the necessary foundations for the combinatorics of minimal expressions. The main results are contained in Section \ref{II-fourier} and Appendix \ref{defects} (in their explicit form).

Let $\l$ be a weight. A minimal expression for $\l$ is an expression of $\l$ as a sum of roots that contains the fewest possible number of terms. The number of terms in a minimal expression for $\l$ is called the co-length of $\l$ (denoted by $\ell^*(\l)$). For small dominant weights we specify a certain distinguished expression (and in fact a distinguished way of ordering the terms in the expression) which we call the canonical expression. It coincides to the canonical decomposition of dimension vectors for representations of quivers \cite[$\S$ 2.8]{kac}. Furthermore, on the set of roots contributing to a fixed expression we define an equivalence relation whose equivalence classes we call blocks. The canonical block decomposition of a small dominant weight is the block structure of the canonical expression of that weight. 

It turns out that recording the block sizes in canonical block decompositions gives a natural parametrization of small dominant weights that also contains information about the inclusion relation between the closures of the convex hulls of their Weyl group orbits and the canonical block decomposition  of weights in such a convex hull. Furthermore, together with some information about co-length two small weights it leads to a complete (weighted) enumeration of minimal expressions of any small weight. 

The symbol of $\l$ (denoted by $[\l]$) is a sequence of positive integers each representing the size of a block in the canonical block decomposition of $\l$. Sometimes these integers have to be adorned with extra information such as the length of the roots in particular blocks (if the root system is not simply laced) or the irreducible component of a root system (when reducible root systems have to be considered). Below, we give the details of how to recover a dominant weight from such information.  

\begin{convention}
Since the root systems of type $A_n$ where treated in \cite{ion-geI} and the non-zero small weights in type $G_2$ are the dominant roots, we henceforth assume that $R$ is a finite irreducible root system not of type $A_n$ or $G_2$. 
\end{convention}

Nevertheless, the root systems excluded above do not deviate from the behaviour we will describe. The convention is motivated by the fact that the analysis for $G_2$  is essentially empty and that for $A_n$ complicates the arguments by a significant factor.

Assume that $\Gamma$ is a parabolic root sub-system of $R$. We use the symbol $[1^{\Gamma, X}_{\ell}]$, and $[1^{\Gamma, X}_{s}]$, to refer to the dominant long and, respectively, short root of the irreducible component of $\Gamma$ of Dynkin type $X$. If there are several components with the same Dynkin type we might need to distinguish them further by adding 
subscripts $a, b, c,\dots$. It is important to note that the reference to length refers to the long and respectively short of roots  in $R$. For example, if $R$ is simply laced  we drop the reference to length from  the subscript but the reference to length remains relevant when $\Gamma$ is simply laced and $R$ is not simply laced.

 When $X$ refers to a simply laced Dynkin type, we denote by $[2^{\Gamma, X}]$ the dominant weight $2\th_{\Gamma,X}-\a_{\th_{\Gamma,X}}$, where $\theta_{\Gamma,X}$ is the dominant root in the irreducible component of $\Gamma$ of type $X$ and $\a_{\th_{\Gamma,X}}$ is a simple root in the indicated irreducible component that is not orthogonal on $\th_{\Gamma,X}$. When $\Gamma$ is of type $A$ there are two simple roots that are not orthogonal on the dominant root and one needs to specify further from which simple root $[2]$ is constructed. It turns out that only one of the choices produces a dominant small weight for $R$. This situation appears only for the root system  $E_6$ so we will prefer to  spell out these details when necessary rather than introduce more notation.

When $\Gamma$ is implicit we may drop it from the notation. When $\Gamma$ is irreducible we do not need to specify its Dynkin type and we may also drop it from the notation. For consistence we use $[\emptyset]$ to refer to the zero weight. 

By $\Gamma_{[1_s^X]}$, $\Gamma_{[1^X_\ell]}$, and $\Gamma_{[2^X]}$, we refer to the parabolic sub-systems of the irreducible component of $\Gamma$ of Dynkin type $X$ spanned by the simple roots that are orthogonal on the dominant short root, the dominant long root, and on both $\theta_{\Gamma,X}$ and $\theta_{\Gamma,X}-\a_{\th_{\Gamma,X}}$, respectively. 

Let $[a^{X_1}_1, a^{X_2}_2,\dots,a^{X_N}_N]$ be a sequence of symbols with $a_i$ from the set $\{1_s,1_\ell,2\}$ and $X_i$ Dynkin types. By convention,  the root system $\Gamma_{[\emptyset]}$ is $\Gamma$. Define inductively the root systems
\begin{equation*}
\Gamma_{[a^{X_1}_1, a^{X_2}_2,\dots,a^{X_N}_N]}:=\left({\Gamma_{[a^{X_1}_1, a^{X_2}_2,\dots,a^{X_{N-1}}_{N-1}]}}\right)_{[a_N^{X_N}]}
\end{equation*}

The same sequence $[a^{X_1}_1, a^{X_2}_2,\dots,a^{X_N}_N]$ will also  refer to the following weight of $R$
\begin{equation*}
\l:=\sum_{i=1}^N \left[a_i^{\Gamma_i,X_i}\right]
\end{equation*}
where $\Gamma_i$ above is an abbreviation for $\Gamma_{[a^{X_1}_1, a^{X_2}_2,\dots,a^{X_{i}}_{i}]}$. The symbol arising from the canonical block decomposition of  $\l$ turns out to be 
$[a^{X_1}_1, a^{X_2}_2,\dots,a^{X_N}_N]$. When $2$ does not appear in the sequence the construction of weights described here is known as the Kostant cascade construction. 
We are now ready to state our first result.

\begin{state}\label{II-state1} The symbols arising from the canonical block decomposition of non-zero small dominant weights of $R$ are the following
$$
\begin{aligned}
B_n:&  && [1_s^{B_n}]\\
        &  &&[1^{B_n}_\ell,1^{B_{n-2}}_\ell,\dots,1^{B_{n-2k}}_\ell],~ 0\leq k\leq (n-2)/2\\ 
        &  &&[1^{B_n}_\ell,1^{B_{n-2}}_\ell,\dots,1^{B_{n-2k}}_\ell, 1_s^{B_{n-2k-2}}],~ 0\leq k\leq (n-3)/2 \\
C_n:&  &&[1_\ell^{C_n}]\\
         &  &&[1_s^{C_n},1_s^{C_{n-2}},\dots,1_s^{C_{n-2k}}],~ 0\leq k\leq (n-2)/2\\
         &  &&[1_\ell^{C_n}, 1_s^{C_{n-1}},\dots, 1_s^{C_{n-2k-1}}],~ 0\leq k\leq (n-3)/2\\
D_n:&  &&[1^{D_n},1^{A_1}], [2^{D_n}]\\
         &  &&[1^{D_n},1^{D_{n-2}},\dots,1^{D_{n-2k}}],~ 0\leq k\leq (n-2)/2\\
         &  &&[2^{D_n},  1^{D_{n-3}},1^{D_{n-5}},\dots,1^{D_{n-2k-3}}],~ 0\leq k\leq (n-5)/2\\
E_6:&  &&[1^{E_6}], [1^{E_6},1^{A_5}], [2^{E_6}], [2^{E_6},1^{A_2}], [2^{E_6},2^{A_2}] \\
E_7:&  &&[1^{E_7}], [1^{E_7},1^{D_6}], [1^{E_7},1^{D_6},1^{D_4}], [2^{E_7}], [2^{E_7},1^{A_5}]\\
E_8:&  &&[1^{E_8}], [1^{E_8},1^{E_7}], [2^{E_8}], [2^{E_8},1^{E_6}]\\
F_4:&  && [1_s^{F_4}], [1^{F_4}_\ell], [1^{F_4}_\ell, 1^{C_3}_s]\\
G_2:& && [1_s^{G_2}], [1_\ell^{G_2}]
\end{aligned}
$$
If there are weights having the same symbol they are conjugate under the action of the automorphism group of the Dynkin diagram of $R$.
\end{state}
Note that by  \cite[(1.4)]{ion-geI}  the representations with highest weights that are conjugate under Dynkin diagram automorphisms have identical generalized exponents. 

Let us specify the symbols that parametrize more than one small dominant weight. In type $D_n$ with $n\geq 6$ even and $k=(n-2)/2$, we get $D_{n-2k}=A_1+A_1$ and therefore we have two weights parametrized by the symbol $[1^{D_n},1^{D_{n-2}},\dots,1^{D_4}, 1^{A_{1}}]$. If in addition $n=4$ there are three weights parametrized by $[1^{D_n},1^{A_1}]$ (the root system $R_{[1]}$ has three connected components, all of type $A_1$).  In type $D_n$ with $n\geq 5$ odd and $k=(n-5)/2$, we get $D_{n-2k-3}=A_1+A_1$ and therefore we have two weights parametrized by the symbol $[2^{D_n},  1^{D_{n-3}},1^{D_{n-5}},\dots,1^{D_4}, 1^{A_{1}}]$.

In type $E_6$, the root system $R_{[2]}$ has two connected components, both of type $A_2$ and hence there are two dominant weights parametrized by the symbol $[2^{E_6},1^{A_2}]$ and again two parametrized by the symbol $[2^{E_6},2^{A_2}]$. Regarding the latter case, there are two simple roots that are not orthogonal on the dominant root in $A_2$; the weight $[2^{A_2}]$ is the weight constructed from the simple root which is furthest away from the degree three node in the diagram of $E_6$. 

Let us briefly mention what is corresponding result for root systems of type $A_n$. A small dominant weight and its contragredient have the same symbol so it is enough to consider first layer weights. If $\l=(\l_1,\dots,\l_{n+1})$ is a first layer dominant weight, its symbol is $$[\l_1+1,\l_2+1,\dots,\l_{n+1}+1] $$  which is a partition of $n+1$ of length $n+1-\ell^*(\l)$. As it can be seen, the canonical block decomposition  can be a lot more complicated in type $A_n$ than in any other type. From this point of view, the symbol of a small dominant weight is a natural root system analogue of a partition.  

With the list available it makes sense to try to remove redundant information from our notation. To make the notation for symbols more compact we will write $a^k$ if $a$ appears $k$ times consecutively in $[\l]$. In a symbol $[a^{X_1}_1, a^{X_2}_2,\dots,a^{X_N}_N]$ we will drop any reference to the $X_i$'s unless it is absolutely necessary. By inspecting the list in Theorem \ref{II-state1} we observe that the only symbol that requires information about components is $[1^{D_n},1^{A_1}]$. In this case we will still choose drop the reference to roots systems but we will write $[1,1]$ rather than $[1^2]$ which refers to $[1^{D_n},1^{D_{n-2}}]$.

Using Theorem \ref{II-state1} it is rather easy to describe the cover relations for the partial order relation on dominant small weights given by the inclusion of closures of convex hulls of Weyl group orbits. This is done in Theorem \ref{II-thm3}. 

Furthermore, Theorem \ref{II-state1} leads to a computation of $\E_{[\l]}$, the number of minimal expressions of $\l$ weighted by the size of the blocks that appear in each minimal expression (see Section \ref{II-minimalexp} for definitions). It turns out that there are two kinds of weights. Some have the property that $\E_{[\l]}$ depends only on $\ell^*(\l)$ and are characterized by the fact that no 2 appears in $[\l]$. We call this normal (i.e. orthogonal) weights motivated by the fact that all roots in their canonical decomposition are mutually orthogonal. In type $A_n$ however, all small weights are normal. When referring to $\E_{[\l]}$ for normal weights, we can therefore drop the reference to length from the symbol of $[\l]$. For $\l$ of co-length at most 2 we denote
$$
\eta_{[\l]}:=\E_{[\l]}
$$
It is easy to see that $\eta_{[1]}=1$ and $\eta_{[2]}=2$. 
\begin{state}\label{II-state2}  With the notation above we have
\begin{enumerate}
\item[i)] $\E_{[1^{N}]}=\eta_{[1]}\eta_{[1^2]}\cdots\eta_{[1^{N]}}$,  where $\eta_{[1^N]}=(\eta_{[1^2]}-1)(N-1)+\eta_{[1]}$, $N\geq 1$.
\item[ii)] $\E_{[2,1^N]}=\eta_{[1]}\eta_{[1^2]}\cdots\eta_{[1^{N+1}]}\eta_{[2, 1^{N}]}$, where $\eta_{[2,1^N]}=(\eta_{[1^2]}-1)N+\eta_{[2]}$, $N\geq 0$.
\item[iii)] $\E_{[1,1]}=\eta_{[1]}\eta_{[1,1]}$.
\item[iv)] $\E_{[2^2]}=\eta_{[1]}\eta_{[1^2]}\eta_{[2,1]}\eta_{[2^2]}$, where $\eta_{[2^2]}=\eta_{[1^2]}\eta_{[2]}$.
\end{enumerate}
\end{state}
Therefore, the  enumeration of minimal expressions depends only on the enumeration for weights of co-length two (i.e. on the constants $\eta_{[1^2]}$, $\eta_{[1,1]}$, and $\eta_{[2]}$) which must be computed separately for each root system. Often, a more convenient invariant is $\delta=\eta_{[1^2]}-\eta_{[2]}$. 
The numerical values are given at the end of Section \ref{invariants}. In type $A_n$, $\eta_{[1^2]}=2$, $\delta=0$, and $\E_{[\l]}=\E_{[1^{\ell^*(\l)}]}=\ell^*(\l)!$

We are now ready to discuss our main result. We will use the conventions and notation spelled out in Section \ref{II-vectorconventions}. The Fourier coefficient $c_\l(t)$ of a small weight $\l$ is expressed in terms of the height $\hght(\l)$ of $\l$ and some data $\boldaleph_\l$ that contains information about the negative roots that appear in minimal expressions of $\l$. For normal weights, $\boldaleph_\l$ is a vector $\ag_\l\in \Z^{\ell^*(\l)}$, which in turn can be expressed as
$$
\ag_\l=\de_\l-\boldeta_{[\l]}
$$
with $\de_\l\in\Z_{\geq 0}^{\ell^*(\l)}$ (called defect vector) and  $\boldeta_{[1,1]}:=(\eta_{[1,1]},\eta_{[1]})$ and $$\boldeta_{[\l]}:=(\eta_{[1^{\ell^*(\l)}]},\dots,\eta_{[1]})$$  for all the other normal weights (called cut-off vector). The terminology is motivated by the fact that $c_\l(t)$ is zero as soon as one of the components of $\ag_\l$ becomes non-negative. The defect vector counts the number of negative roots in $R$ that participate in the minimal expressions of $\l$ (it is therefore the zero vector for dominant weights). In this case, $(1-t^{\boldaleph_\l})$ denotes 
$$
(1-t^{\ag_\l})
$$

For non-normal weights  other than those with symbols $[2^2]$ in type $E_6$ and $[2,1]$ in types $E_7$ and $E_8$, the data  $\boldaleph_\l$ is a 4-tuple $(\ag_\l,\tilde{a}_\l,\delta_\l,\bar{\ag}_\l)$, with $\ag_\l\in \Z^{\ell^*(\l)}$, $\tilde{a}_\l,\delta_\l\in\Z$, and $\bar{\ag}_\l\in\Z^{\ell^*(\l)-1}$ with each of the four parts being  a difference of vectors as for normal weights.  In this case, $(1-t^{\boldaleph_\l})$ denotes 
$$
(1-t^{\ag_\l})-t^{\tilde{a}_\l}(1-t^{\delta_\l})(1-t^{\bar{\ag}_\l})
$$

For the remaining weights  the data  $\boldaleph_\l$ is a 7-tuple  $(\ag_\l,\tilde{a}^\prime_\l,\delta^\prime_\l,\bar{\ag}^\prime_\l,\tilde{a}^\dprime_\l,\delta^\dprime_\l,\bar{\ag}^\dprime_\l)$ with $\ag_\l\in \Z^{\ell^*(\l)}$, $\tilde{a}^\prime_\l,\tilde{a}^\dprime_\l,\delta^\prime_\l,\delta^\dprime_\l\in\Z$, and $\bar{\ag}^\prime_\l,\bar{\ag}^\dprime_\l\in\Z^{\ell^*(\l)-1}$ satisfying similar properties. In this case, $(1-t^{\boldaleph_\l})$ denotes 
$$
(1-t^{\ag_\l})-t^{\tilde{a}_\l^\prime}(1-t^{\delta_\l^\prime})(1-t^{\bar{\ag}_\l^\prime})-t^{\tilde{a}_\l^\dprime}(1-t^{\delta_\l^\dprime})(1-t^{\bar{\ag}_\l^\dprime})
$$
We refer to Section \ref{II-nonnormal} for the details. The following statement sums up the results proved in Theorem \ref{II-thm6}, Theorem \ref{II-thm7}, Theorem \ref{II-thm10}, and Theorem \ref{II-thm8}.
\begin{state}\label{II-state3}
Let $\l$ be a small weight. Then,
$$
c_\l(t)=t^{\hght(\l)}(1-t^{\boldaleph_\l})
$$
\end{state}
The corresponding result in type $A_n$ is \cite[Theorem 3.5]{ion-geI}. The explicit definition of defect vectors is contained in Appendix \ref{defects}. For normal weights there should be a more canonical way to specify the defect vector by investigating some strongly regular graphs that are naturally constructed from the combinatorics of minimal expressions. For non-normal weights it is  arguable to what degree the defect vectors are canonical since  for such weights $(1-t^{\boldaleph_\l})$ can be written in several ways as a sum of two terms, respectively three terms. In Section \ref{II-symgroups} we define  symmetry groups and the notion of constrained orbit of a group. This is relevant for the understanding of the different expressions of $(1-t^{\boldaleph_\l})$ as a sum of two, and respectively three terms.

Below, $E(V_\l)$ are the generalized exponents, $\wt(\l)$ is the set of weights, and $m_{\l\mu}$ are the weight multiplicities of $V_\l$, the irreducible representation with highest weight $\l$. The following is the analogue of \cite[Theorem 3.7]{ion-geI} and an immediate consequence of Theorem \ref{II-state3}.
\begin{state}\label{II-state4}
Let $\l$ be a small dominant weight. Then,
$$
E(V_\l)=\sum_{{\mu\in\wt(\l)}}m_{\l\mu}t^{\hght(\mu)}(1-t^{\boldaleph_\mu})
$$
\end{state}
For the adjoint representation the above formula becomes the Shapiro-Steinberg formula for classical exponents. As revealed in \cite[$\S$ 4]{ion-geI} a further refinement  is possible. It is based  on the notion of quasisymmetric function which  naturally emerges from the combinatorics of the defect data.  We pursue this in  \cite{ion-geIII}.

\begin{ack}
 Work supported in part by the NSF grant DMS--0536962 and  the  CNCSIS grant nr. 24/28.09.07 (Groups, quantum groups, corings, and representation theory).  A substantial part of this work was carried out  while visiting the Max Planck Institute for Mathematics, Bonn during June--July 2008.  I thank MPIM for their hospitality and support. 
\end{ack}
\section{Parametrization}\label{II-parametrization}
\subsection{Notation}  We adopt the notation and conventions from \cite[$\S$ 1]{ion-geI}. In addition, for $S$ a set consisting of  dominant weights consider $$R_S:=\{\a\in R~|~(\l,\a)=0 \text{ for all } \l\in  S\}$$ It is well known that $R_S$ is root system that is spanned by the simple roots contained in $R_S$ which, in consequence, form a basis of $R_S$.
In general, $R_S$ is not irreducible. Pictorially, the irreducible components can be represented by the connected components of the graph obtained by 
removing from the Dynkin diagram of $R$ the nodes corresponding to simple roots that are not in $R_S$. 

For $S=\{\th_\ell\}$, or $S=\{\th_s\}$,   there is a unique simple root which is not in $R_S$. This simple root will be denoted by $\a_{\th_\ell}$ and,  respectively $\a_{\th_s}$.

\subsection{Co-length}
\begin{Def} Let $\l$ be a weight. 
\begin{enumerate}
\item[(a)] The {{co-length}} of $\l$, denoted $\ell^*(\l)$,  is  the minimal number of roots necessary to write $\l$ as a sum of roots. 
\item[(b)] A  { minimal expression} for $\l$ is a sum 
\begin{equation}\label{II-eq1}
\l=\sum_{i=1}^{\ell^*(\l)} \b_i
\end{equation}
where $\b_i$ are roots.  Two minimal expressions coincide if their terms coincide up to order.
\item[(c)] A root that appears in a minimal expression for $\l$ will be called $\l$-relevant. 
\item[(d)] A sub-expression of $\l$ is a partial sum in a minimal expression of $\l$. 
\end{enumerate}
\end{Def}

The co-length is constant on Weyl group orbits. A weight $\l$ has co-length zero if and only if it is the zero weight and it has co-length one if and only if it is a root.

Note that in general some terms of the sum \eqref{II-eq1} might appear with multiplicity.  We will show in Section \ref{II-multiplicities}  that for small weights any minimal expression must be multiplicity free.
\begin{Lm}\label{II-lemma1}
 Let $\l$ be a weight  and consider a fixed minimal expression for $\l$ as in \eqref{II-eq1}. Then, all the scalar products $(\b_i,\b_j)$ are non-negative. In particular, 
 \begin{equation}\label{II-eq2}
 (\l,\b^\vee)\geq 2
 \end{equation}
 for any $\l$-relevant root $\b$.
\end{Lm}
\begin{proof}
Indeed, assume that for some $i$ and $j$ we have $(\b_i,\b_j)<0$. Then, either 
$$
(\b_i^\vee,\b_j)=-1\quad \text{or}\quad (\b_i,\b_j^\vee)=-1 \quad \text{(or both)}
$$
Assuming that $(\b_i^\vee,\b_j)=-1$ we see that $$\b_i+\b_j=s_{\b_i}(\b_j) $$ which is again a root. Therefore, the expression \eqref{II-eq1} cannot be minimal.  

The remaining claim is an immediate consequence.
\end{proof}
\begin{Prop}\label{II-lemma2} Let $\l$ be a dominant weight.  Then, all the $\l$-relevant roots are positive.
\end{Prop}
\begin{proof}
Consider a minimal expression as in \eqref{II-eq1}. Then, by Lemma \ref{II-lemma1}, we have that $(\l,\b_i^\vee)\geq 2$ for any $1\leq i\leq \ell^*(\l)$. The weight $\l$ being dominant, the roots $\b_i$ must positive.
\end{proof}

\begin{Prop}\label{II-lemma10}
Let $\l$ and $\mu$ be  weights such that $\mu$ is a convex linear combination of  $\l$ and $s_\a(\l)$ for some root $\a$. Then,
$$
\ell^*(\mu)\leq \ell^*(\l)
$$
\end{Prop}
\begin{proof} If $\l$ and $\a$ are orthogonal the claim is trivial. By exchanging the role of $\l$ and $s_\a(\l)$ if needed, we may assume that $$k:=(\l,\a^\vee)>0$$
The weight $\mu$ must therefore be of the form
$$
\mu=\l-t\a
$$
for some integer $0\leq t\leq k$. Let 
$$
\l=\sum_{i=1}^{\ell^*(\l)} \b_i
$$
be a minimal expression for $\l$. By rearranging the terms in the sum we may assume that the first $s$ terms are those which have strictly positive scalar product with $\a$.
Therefore, $\mu$ can be written as 
$$
\sum_{i=1}^s (\b_i-j_i\a)
+\sum_{i=s+1}^{\ell^*(\l)}\b_i$$ for some integers $j_i$ such that $0\leq j_i\leq (\b_i,\a^\vee)$. But for all such integers $\b_i-j_i\a$ are roots so we managed to write $\mu$ as a sum of $\ell^*(\l)$ roots. The conclusion immediately follows.
\end{proof}


\subsection{Canonical expressions} 
\begin{Prop}\label{II-prop2} Let $\l$ be a dominant weight and let $\a$ be a 
$\l$-relevant root of smallest possible height. Then, $\l-\a$ is  dominant. 
\end{Prop}
\begin{proof}
Indeed, if there is a simple root $\a_i$ such that 
\begin{equation}\label{II-eq3}
(\l-\a,\a_i^\vee)<0
\end{equation}
then   $2\geq (\a,\a_i^\vee)\geq 1$ so $s_i(\a)$ is a root of height strictly smaller than that of $\a$. This forces $(\l,\a_i^\vee)$ to be 0 or 1. In the former
case, $s_i(\l)=\l$  so $s_i(\a)$ is $\l$-relevant which contradicts the choice of $\a$. In the latter case, we must have $(\a,\a_i^\vee)=2$.  Now $\a\neq\a_i$ because $\a$, being 
$\l$-relevant, must satisfy \eqref{II-eq2} which can be written as 
$$
(\l-\a,\a^\vee)\geq 0
$$
Therefore, $s_i(\a)=\a-2\a_i$ and $\a-\a_i$ is again a root.  If 
$$
\l=\a+\sum_{j=2}^{\ell^*(\l)}\b_j 
$$
is a minimal expression containing $\a$ then by applying $s_i$ to this equality we obtain
$$
\l-\a_i=(\a-2\a_i)+\sum_{j=2}^{\ell^*(\l)}s_i(\b_j) 
$$
In conclusion
$$
\l=(\a-\a_i)+\sum_{j=2}^{\ell^*(\l)}s_i(\b_j) 
$$
is another minimal expression for $\l$ which means that $\a-\a_i$ is $\l$-relevant and this contradicts the choice of $\a$. Hence, there is no simple root satisfying \eqref{II-eq3}.
\end{proof}

\begin{Def}
 A canonical expression of a dominant weight is defined  inductively on co-length as follows.
\begin{enumerate}
\item[(a)] The canonical expression of the zero weight is by definition  $\l=0$.
\item[(b)] Let $\l$ be a dominant weight of co-length $N\geq 1$. The minimal expression $$
\l=\b_1+\cdots+\b_N
$$ is a canonical expression of $\l$ if $\b_N$ is a $\l$-relevant root  of smallest possible height and $$\l-\b_N =\b_1+\cdots+\b_{N-1}$$ is a canonical expression of  $\l-\b_N$.  
\end{enumerate}
\end{Def}

Note that by Proposition \ref{II-prop2}, the weight $\l-\b_N$ is a  dominant weight of co-length $N-1$. 
Proposition \ref{II-prop2} assures in fact that any dominant weight has at least one canonical expression. The fact that small dominant weights have exactly one canonical expression is a consequence of the unicity of their canonical block decomposition which is proved in Section \ref{thelist}. Until then, remark that the unicity can be easily established for small dominant weights of co-length two. Indeed, if $\l$ is a small dominant weight of co-length two then the first term in a canonical expression must be a dominant root. In case there are two dominant roots that are $\l$-relevant then the longest one must participate in the canonical expression. This uniquely identifies the first term of the canonical expression and consequently the entire expression. 


\subsection{Multiplicities}\label{II-multiplicities}
\begin{Prop}\label{II-prop4}
Let $\l$ be a small weight. Then any sub-expression of $\l$ is  small.
\end{Prop}
\begin{proof} It is safe to assume that $\l$ is dominant. We prove by induction on $\ell^*(\l)$ that any sub-expression of $\l$ is small.

The claim is easily verified for dominant weights of co-length one. We assume that $\ell^*(\l)\geq 2$ and that the claim is true for all small dominant weights of strictly smaller co-length.
By making use of the induction hypothesis it is enough to show that any sub-expression of $\l$ of co-length $\ell^*(\l)-1$ is small. Let 
$$
\l=\sum_{i=1}^{\ell^*(\l)} \b_i
$$
be a minimal expression for $\l$ and let $j$ be an arbitrary element of $[\ell^*(\l)]$. We will show that 
$$
\l-\b_j=\sum_{i\neq j} \b_i
$$
lies in the convex hull of the $W$ orbit of $\l$ and, in consequence, it is a small weight.

The integer $k:=(\l,\b_j^\vee)$ is at least 2
as guaranteed by Lemma \ref{II-lemma1}. Then, $$\l-\b_j=(1-\frac{1}{k})\l+\frac{1}{k}s_{\b_j}(\l)$$
is a weight  which is a nontrivial convex linear combination of $\l$ and $s_{\b_j}(\l)$. Therefore, $\l-\b_j$ lies in the convex hull of the $W$ orbit of $\l$.
\end{proof}

\begin{Prop}\label{II-prop1}
Let $\l$ be a small weight. Then, any minimal expression of $\l$ is multiplicity free.
\end{Prop}
\begin{proof} The claim is an immediate consequence of the previous Proposition. Indeed, let
$$
\l=\sum_{i=1}^{r} m_i\b_i
$$
 a minimal expression for $\l$ such that the roots $\b_i$, $1\leq i\leq r$ are  distinct and the positive integers $m_i$ are their multiplicities. Then every partial sum $m_i\b_i$ is a sub-expression of $\l$ and therefore small. This immediately implies that all the multiplicities $m_i$ are necessarily equal to one.
\end{proof}

\subsection{Small weights of co-length two} The small weights of co-length two play an especially important role in the forthcoming arguments. We state here some facts about 
their possible minimal expressions. Full proofs can be found in Appendix \ref{appA}.

\begin{Def} Let  $\l$ be a small dominant weight of co-length two such that all its minimal expressions consist of orthogonal roots. Define $\eta_2(\l)$ to be the total number of minimal expressions of $\l$. 
\end{Def}

Using Lemma \ref{II-lemma9} it is very easy to list the values of $\eta_2(\l)$ for all weights $\l$ that satisfy the conditions (they are as many as connected components of $R_\th$). 
\begin{description}
\item[$A_n$] $\eta_2(\l)=2$
\item[$B_n$] $\eta_2(\l)=3$
\item[$C_n$] $\eta_2(\l)=3$
\item[$D_n$] $\eta_2(\l)\in\{3, n-1\}$
\item[$E_6$] $\eta_2(\l)=4$
\item[$E_7$] $\eta_2(\l)=5$
\item[$E_8$] $\eta_2(\l)=7$
\end{description}
The values were also computed in \cite[Section 5.6]{ion-ge}. For $F_4$ there are no such weights. What is important for us is that in all cases $\eta_2(\l)\geq 3$ and this is one crucial aspect that differentiates root systems of type $A_n$ from the other roots systems. As we will see, this fact forces blocks (defined in Section \ref{blocks}) to have rather small size, as opposed to the type $A$ case where they can have the size as large as the rank.

\begin{Prop}\label{II-prop3}
Let $\l$ be a small weight of co-length two. 
\begin{enumerate}
\item[(a)]  If $R$ is simply laced and $\l$ has a minimal expression consisting of non-orthogonal roots then this is the unique minimal expression of $\l$. If $\l$ is in addition dominant then this unique minimal expression is $$\l=\th+(\th-\a_\th)$$
\item[(b)] If $R$ is not simply laced and $\l$ has a minimal expression consisting of non-orthogonal roots then this expression is the unique one consisting of non-orthogonal roots, it is not canonical, its constituents are short roots, and $\l$ has other minimal expressions, each consisting of orthogonal roots of different length. If $\l$ is in addition dominant then $$\l=\th_s+(\th_s-\a_{\th_s})$$ Its canonical expression is $$\l=\th_\ell+(2\th_s-\th_\ell-\a_{\th_s})$$ and any root having positive scalar product with both roots in the canonical expression is $\l$-relevant.
\item[(c)] If  all minimal expressions of $\l$ consist of orthogonal roots then the roots that have positive scalar product with both terms in a fixed minimal expression are $\l$-relevant.
Moreover, $\eta_2(\l)$ is at least 3. If $\l$ is in addition dominant then its canonical expression is $$\l=\th_1+\th_2$$ where $\th_1$ is a dominant root in $R$ and $\th_2$ is a dominant root in $R_{\th_1}$ whose length does not exceed the length of $\th_1$.
\end{enumerate}
\end{Prop}


\subsection{Blocks}\label{blocks}  Let $\l$ be a small weight and 
$$
\l=\sum_{i=1}^{\ell^*(\l)}\b_i 
$$
a minimal expression for $\l$. We define an equivalence relation $\not\perp$ on the set of roots that appear in the above expression as follows
$$
\b_i\not\perp\b_j\quad \text{if and only if}\quad (\b_i,\b_j)>0
$$
The only thing that needs justification is transitivity. Assume that $\b_i\not\perp\b_j$ and $\b_j\not\perp\b_k$. From Lemma \ref{II-lemma1} the scalar product between $\b_i$ and $\b_k$ is non-negative. If $\b_i$ and $\b_k$ are orthogonal then, by Proposition \ref{II-prop3}{(c)}
$$
\b_j+(\b_i+\b_k-\b_j)
$$
is a minimal decomposition of $\b_i+\b_k$. Therefore,
$$
\l=\b_j+(\b_i+\b_k-\b_j)+\sum_{s\neq i,k}\b_s 
$$
is a minimal expression of $\l$ in which $\b_j$ appears with multiplicity $2$; this is in contradiction with Proposition \ref{II-prop1}. In consequence, $\b_i\not\perp\b_k$ and $\not\perp$ is an equivalence relation. It is worth mentioning that from Proposition \ref{II-prop3}{(a,b)} the scalar products $(\b_i,\b_j^\vee)$ between equivalent roots are all 1.
\begin{Def} Let $\l$ be a small weight,
$$
\l=\sum_{i=1}^{\ell^*(\l)}\b_i 
$$
a minimal expression for $\l$, and $\not\perp$ the equivalence relation on the set of roots that appear in this expression, defined as above.
The equivalence classes with respect to this equivalence relation will be called blocks.  A $\l$-block is a block for some minimal expression of $\l$.
\end{Def}

If $B$ is a $\l$-block, we will use $\{B\}$ and $B$ to refer to the block as a set of roots and, respectively, the weight 
$$\sum_{\a\in \{B\}}\a$$ 
For example, the number of elements in $\{B\}$ is $\ell^*(B)$ and the sum of the heights the elements in $\{B\}$ is $\hght(B)$.


\subsection{Canonical block decompositions}

\begin{Prop}\label{II-prop5} Let $\l$ be a small dominant weight and denote by  $B$ a $\l$-block of smallest possible height. Then, $\l-B$ is  dominant. 
\end{Prop}
\begin{proof} If $\{B\}$ has only one element then this is exactly Proposition \ref{II-prop2}. Hence, we may assume that $\{B\}$ has more than one element. From Proposition \ref{II-prop3}{(a,b)} we know that either $R$ is simply laced or $R$ is not simply laced and the elements of $\{B\}$ are short roots. What is important for us is that if 
$$
(\b,\a^\vee)=\pm 2
$$
for $\b$ an element of $\{B\}$ and $\a$ an arbitrary root then $\b=\pm \a$.

Assume that there is a simple root $\a_i$ such that 
\begin{equation}\label{II-eq4}
(\l-B,\a_i^\vee)\leq -1
\end{equation}
With the notation $s:=(\l,\a_i^\vee)$, $t:=(B,\a_i^\vee)$, this can be restated as 
$$
t\geq s+1\geq 1
$$
Let $$\l=B_1+\cdots+B_k+B$$
be the block decomposition of a minimal expression that contains $B$. We claim that
\begin{equation}\label{II-eq5}
\l=s_i(B_1)+\cdots+s_i(B_k)+(s_i(B)+s\a_i)
\end{equation}
is another minimal expression for $\l$. Since $\ell^*(B_j)=\ell^*(s_i(B_j))$ for any $1\leq j\leq k$ and $$\ell^*(B-(t-s)\a_i)\leq \ell^*(B) $$ as assured by Lemma \ref{II-lemma10},
we see that \eqref{II-eq5} is indeed a minimal expression for $\l$.

The elements of each $s_i(B_j)$, $1\leq j\leq k$, stay in the same equivalence class of the minimal expression \eqref{II-eq5}. Since
$$
s_i(B)+s\a_i=B-(t-s)\a_i
$$ 
its height is strictly smaller than that of $B$. All the elements in any minimal expression of $s_i(B)+s\a_i$ are therefore forced to be equivalent to elements of $s_i(B_j)$, $1\leq j\leq k$, since the ones that would not satisfy this property would form  equivalence classes which have smaller height that of $B$.

Let $\{B^\prime\}$ and $\{B^{\dprime}\}$ be non-empty complementary  subsets of $\{B\}$. We keep the same convention for subsets of $\{B\}$  that we use for blocks, leaving out the parentheses to refer to the sum of elements. If 
$$
(B^\prime,\a_i^\vee)\geq s
$$ then the weight $B^\prime-s\a_i$ is a convex linear combination of $B^\prime$  and $s_i(B^\prime)$ and by Lemma \ref{II-lemma10} 
$$
\ell^*(s_i(B^\prime)+s\a_i)\leq \ell^*(s_i(B^\prime))
$$
In consequence,
$$
\ell^*(s_i(B^\prime)+s\a_i)+\ell^*(s_i(B^{\prime\prime}))\leq \ell^*(s_i(B^\prime))+\ell^*(s_i(B^{\prime\prime}))=\ell^*(s_i(B))
$$
which forces equality in the above equation. This simply means that if we have a minimal expression for $s_i(B^\prime)+s\a_i$ and a minimal expression for $s_i(B^{\prime\prime})$, we obtain a minimal expression for $s_i(B)+s\a_i$ merely by summing these two expression. In such a scenario the elements of $s_i(B^{\prime\prime})$ are not equivalent to any element of $s_i(B_j)$ for any $1\leq j\leq k$, which is not acceptable. From this discussion we conclude that 
$$
(B^\prime,\a_i^\vee)< s
$$
for any subset $\{B^\prime\}$ of $\{B\}$.

Now, let $\b$ be an element of $\{B\}$, $\{B^{\dprime}\}=\{\b\}$ and let $\{B^\prime\}$ be the complement of $\{B^{\dprime}\}$ inside $\{B\}$. By the above considerations
$$
s>(B^\prime,\a_i^\vee)=t-(\b,\a_i^\vee)\geq s+1-(\b,\a_i^\vee)
$$
which means that $(\b,\a_i^\vee)=2$ and, as remarked at the beginning of the proof, $\b=\a_i$. As the chosen element of $\{B\}$ was arbitrary and $\{B\}$ has at least two elements we obtain 
from Proposition \ref{II-prop4} that $\l$ is not small. This contradicts the hypothesis. The assumption \eqref{II-eq4} must be false and the conclusion follows.
\end{proof}

\begin{Def}
 A canonical block decomposition of a small dominant weight is defined  inductively on co-length as follows.
\begin{enumerate}
\item[(a)] The canonical block decomposition of the zero weight is by definition  $\l=0$.
\item[(b)] Let $\l$ be a dominant weight of co-length $N\geq 1$. The block decomposition $$
\l=B_1+\cdots+B_N
$$ is a canonical block decomposition of $\l$ if $B_N$ is a $\l$-relevant block of smallest possible height (among all  $\l$-blocks) and $$\l-B_N =B_1+\cdots+B_{N-1}$$ is a canonical block decomposition of  $\l-\b_N$.  
\end{enumerate}
\end{Def}
The following result is a straightforward consequence of Proposition \ref{II-prop3}.
\begin{Lm}\label{II-lemma19}
Let $\l$ be a small dominant weight of co-length 2. Then a canonical block decomposition of $\l$ is the block decomposition of a canonical expression. In particular, it is unique.
\end{Lm}
\subsection{Blocks of size three} In this section we show that a small weight does not have any blocks of size three. Here our convention on the root system $R$ is especially 
important since the root systems of type $A_n$ are the only ones for which a small weight can have blocks of arbitrary size (less than $n$). The reason is different for simply laced and non-simply laced root systems and we treat them separately.
\begin{Lm}\label{II-lemma11}
Assume that $R$ is not simply laced and let $\l$ be a small weight of co-length 3. Then there is no $\l$-block of size three.  
\end{Lm}
\begin{proof} It is enough to prove the claim for dominant weights. Assume that $\l$ has a minimal expression 
\begin{equation}\label{II-eq6}
\l=\b_1+\b_2+\b_3
\end{equation}
with only one block. From Proposition \ref{II-prop3}{(b)} we know that the weight $\b_1+\b_2$ has at least one other minimal expression 
$$
\b_1+\b_2=\a+\b_4
$$
consisting of orthogonal roots of different length. We denoted by $\a$ the long root. From \eqref{II-eq6}  we obtain that $$(\l,\b_3^\vee)=4$$
Hence,
$$
4=(\a,\b_3^\vee)+(\b_4,\b_3^\vee)+2
$$
Analysing the possible values of the scalar product it is easy to see that we must have 
$$
(\a,\b_3^\vee)=2\quad \text{and} \quad (\b_4,\b_3^\vee)=0
$$
In other words, the set $\{\a,\b_3\}$ is a block of size two of the minimal expression
$$
\l=\a+\b_4+\b_3
$$
In particular, $\a+\b_3$ must be a small weight of co-length two but the lengths of $\a$ and $\b_3$ do not satisfy the conditions in Proposition \ref{II-prop3}{(b)}. Our assumption  on the existence of a $\l$-block of size three is therefore false.
\end{proof}
\begin{Lm}\label{II-lemma12}
Assume that $R$ is simply laced and let $\l$ be a small weight of co-length 3. Then there is no $\l$-block of size three.  
\end{Lm}
\begin{proof} It is enough to prove the claim for dominant weights. Assume that $\l$ has a minimal expression 
\begin{equation}\label{II-eq7}
\l=\b_1+\b_2+\b_3
\end{equation}
with a single block. We may assume that the terms of the sum are written in the descending order of their heights. Let us remark first that this must be the unique minimal expression for $\l$. Indeed, if 
\begin{equation}\label{II-eq8}
\l=\b_4+\b_5+\b_6
\end{equation}
is another minimal expression, then
$$
4=(\l,\b_1)=(\b_4,\b_1)+(\b_5,\b_1)+(\b_6,\b_1)
$$
and at least one of the scalar product on the right hand side must equal two. Hence, the root $\b_1$ appears also in the second expression and by keeping in mind that from Proposition \ref{II-prop3}{(a)} the weight $\-\b_1$ has a unique minimal expression we deduce that the expressions \eqref{II-eq7} and \eqref{II-eq8} coincide.

The expression \eqref{II-eq7} being the unique minimal expression for $\l$ is a canonical expression. Again, from Proposition \ref{II-prop3}{(a)} we know that
$$
\b_1=\th\quad \text{and}\quad \b_2=\th-\a_\th
$$
The fact that the minimal expression has a single block allows the third term to be written as
$$
\b_3=\th-\a_\th-\a
$$
with $\a$ a positive root satisfying
$$
(\a,\th)=0\quad\text{and}\quad (\a,\a_\th)=-1
$$
Now, we remark that the weight $\l-2\th$ must have co-length at least 2, otherwise $\l$ would not be small. Also
\begin{eqnarray*}
\l&=& -\a_\th+(\th-\a_\th-\a)\\
&=& -(\a_\th+\a)+(\th-\a_\th)
\end{eqnarray*}
are two expressions of $\l$ as  sum of two orthogonal roots, so $\ell^*(\l-2\th)=2$. 

We claim that $-\a_\th$ and $-(\a_\th+\a)$ are the only $(\l-2\th)$-relevant negative roots. Indeed, assume that $\b$ is a $(\l-2\th)$-relevant negative root. Applying Lemma \ref{II-lemma8} for both expressions we find that either $\b$ appears in one of these expressions or
$$
(\b,-\a_\th)=(\b,\th-\a_\th-\a)=1\quad\text{and}\quad (\b,-\a_\th-\a)=(\b,\th-\a_\th)=1
$$
In the latter case, the above equalities imply that 
$$
(\b,2\th-\a_\th)=2
$$
But, $2\th-\a_\th$ is a dominant weight (again, Proposition \ref{II-prop3}{(a)}) and $\b$ is a negative root. We arrived at a contradiction so $\b$ appears necessarily in  one of the two indicated minimal expressions of $\l-2\th$.

From Proposition \ref{II-prop3}{(c)} we know that there must be at least one other minimal expression of $\l-2\th$ and, from the above considerations, this minimal expression consists of positive roots. 

We have proved that $\l-2\th$ is a sum of two positive roots which means that $\l$ is not small, in contradiction with the hypothesis. In conclusion, there is no $\l$-block of size three. 
\end{proof}
\begin{Prop}\label{II-prop6}
Small weights  have $\l$-blocks of size at most two.
\end{Prop}
\begin{proof} Immediate from Lemma \ref{II-lemma11}, Lemma \ref{II-lemma12}, and Proposition \ref{II-prop4}.
\end{proof}

\subsection{Further constraints}\label{constraints} We collect here some constraints on the block sizes of a small weight as well as on their distribution in a canonical block decomposition of a small dominant weight. 
\begin{Lm}\label{II-lemma13} Assume that $R$ is not simply laced and let $\l$ be a small weight. Then, any block decomposition  of $\l$ has at most one block of size two.
\end{Lm}
\begin{proof} The claim can be checked directly.  In type $B$ all short roots are orthogonal so, by Proposition \ref{II-prop3}{(b)},  all blocks must have one element. In type $C$ a block of size two can be written as a sum between a long root and a short root, and the sum of any two long roots is twice a root hence there cannot be two blocks with two elements. 

In type $F$, assume that $\l$ is a small weight of co-length 4, that has a minimal expression 
$$
\l=B_1+B_2
$$
consisting of two blocks, each with two elements. By replacing $\l$ with an element in its orbit we may assume that $B_1$ is a dominant weight and hence, by Proposition \ref{II-prop3}{(b)}
$$
B_1=\th_s+(\th_s-\a_{\th_s})
$$
The weight $B_2$ is a weight in the root system $R_{\{\th_s, B_1\}}$ which, in this case, is the root system of type $A_2$ spanned by the two simple long roots in $R$. From Lemma \ref{II-lemma5} we know that $B_2$ cannot be a small weight, contradicting the smallness of $\l$.
\end{proof}
The case of simply laced root systems is sightly more delicate. Assume that $R$ is simply laced and denote by $\a_j$ the unique simple root that is not orthogonal on $\th$. 
Consider now the root system $R_{\{\th, 2\th-\a_j\}}$ and denote by $I$ the set of nodes of the Dynkin diagram of $R$ corresponding to $\a_j$ and its neighbours. 
 It is clear that the nodes  corresponding to the basis of $R_{\{\th, 2\th-\a_j\}}$ are exactly those in the complement of $I$.  

Let $\th^\prime$ be the dominant root in an irreducible component of $R_{\{\th, 2\th-\a_j\}}$, and assume that $\a_k$ is a simple root in $R_{\{\th, 2\th-\a_j\}}$ that is not orthogonal on $\th^\prime$ and orthogonal on all simple roots indexed by elements of $I$. 
\begin{Lm}\label{II-lemma14} With the notation above, the weight 
$$
\l=(2\th-\a_j)+(2\th^\prime-\a_k)
$$
is not small.
\end{Lm}
\begin{proof} Remark first that $\l$ is dominant and that it can be written as a sum of four roots
\begin{equation}\label{II-eq9}
\l=\th+(\th-\a_j)+\th^\prime+(\th^\prime-\a_k)
\end{equation}
We claim that $\l$ has co-length four so the above expression is a minimal expression for $\l$. Indeed, consider 
$$
\l=B_1+B_2+\cdots
$$
a canonical block decomposition for $\l$. If $\{B_1\}$ has one element, then $$B_1=\th\quad \text{and}\quad (\l,\th)=2$$ But, from \eqref{II-eq9} we must have $(\l,\th)=3$ which contradicts the above equality. Therefore, $\{B_1\}$ must have size two and, consequently,
$$
B_1=\th+(\th-\a_j)
$$
Now,
$$
\l-B_1=\th^\prime+(\th^\prime-\a_k)
$$
which clearly has co-length two. The claim is therefore proved.

Let $\a_t$ be the simple root defined as the neighbour of $\a_j$ that is closest to $\a_k$ in the Dynkin diagram. It should be noted that 
$$
(\th^\prime,\a_t)=-1
$$
otherwise $\th^\prime$ would be dominant in a root system that contains  $R_{\{\th, 2\th-\a_j\}}$ strictly.
Let
$$
\a=\th-\a_j-\a_t=-s_\th s_j(\a_t)
$$
Then,
$$
(\th,\a)=(\th^\prime,\a)=(\th-\a_j,\a)=(\th^\prime-\a_k,\a)=1
$$
From, Proposition \ref{II-prop3}{(c)} we know that $\a$ is both $(\th+\th^\prime)$-relevant and $(\th-\a_j+\th^\prime-\a_k)$-relevant. In consequence,
$$
\l=2\a+(\th+\th^\prime-\a)+(\th-\a_j+\th^\prime-\a_k-\a)
$$
is a minimal expression of $\l$ which is not multiplicity free and therefore $\l$ is not small.
\end{proof}
\begin{Cor}\label{II-cor3}
Assume that $R$ is simply laced. A small weight has a minimal expression with two blocks of size two only if $R$ is of type $E_6$.
\end{Cor}
\begin{proof} Let $\l$ be a small weight which has a minimal expression with two blocks of size two. By considering the weight given by the sum of the the two blocks of size two we may assume that $\l$ has co-length four. Moreover, by replacing $\l$ with an element in its orbit we may assume that the first block is dominant (and hence it equals $2\th-\a_{\th_s}$) and the second is dominant in an irreducible component of  $R_{\{\th, 2\th-\a_{\th_s}\}}$. It is  straightforward check that the hypotheses of Lemma \ref{II-lemma14} are satisfied for $R$ of type $D_n$, $E_7$ and $E_8$. Therefore, $R$ is forced to be of type $E_6$. 
\end{proof}
We should stress the if $R$ is of type $E_6$ then there is indeed a small dominant weight of length four which has a minimal expression consisting of two blocks. With the notation as in the discussion before Lemma \ref{II-lemma14}, the root system $R_{\{\th, 2\th-\a_j\}}$ has two irreducible components, each of type $A_2$. The root $\th^\prime$ is dominant in one of the components and as such it is not orthogonal on any of the two simple roots which form the basis.  One of the two simple roots, let us call it $\a_l$, does not satisfy the conditions required for $\a_k$ and 
$$
2\th-\a_j+2\th^\prime-\a_l
$$
is a small dominant weight of length four which has a minimal expression consisting of two blocks. The other dominant weight, constructed from the second irreducible component of  $R_{\{\th, 2\th-\a_j\}}$ is in fact the image of the first one via the non-trivial automorphism of the Dynkin diagram. 

Before moving on to revisit the canonical block decomposition of small dominant weights we need one last, rather empirical, observation. Assume that $R$ is simply laced, consider the root system $R_{\th}$, and  let $\th^\prime$ be the dominant root in an irreducible component  of rank at least two of $R_{\th}$. Denote by $\a_k$  a simple root in $R_{\th}$ that is not orthogonal on $\th^\prime$. 
\begin{Lm}\label{II-lemma15} With the notation above, the weight 
$$
\l=\th+(2\th^\prime-\a_k)
$$
is not dominant.
\end{Lm}
\begin{proof} We observe that $\a_{\th}$ and $\a_k$ are in fact orthogonal. Also, $(\th^\prime,\a_{\th})=-1$ otherwise $\th^\prime$ would be dominant in $R$. 
Now,
$$
(\l,\a_{\th})=-1
$$
so $\l$ is not dominant.
\end{proof}
\begin{Lm}\label{II-lemma17}
Assume that $R$ is simply laced, let $\l$ be a small dominant weight and
\begin{equation}\label{II-eq11}
\l=B_1+\cdots+B_k
\end{equation}
a canonical block decomposition. The size of the blocks that appear is decreasing.
\end{Lm} 
\begin{proof} If there is a block of size one (say $B_i$) followed by a block of size two then apply Lemma \ref{II-lemma15} for $B_i+B_{i+1}$ which is dominant in $R_{\{B_1\}\cup\dots \{B_{i-1}\}}$.
\end{proof}
\begin{Lm}\label{II-lemma20}
Assume that $R$ is simply laced, let $\l$ be a small  weight which does not have more than two blocks of size two in any minimal expression. Then any minimal expression of $\l$ has the same number of blocks of size two.
\end{Lm} 
\begin{proof} The claim obvious for weights of co-length one and two. In what follows we assume that $\l$ has co-length at least three. Let \begin{eqnarray*}
\l&=& \b_1+\cdots+\b_N\\
&=&  \gamma_1+\cdots+\gamma_N
\end{eqnarray*}
two minimal expressions for $\l$. Let us assume that $\b_1$ and $\b_2$ form  block of size 2 in the first expression. We will show that the second expression must have also a block of size two. Indeed, if all the roots that appear in the second expression are mutually orthogonal, then either $\b_1$ appears also in the second expression or there exist three roots in the second expression (say $\gamma_1, ~\gamma_2$ and $\gamma_3$) such that 
$$
(\b_1,\gamma_i)=1, \quad 1\leq i\leq 3
$$
In the latter situation, we know from Lemma \ref{II-lemma8} that 
$$
\l=\b_1+(\gamma_1+\gamma_2-\b_1)+\gamma_3+\cdots+\gamma_N
$$
is a minimal expression for $\l$. But $(\gamma_1+\gamma_2-\b_1,\gamma_3)=-1$, contradicting Lemma \ref{II-lemma1}. Therefore, this situation cannot occur and $\b_1$ must also appear in the second expression. The same argument shows that $\b_2$ must also appear in the second expression and hence together with $\b_1$ forms a block of size two, contradicting our assumption. 

We have shown that the second expression must also have a block of size two. In fact we have shown a little bit more: if $\b_1$ is part of a bock of size two in the first expression then two of the roots in the second expression that are not orthogonal on $\b_1$ must form a block of size two.

Let us assume now that there are two blocks of size two (say $\{\b_1,~ \b_2\}$ and $\{\b_3,~\b_4\}$) in the first expression. We know from the above argument that,  for each $\b_i$, $1\leq i\leq 4$, there must exist a block of size two in the second expression whose constituents are not orthogonal on $\b_i$. If the second expression has a single block of size two (say $\{\gamma_1,~\gamma_2\}$) then 
$$
(\gamma_1,\b_i)=1, \quad 1\leq i\leq 4
$$
Invoking again Lemma \ref{II-lemma8} we obtain that 
$$
\l=\gamma_1+(\b_1+\b_3-\gamma_1)+\b_2+\b_4+\cdots+\b_N
$$
is a minimal decomposition of $\l$. However,  $\b_2\not\perp \gamma_1\not\perp \b_4$, contradicting the fact that $\b_2$ and $\b_4$ are orthogonal. In conclusion, there must be a second block of size two in the second expression. 
\end{proof}
\begin{Lm}\label{II-lemma16}
Assume that $R$ is not simply laced, let $\l$ be a small dominant weight and
\begin{equation}\label{II-eq10}
\l=B_1+\cdots+B_k
\end{equation}
a canonical block decomposition. Then, all the blocks  above have size one.
\end{Lm} 
\begin{proof} We prove this by induction on $k$. If $k=1$ then $\l$ cannot have co-length two, as assured by Proposition \ref{II-prop3}{(b)}.

Assume now that $k$ is at least two and that the claim holds for all small dominant weights which have strictly fewer than $k$ blocks in a canonical block decomposition. For example, the conclusion holds for $\l-B_k$ so the first $k-1$ blocks in \eqref{II-eq10} have exactly one element $B_i=\b_i$. We only need to examine $\{B_k\}$. If it has more than one element, then by Proposition \ref{II-prop6}
it must have two elements and Proposition \ref{II-prop3}{(b)} assures that $B_k$ can also be written as 
$$
B_k=\b_k+\b_{k+1}
$$
where $\b_{k}$ and $\b_{k+1}$ are orthogonal positive roots, $\b_{k}$ is short and $\b_{k+1}$ is long. In the minimal expression
$$
\l=\b_1+\cdots+\b_{k-1}+\b_{k}+\b_{k+1}
$$
the root $\b_{k+1}$ is not equivalent to any other root because it is long (we are again appealing to Proposition \ref{II-prop3}{(b)}). In other words, $\{\b_{k+1}\}$  is a block  in the above expression. Furthermore, its height is strictly smaller than the height of $B_k$ contradicting thus the fact that $B_k$ has minimal height among all $\l$-blocks. Our assumption was therefore false, so $\{B_k\}$ has also size one.
\end{proof}
\begin{Lm}\label{II-lemma18}
Assume that $R$ is not simply laced, let $\l$ be a small weight and let
\begin{eqnarray*}
\l&=& \b_1+\cdots+\b_N\\
&=&  \b^\prime_1+\cdots+\b^\prime_N
\end{eqnarray*}
two minimal expressions for $\l$, each consisting of orthogonal roots. If the first expression contains a long root then the second expression contains a long root.
\end{Lm}
\begin{proof} 
Let is assume that $N\geq 2$, $\b_1$ is a long root and all the roots appearing in the second expression are short. Since $(\l,\b_1^\vee)$=2, there must be two roots (say $\b_1^\prime$ and $\b_2^\prime$) such that
$$
(\b_1^\prime,\b_1)=(\b_2^\prime,\b_1)=1
$$
But then the norm of $\b_1-\b^\prime_1-\b^\prime_2$ equals zero, which implies that $\b_1=\b^\prime_1+\b^\prime_2$. In particular, the second expression is not minimal contradicting the hypothesis.
\end{proof}
\begin{Lm}\label{II-lemma21}
Assume that $R$ is not simply laced, let $\l$ be a small dominant weight and
\begin{equation}\label{II-eq12}
\l=\b_1+\cdots+\b_N
\end{equation}
a canonical block decomposition. The length of the roots that appear is decreasing.
\end{Lm} 
\begin{proof} If there is a short root (say $\b_i$) followed by a long root then apply Lemma \ref{II-lemma9} for $\b_i+\b_{i+1}$ which is a dominant weight of co-length two in $R_{\{\b_1,\b_1+\b_2,\dots,\b_1+\cdots+\b_{i-1}\}}$.
\end{proof}

\subsection{Uniqueness of canonical decompositions}\label{thelist}
 We are now in position to prove that small dominant weights have unique canonical block decompositions and canonical expressions. Furthermore, we show that small dominant weights are essentially parametrized by the block sizes in their canonical block decomposition. 

\begin{Thm}\label{II-thm1}
Let $\l$ be a small dominant weight. Then $\l$ has a unique canonical block decomposition. 
\end{Thm}
\begin{proof} We prove the claim by induction on the rank of $R$.

If $R$ is not simply laced then Lemma \ref{II-lemma16} assures that any canonical block decomposition of $\l$ has only blocks of size one. The first roots that appear in any two canonical block decompositions are equal since they are dominant and they have the same length as assured by Lemma \ref{II-lemma18} and Lemma \ref{II-lemma21}. Let 
$$
\l=\b_1+\cdots+\b_N
$$
be a canonical block decomposition. The weight $\l-\b_1$ is a small dominant weight in the roots system $R_{\b_1}$ (which has smaller rank than the rank of $R$) and 
$$
\l-\b_1=\b_2+\cdots+\b_N
$$
is a canonical block decomposition which is unique by the induction hypothesis. In conclusion there is a unique canonical block decomposition for $\l$.

If $R$ is simply laced then the first blocks that appear in any two canonical block decompositions are equal by Lemma \ref{II-lemma20},  Lemma \ref{II-lemma17}, and Proposition \ref{II-prop3}{(a)}. Let 
$$
\l=B_1+\cdots+B_N
$$
be a canonical block decomposition. The weight $\l-B_1$ is a small dominant weight in the root system $R_{\{B_1\}}$ (which has smaller rank than the rank of $R$) and 
$$
\l-B_1=B_2+\cdots+B_N
$$
is a canonical block decomposition which is unique by the induction hypothesis. In conclusion there is a unique canonical block decomposition for $\l$.
\end{proof}
\begin{Thm}\label{II-thm2}
Let $\l$ be a small dominant weight. Then $\l$ has a unique canonical expression. 
\end{Thm}
\begin{proof} We prove that a canonical expression is necessarily the expression giving the canonical block decomposition. 
Indeed, assume that this claim is false and let $R$ the smallest rank root system for which it fails. Let $\l$ be a small dominant weight of minimal co-length such that $\l$ has a minimal expression 
$$
\l=\gamma_1+\cdots+\gamma_N
$$
that is different from the expression 
$$
\l=\b_1+\cdots+\b_N
$$
giving the canonical block decomposition. In such a case, it is easy to verify that $N\geq 3$, that $\gamma_N$ is part of a block of size two in the first expression and that 
$\gamma_N$ must have non-zero scalar product with $\b_1$ and $\b_N$. But then, $\hght(\gamma_N)$ is strictly smaller than $\hght(\b_N)$ and therefore the root $$\b_1+\b_N-\gamma_N$$ has height strictly larger than the height of $\b_1$. Taking into account Lemma \ref{II-lemma21} we see that this is a contradiction.
\end{proof}

Now, it is easy to see that the block sizes and information about the length of roots in blocks of size one in the canonical block decomposition of a small dominant weight determine that weight uniquely. Indeed, if $R$ is simply laced then the size of the first block determines that block uniquely. Also, the weight minus its first block is a small dominant weight in the root system $R_{\{B_1\}}$ and its canonical block decomposition is the one given. If $R_{\{B_1\}}$ is not irreducible then we need to specify the irreducible components in which $B_2$ lays. 
repeating this process, we reconstruct the dominant weight in question. For $R$ not simply laced all the blocks in the canonical block decomposition have size one but they might contain a short root or a long root. If this information is provided we can reconstruct the weight exactly as before. The notation from Introduction makes all this procedure precise. 

\subsection{Proof of Theorem \ref{II-state1}} From the above discussion above it is clear that the canonical block decomposition of a small weight must be constructed following the procedure described in Introduction. We also take into account the restrictions from Section \ref{constraints}. Eliminating from the list the symbols that correspond to weights that are not small (e.g. $[1_\ell,1_\ell]$ in type $C$) we obtain the list specified in the statement. \qed

\section{Numerical invariants}\label{invariants}

\subsection{Scalar products} The analysis of  the block structure of small weights allows us to easily study the  scalar product values $(\l,\a^\vee)$ for $\a$ a root. Let us remark first that as an immediate consequence of Theorem \ref{II-state1} we know that if the root $\a$ is $\l$--relevant then $(\l,\a^\vee)\in \{2,3\}$. Motivated by this fact, for $i$ a positive integer,  we define 
$$
\A_i(\l):=\{\a\in R~|~(\l,\a^\vee)=i\}\quad\text{and}\quad \A_i^\pm(\l):=\{\a\in R^\pm~|~(\l,\a^\vee)=i\}
$$
\begin{Prop}\label{II-prop7}
Let $\l$ be small and let $\a$ be an element of $\A_3(\l)$. Then, $\a$ is $\l$-relevant. Moreover, if a minimal expression for $\l$ is fixed, then $\a$ has positive scalar product with at most two blocks in the expression and it is orthogonal on all other roots participating in the expression.
\end{Prop}
\begin{proof}
Fix a minimal expression of $\l$. Then, either $\a$ participates in this expression or,  since blocks have at most size two and consist of short roots, $\a$ has positive scalar product with at least two orthogonal roots (say $\b_1$ and $\b_2$) that appear in the expression. But then, by Proposition \ref{II-prop3}{(c)}, $\a$ is $(\b_1+\b_2)$-relevant and hence $\l$-relevant.
If there exists a third root $\b_3$ participating in the fixed minimal expression for $\l$ such that $\b_3$ is orthogonal on $\b_1,\b_2$ and with positive scalar product with $\a$ then, after replacing $\b_1+\b_2$ with $\a+(\b_1+\b_2-\a)$ we obtain that the roots $\b_3$ and $\b_1+\b_2-\a$ participate in the same minimal expression for $
\l$ and have negative scalar product, contradicting Lemma \ref{II-lemma1}.
\end{proof}
\begin{Prop}\label{II-prop8}
For $\l$  small and $i\geq 4$ the set $\A_i(\l)$ is empty. 
\end{Prop}
\begin{proof}
Proceeding as in the proof of the previous result, we find that a root in $\A_i(\l)$ must be $\l$-relevant. But then, $(\l,\a^\vee)\leq 3$.
\end{proof}
\begin{Prop}\label{II-prop9}
Let $\l$ be small and let $\a$ be an element of $\A_2(\l)$. Then, $\a$ has non-negative scalar product with any $\l$-relevant root. Moreover, $\a$ cannot have positive scalar product with more than two blocks that appear in the same minimal expression for $\l$. If $\a$ is not relevant then it has positive scalar product with exactly one block in any given minimal expression. 
\end{Prop}
\begin{proof}
Fix a minimal expression and proceed as in the proof of Proposition \ref{II-prop7}.
\end{proof}
These results allow for a strengthening of Lemma \ref{II-lemma1}.
\begin{Cor}\label{II-cor4}
For a small weight $\l$, any two $\l$-relevant roots have non-negative scalar product.
\end{Cor}
As we have just seen, all the elements of $\A_3(\l)$ are $\l$-relevant but the situation with $\A_2(\l)$ can be different. Hence, we define
$$
\A^r_2(\l):=\{\a\in \A_2(\l)~|~\a\text{ is }\l\text{-relevant}\}\quad \text{and}\quad \A^{nr}_2(\l):=\A_2(\l)\setminus \A^r_2(\l)
$$
The subsets $\A_2^{\pm,r}(\l)$ and $\A_2^{\pm,nr}(\l)$ of $\A_2^{\pm}(\l)$ are defined in the same fashion. Of course the set of $\l$-relevant roots is $\A^r_2(\l)\cup \A_3(\l)$.
\subsection{Convex hulls} With such information available it is quite easy to specify the inclusion relations between the convex hulls of orbits of small weights. 
Indeed, it is enough to observe that if $\a$ is a $\l$-relevant root then $\l-\a$ is convex linear combination of $\l$ and $s_\a(\l)$ and 
$$
\ell^*(\l-\a)=\ell^*(\l)-1
$$
Similarly, if $\a$ is in $\A^{nr}_2(\l)$ then $\l-\a$ is convex linear combination of $\l$ and $s_\a(\l)$ and 
$$
\ell^*(\l-\a)=\ell^*(\l)
$$
Proposition \ref{II-prop7} and Proposition \ref{II-prop9} also allow us to specify exactly how the block decompositions of $\l$ and $\l-\a$ are related. These observations allow us to describe concisely  the inclusion relations between convex hulls of small weights. 
\begin{Thm}\label{II-thm3} 
In the partial order relation on orbits of small weights given by the inclusion of convex hulls induces a partial order on symbols whose covers are 
$$
\begin{aligned}
B_n: &  &&[1_\ell^k,1_s]\gtrdot [1_\ell^k]\gtrdot [1_\ell^{k-1},1_s]\\
         &  &&[1_\ell]\gtrdot [1_s]\gtrdot [\emptyset] \\
C_n: & &&[1_\ell, 1^k_s]\gtrdot [1_\ell, 1_s^{k-1}]\\
          & &&[1_\ell, 1_s^k]\gtrdot [1^{k+1}_s]\gtrdot [1_s^k]\\
          & &&[1_\ell]\gtrdot [1_s]\gtrdot [\emptyset]\\
D_n: & &&[2,1^{k}] \gtrdot [2, 1^{k-1}]\\
          &  &&[2,1^{k}] \gtrdot [1^{k+2}]\gtrdot [1^{k+1}]\\
          &  &&[2]\gtrdot [1,1]\gtrdot [1]\gtrdot [\emptyset]\\
E_6:  &  &&[2^2] \gtrdot [2,1] \gtrdot [2] \gtrdot [1^2] \gtrdot [1] \gtrdot [\emptyset]\\
E_7:  & &&[2,1] \gtrdot [1^3] \gtrdot [1^2] \gtrdot [1] \gtrdot [\emptyset] \\
           & &&[2,1] \gtrdot [2] \gtrdot [1^2]\\
E_8:  & &&[2,1] \gtrdot [2] \gtrdot [1^2] \gtrdot [1] \gtrdot [\emptyset]\\
F_4:   & &&[1_\ell, 1_s] \gtrdot [1_\ell] \gtrdot [1_s] \gtrdot [\emptyset]
\end{aligned}
$$
\end{Thm}
\begin{proof}
Straightforward from the above remarks.
\end{proof}
For the root system of type $A_n$  the symbols of small dominant weights are partitions of $n+1$ and the partial order relation described in Theorem \ref{II-thm3} is the usual partial order on partitions. In fact, it is easy to see the partial order relations described in Theorem \ref{II-thm3} can be described combinatorially in the same fashion as the partial order on partitions.

\subsection{Counting minimal expressions} \label{II-minimalexp}
Let $\l$ be a small weight. For $\mathfrak e$ a minimal expression for $\l$ denote by $\E_{[\l],\mathfrak e}$ the product of  block sizes that appear in $\mathfrak e$. From Theorem \ref{II-state1} it is clear that the possible values for $\E_{[\l],\mathfrak e}$ are 1, 2, and 4. Define,
$$
\E_{[\l]}=\sum_{\mathfrak e} \E_{[\l],\mathfrak e}
$$
where the sum is over all minimal expressions of $\l$. This positive integer counts the total number of minimal expressions for $\l$ (or any other weight in its orbit) weighted by the size of the blocks that appear in each expression.
\begin{Lm}\label{II-lemma22}
Let $\l$ be a small weight. Then,
\begin{equation}\label{II-eq13}
\sum_{\a\in\A^r_2(\l)\cup\A_3(\l)} \left((\l,\a^\vee)-1\right)\E_{[\l-\a]}=\E_{[\l]}\ell^*(\l)
\end{equation}
\end{Lm}
\begin{proof}
Keep in mind that $\A^r_2(\l)\cup\A_3(\l)$ is the set of $\l$-relevant roots. Fix $\a$ a $\l$-relevant root. Then,
$$
\left((\l,\a^\vee)-1\right)\E_{[\l-\a]}=\sum_{\a\in\mathfrak e} \E_{[\l],\mathfrak e}
$$
where the sum is over minimal expressions of $\l$ in which $\a$ participates. Therefore, the left hand side of \eqref{II-eq13} counts all minimal expressions for $\l$, each with multiplicity $\ell^*(\l)$ which is exactly the right hand side.
\end{proof}
\begin{Def}
For $\l$ small define  the following quantities
$$
\Do_{[\l]}=\sum_{\a\in\A^r_2(\l)\cup\A_3(\l)} \left((\l,\a^\vee)-1\right)
\quad\text{and}\quad
\Do^\pm_{\l}=\sum_{\a\in \A^{\pm,r}_2(\l)\cup\A^\pm_3(\l)} \left((\l,\a^\vee)-1\right)
$$
The positive integer $\Do^-_\l$ will be called the total defect of $\l$.
\end{Def}
The integer $\Do_{[\l]}$ is a weighted count of all the $\l$-relevant roots. As such, it depends only on the orbit of $\l$. The total defect of $\l$ is a weighted count of all the negative  $\l$-relevant roots. Of course,
$$
\Do_{[\l]}=\Do^+_{\l}+\Do^-_{\l}
$$
\begin{Def}
 A small dominant weight is said to be normal if all the blocks in its canonical block decomposition have size one. A small weight is said to be normal if the dominant weight in its orbit is normal.
 \end{Def}
 The terminology is motivated by the fact that the roots that appear in the canonical decomposition of such a weight are mutually orthogonal. However, they are distinguished by a more important property: if $R$ is fixed then $\E_{[\l]}$ for such weights depends only on the co-length of $\l$. One of the reasons the type $A_n$ is so well-behaved  is that all small weights have this property; it is therefore natural to say that all small weights in type $A_n$ are normal.

In our situation, it is easy to verify directly that for normal weights of co-length two  $\E_{[\l]}$ depends only of the root system. We use $\E_{[1^2]}$ to refer to this integer even for non-simply laced root systems. The only exception is in type $D$ where there are two normal weights of co-length two: $[1^2]$ and $[1,1]$. In this case we stick to the usual notation $\E_{[1^2]}$ and $\E_{[1,1]}$. If $\l$ has co-length one then $\E_{[\l]}$ equals one regardless of the root system.

From Lemma \ref{II-lemma22} it is clear that for normal weights of co-length two  $\Do_{[\l]}$ depends only of the root system, except for the situation noted above. We adopt the same convention and use $\Do_{[1^2]}$ to refer to this integer even for non-simply laced root systems and $\Do_{[1,1]}$ for the exception.

For non-simply laced root systems all small weights are normal. For simply laced root systems a small weight is small if and only if it has one expression consisting of mutually orthogonal roots as assured by Lemma \ref{II-lemma20}.

\begin{Prop}\label{II-prop10}
Let $\l$ be a normal weight of co-length at least three. Then, $\E_{[\l]}$ depends only on $R$ and $\ell^*({\l})$.
\end{Prop}
\begin{proof} Let $\l$ be a normal dominant weight of co-length $N\geq 3$ and 
$$
\b_1+\b_2+\cdots+\b_N
$$
the canonical block decomposition of $\l$. For each $2\leq i\leq N$ denote by $S_i$ the set of $(\b_1+\b_i)$-relevant roots different from $\b_1$ and $\b_i$. From Proposition \ref{II-prop7} and Proposition \ref{II-prop9} we know that $\a$ an element of $S_i$ is orthogonal on $\b_j$ for all $j\neq 1, i$ and hence $\l-\a$ is a normal weight of co-length $N-1$. It is important to note that $(\a,\b_1^\vee)=1$ for any element of $S_i$.

Let $S$ denote the set of minimal expressions that do not contain $\b_1$. 

Now, we claim that 
\begin{equation}\label{II-eq14}
\sum_{i=2}^N\sum_{\a\in S_i} \left((\l,\a^\vee)-1\right) \E_{[\l-\a]}=2\sum_{\mathfrak e\in S} \E_{[\l],\mathfrak e}
\end{equation}
For $\a$ in $S_i$ and $\tilde{\mathfrak e}$ a minimal expression of $\l-\a$
$$
 \left((\l,\a^\vee)-1\right) \E_{[\l-\a], \tilde{\mathfrak e}}= \E_{[\l],\mathfrak e}
$$
where $\mathfrak e$ is the minimal expression of $\l$ which is obtained by adjoining $\a$ to $\tilde{\mathfrak e}$. We need to investigate what kind of expressions $\mathfrak e$ appear in this fashion. Our claim \eqref{II-eq14} is that the expressions $\mathfrak e$ that appear are exactly those in $S$ and each appears precisely twice.

Let $$\gamma_1+\cdots+\gamma_N $$ be an expression in $S$. Since $(\l,\b_1^\vee)=2$ we deduce that there are exactly two roots in this expression (say $\gamma_1$ and $\gamma_2$) such that
$$
(\gamma_1,\b_1^\vee)=(\gamma_2,\b_1^\vee)=1
$$ 
Therefore, $\gamma_1$ and $\gamma_2$ are the only roots that can potentially lead to $\mathfrak e$ as described in the previous paragraph. 
Of course, $(\l,\gamma_1^\vee)\geq 2$. If there is a $\b_i$ such that $(\b_i,\gamma_1)>0$ then $\gamma_1$ is in $S_i$. The only other possibility would be $(\b_1,\gamma_1^\vee)=2$ and $(\b_i,\gamma_1)=0$ for all $2\leq i\leq N$. In this case, $\gamma_1$ is a short $\l$-relevant root and $\b_1$ is the long dominant root. Hence $\b_N$ is also short and $\gamma_1$ is in $S_N$. The same story is valid for $\gamma_2$. In any case, $\mathfrak e$ appears as described above and so exactly twice: by adjoining $\gamma_1$ to $\l-\gamma_1$ and by adjoining $\gamma_2$ to $\l-\gamma_2$. The claim \eqref{II-eq14} is therefore proved. 

Assume now that $\E_{[\l-\a]}$ depends only on co-length and denote the common value by $\E_{[1^{N-1}]}$. The equality \eqref{II-eq14} now reads
\begin{equation}\label{II-eq15}
\left(\frac{\Do_{[1^2]}}{2}-1\right)(N-1)=\frac{\E_{[\l]}}{\E_{[1^{N-1}]}}-1
\end{equation}
which shows that indeed $\E_{[\l]}$ depends only on $N$. To conclude, the case $N=3$ of our statement follows from the verification for $N=2$ and the fact that for the root system of type $D$ and $[\l]=[1^3]$, all the roots $\a$ that appear in \eqref{II-eq14} satisfy $[\l-\a]=[1^2]$. We can therefore proceed by induction on co-length and use \eqref{II-eq15} to validate the induction step.
\end{proof}

It is therefore natural to denote by $\E_{[1^N]}$ the common value of  $\E_{[\l]}$ for all normal weights $\l$ of co-length $N$. We also define the following integer
$$
\eta_{[1^N]}:=\frac{\E_{[1^N]}}{\E_{[1^{N-1}]}}
$$
For consistency we also use $\eta_{[1]}=\E_{[1]}=1$.

The non-normal weights behave similarly. For $N\geq 0$ denote 
$$
\eta_{[2,1^N]}:=\frac{\E_{[2,1^N]}}{\E_{[1^{N+1}]}}
$$
and also 
$$
\eta_{[2^2]}:=\frac{\E_{[2^2]}}{\E_{[2,1]}}
$$
Then, by the proceeding as in the proof of Proposition \ref{II-prop10} one can obtain the analogues of \eqref{II-eq15}
\begin{equation}\label{II-eq16}
2N(\eta_{[1^2]}-1)\E_{[1^{N+1}]}=2(\E_{[2,1^N]}-2\E_{[1^{N+1}]}) , \quad N\geq 0
\end{equation}
\begin{equation}\label{II-eq17}
2(\eta_{[1^2]}-1)\E_{[2,1]}=\E_{[2^2]}-2\E_{[2,1]}
\end{equation}
We will not insist any further on the proof of these equalities which, as already mentioned, are rather routine once the proof of Proposition \ref{II-prop10} is examined. 

\subsection{Proof of Theorem \ref{II-state2}} Straightforward from \eqref{II-eq15}, \eqref{II-eq16}, and \eqref{II-eq17}. \qed
\begin{Cor} Let $\l$ be a normal weight. Then,
$$
\Do_{[\l]}=\eta_{[\l]}\ell^*(\l)
$$
\end{Cor}
\begin{proof}
Straightforward from Theorem \ref{II-state2} and \eqref{II-eq13}.
\end{proof}

For normal weights $\eta_{[1^N]}$ is in fact equal to the total number of minimal expression for $[1^N]$ divided by the total number  a relevant root appears in all these expressions. For non-normal weights there is another quantity that plays a similar role.

We also need to consider the following quantity. For $\l$  weight with symbol $[2]$ we denote by $\eta_{nr}$ the total number of expressions of $\l$ as a sum of \emph{three} mutually orthogonal roots (obviously, non-relevant roots will have to be involved) divided by the number of times a root appears in all such expressions. Since there are only a few non-normal weights, we can compute this integer directly. We also define
$$
\tilde{\eta}:=\eta_{nr}+|\A_3(\l)|
$$
Also, for the root systems that have a weight with symbol $[2]$ denote
$$
\delta:=\eta_{[1^2]}-\eta_{[2]}
$$

Below we list the values of $\eta_{[\l]}$ (and, if applicable of $\tilde{\eta}$ and $\delta$)  for each root system. We leave out  $\eta_{[1]}=1$ and $\eta_{[2]}=2$  which  are independent of root system. For completeness, we also include $A_n$.
\begin{enumerate}
\item[$A_n$:] $\eta_{[\l]}=\ell^*(\l)$; $\delta=0$
\item[$B_n$:] $\eta_{[1^N]}=2N-1$
\item[$C_n$:] $\eta_{[1^N]}=2N-1$
\item[$D_n$:] $\eta_{[1,1]}=n-1$; $\eta_{[1^N]}=2N-1$; $\eta_{[2, 1^N]}=2N+2$;   $\tilde{\eta}=n-1$; $\delta=1$
\item[$E_6$:] $\eta_{[1^2]}=4$; $\eta_{[2, 1]}=5$; $\eta_{[2^2]}=8$; $\tilde{\eta}=5$; $\delta=2$
\item[$E_7$:] $\eta_{[1^2]}=5$; $\eta_{[1^3]}=9$; $\eta_{[2, 1]}=6$; $\tilde{\eta}=7$; $\delta=3$
\item[$E_8$:] $\eta_{[1^2]}=7$; $\eta_{[2, 1]}=8$;  $\tilde{\eta}=11$; $\delta=5$
\item[$F_4$:] $\eta_{[1^2]}=5$
\end{enumerate}
It is interesting to note that all the integers corresponding to normal weights are classical exponents for the root system in question. As we will see in what follows, Fourier coefficients of small dominant weights can be expressed only in terms of  the height of the weight in question and the above integers. 

\section{Fourier coefficients} \label{II-fourier}
We are now in position to describe the Fourier coefficients of small weights. The normal weights have especially nice formulas, entirely similar to the formulas in type $A_n$.  As already mentioned, all small weight in type $A_n$ have the property that $\E_{[\l]}$ depends only on $\ell^*(\l)$ (in fact  $\E_{[\l]}=\ell^*(\l)!$) and this is the property that distinguishes the normal weights from the rest for the root systems considered in this paper. The formulas for non-normal weights are similar but have one or two correction terms. 

\subsection{Notation}\label{II-vectorconventions}

First, let us briefly recall a very important convention introduced in \cite[\S 3.1]{ion-geI}. As a general convention, if $S$ is a subset the integers we use the notation
\begin{equation}\label{II-convention}
(1-t^S)
\end{equation}
to refer to $1$ if $S$ is the empty set and to $\prod_{s\in S}(1-t^{\min\{0,s\}})$ otherwise. The product is zero unless $S$ consists of negative integers. We will use the analogue notation for $(t^S-1)$. If ${\bf v}=(v_1,\dots,v_k)$ is a vector in $\Re^k$ and its coordinates in the usual standard basis are integers then we use $$
(1-t^{\bf v})
$$
to refer to $(1-t^S)$ with $S=\{ v_1,\cdots ,v_k\}$. 

The zero vector in $\Re^k$ will be denoted by $\bf 0$. Note that we suppressed any reference to $k$ from the notation this information being hopefully unambiguous from the context.
If $\bf v$ and $\bf w$ are two vectors in $\Re^k$ we write $$\bf v< w $$ if and only if $v_i<w_i$ for all $1\leq i\leq k$. For any fixed $1\leq i\leq k$ we denote by $\hat{\bf v}^{i}$ the vector in $\Re^{k-1}$ obtained by omitting the $i$-th coordinate from $\bf v$. We use the same notation in the case we need to omit more than one coordinate.

\subsection{Symmetry groups}\label{II-symgroups} For the description of Fourier coefficients of non-normal weights we need some more notation pertaining to finite group actions.
Let $K$ be a finite group acting on the set $Y$. The group consisting of bijective functions from $Y$ to itself, with the multiplication given by composition will be denoted by $S(Y)$.The orbit of the element $y$ is denoted by $K\cdot y$.

Let $\ddotbvs$ be an involution of $Y$. We write $\ddotbvs\cdot y$ for $\ddotbvs(y)$. Define the $\ddotbvs$-extended orbit of $y$ as
$$
K\langle\ddotbvs\rangle K\cdot y:=(K\cdot y)\cup (K\ddotbvs K\cdot y)
$$
Note that generally the subgroup of $S(Y)$ generated by the action of $K$ and $\ddotbvs$ is  not finite and $K\cup K\ddotbvs K$ might  not be a group. 

Assume that $\mathcal{K}=\{k_i\}_{i=1}^u$ is a set of generators for $K$ consisting of order two elements. For any $1\leq i\leq u$ let $C_i$ a (potentially empty) finite set of real-valued functions on $Y$ such that each $k_i$ stabilizes the set $C_i$ (under the induced left action on functions). We refer to the sets $C_i$ as local constraints and we say that $y\in Y$ satisfies the local constraints $C_i$ if all the functions in $C_i$ vanish at $y$. The hypothesis assures that $y$ satisfies $C_i$ if and only if $k_i\cdot y$ satisfies $C_i$.

\begin{Def} With the notation established above, we say that $z\in Y$ is an element of the constrained orbit of $y\in Y$ with respect to the action of $K$ and relative to the set of generators $\mathcal{K}=\{k_i\}_{i=1}^u$ and the local constraints $\mathcal{C}=\{C_i\}_{i=1}^u$ if there exists an element $g\in K$ such that $z=g\cdot y$, $g=k_{i_s}\cdots k_{i_1}$ and  $k_{i_r}\cdots k_{i_1}\cdot y$ satisfies $C_{i_{r+1}}$
for all $0\leq r\leq s-1$. We denote the constrained orbit of $y$ by $K\diamond y$.
\end{Def}
The notion of constrained orbit can be carried on to $\ddotbvs$-extended orbits by taking into account an extra local constraint for $\ddotbvs$.

We will be particularly interested in the situations described below.

\subsubsection{A symmetric group action}\label{II-Sn}
Let $N\geq 2$ be an integer and denote by $S_{N+1}$ the symmetric group on $N+1$ letters. The group $S_{N+1}$ is generated by the $N$ elementary transpositions $\boldsigma_j=(j,j+1)$, $1\leq j\leq N$.

Define
$$
Y(N)=\{({\bf v},\tilde{v},\partial)~|~ {\bf v}\in \Z^N, \tilde{v},\partial\in \Z \}
$$
For $1\leq j\leq N$ and 
$
y=({\bf v},\tilde{v},\partial)\in Y(N)
$ 
let $$
y^{\boldsigma_j}:=({\bf v}^{\sigma_j},\tilde{v}^{\sigma_j},\partial^{\sigma_j})
$$
 the 3-tuple defined by 
\begin{subequations} \label{II-firstequiv}
 \begin{alignat}{2}\label{II-eq44} 
 & \widehat{{\bf v}^{\sigma_j}}^j:=\widehat{\bf v}^j,&\quad & v^{\sigma_j}(j) :=\tilde{v} \\ \label{II-eq45} 
 & \partial^{\sigma_j}:=\tilde{v}+\partial-v(j),&\quad & \tilde{v}^{\sigma_j} :=v(j)
 \end{alignat}
 \end{subequations}
Our claim is that the formula \eqref{II-firstequiv} defines an $S_{N+1}$ action on $Y(N)$. Indeed, it is easy to check the relations satisfied by  the elementary transpositions
\begin{alignat*}{2}
  &y^{\boldsigma_j\boldsigma_{j}}=y&  \quad& \text{ for } 1\leq j\leq N\\ 
  &y^{\boldsigma_j\boldsigma_{k}\boldsigma_j}=y^{\boldsigma_{k}\boldsigma_{j}\boldsigma_{k}} &\quad & \text{ if } |j-k|=1\\ 
 &y^{\boldsigma_j\boldsigma_{k}}=y^{\boldsigma_{k}\boldsigma_{j}}  &  \quad&\text{ if } |j-k|\geq 2
 \end{alignat*}
Therefore, there is a group action 
$$
S_{N+1}\times Y(N)\to Y(N),\quad (w,y)\mapsto w\cdot y
$$
such that $\boldsigma_j\cdot y=y^{\boldsigma_j}$ for all $1\leq j\leq N$.

We now consider the following sets 
$$
\mathcal{Y}(N)=\{({\bf v},\tilde{v},\partial,\bar{\bf v})~|~ {\bf v}\in \Z^N, \tilde{v},\partial\in \Z, \bar{\bf v}\in \Z^{N-1} \}
$$
$$
\mathcal{Z}(N)=\{({\bf v},\tilde{v},\partial,\bar{\bf v})~|~ {\bf v}\in \Z_{\leq 0}^N, \tilde{v},\partial\in \Z_{\leq 0}, \bar{\bf v}\in \Z_{\leq 0}^{N-1} \}
$$
The symmetric group  $S_{N+1}$ acts on $\mathcal{Y}(N)$ by acting trivially on the last entry of a 4-tuple and by \eqref{II-firstequiv} on the first three entries. We say that 
$\boldaleph=({\bf v},\tilde{v},\partial,\bar{\bf v})\in \mathcal{Y}(N)$ satisfies the  the local constraint for the generator $\boldsigma_j$ if
\begin{equation} \label{II-sigmaconditions}
  \widehat{\bf v}^j=\bar{\bf v}
 \end{equation}
 The reason for which such a local constraint is natural is the following. Let us suspend  for a moment from the convention \eqref{II-convention}  the part that replaces exponents of $t$ by zero if they are strictly positive and for $\boldaleph$ an element of $\mathcal{Y}(N)$ as above let us use the notation 
\begin{equation}\label{II-aleph4}
(1-t^{\boldaleph}):= (1-t^{\bf v})-t^{\tilde{v}}(1-t^{\partial})(1-t^{\bar{\bf v}})
\end{equation}
If $\boldaleph$ satisfies the constraint \eqref{II-sigmaconditions} then $(1-t^{\boldaleph})=(1-t^{\boldsigma_j\cdot\boldaleph})$ as it can be easily checked. Therefore, $(1-t^{\boldaleph})$ does not change if we replace $\boldaleph$ by an element of its constrained orbit $S_{N+1}\diamond\boldaleph$. If we restore the convention \eqref{II-convention} and we keep the notation \eqref{II-aleph4} then $(1-t^{\boldaleph})$ does not change if we replace an element $\boldaleph\in\mathcal{Z}(N)$ by an element of $(S_{N+1}\diamond\boldaleph)\cap\mathcal{Z}(N)$.

\subsubsection{A dihedral group action}\label{II-G2} 
We use the symbol $\pcoef{\fcoef{a}{b}{c}}{\scoef{d}{e}{f}{g}{h}}$ to refer to the fact that all the entries  are integers and that 
\begin{equation}\label{II-eq57}
a+b=c,\quad d+e=g,\quad e+f=h
\end{equation}
Denote by $\dotY$ the set of all such symbols. 

Remark that in any equilateral triangle having entries as vertices, it is sufficient to specify only two entries. We take advantage of this fact by marking with $\bullet$  any entry that can be deduced using \eqref{II-eq57}. 

We define two involutions in $S(\dotY)$ by
$$
\dotbvs_1\cdot \pcoef{\fcoef{a}{b}{\bullet}}{\scoef{\bullet}{e}{\bullet}{g}{h}}=\pcoef{\fcoef{e}{b}{\bullet}}{\scoef{\bullet}{a}{\bullet}{h}{g}}
$$ 
$$
\dotbvs_2\cdot \pcoef{\fcoef{a}{b}{\bullet}}{\scoef{\bullet}{e}{\bullet}{g}{h}}=\pcoef{\fcoef{a}{e}{\bullet}}{\scoef{\bullet}{b}{\bullet}{g}{h}}
$$
It is straightforward to verify that $\dotbvs_1$ and $\dotbvs_2$ have indeed order two and that $\dotbvs_1\dotbvs_2$ has order six. Therefore the group generated by them is $D_{12}$, the dihedral group of order 12 and there is a group action 
$$
D_{12}\times \dotY\to \dotY,\quad (w,y)\mapsto w\cdot y
$$
that extends the above action of $\dotbvs_1$ and $\dotbvs_2$.

It is also easy to verify that the group $D_6$ generated by $(\dotbvs_1\dotbvs_2)^2$ and $\dotbvs_2$ is the dihedral group of order six  and it is acting on symbols $\pcoef{\fcoef{a}{b}{\bullet}}{\scoef{\bullet}{e}{\bullet}{g}{h}}$ by permuting the underscored entries $a,b,e$ and fixing $g$ and $h$. One complement of $D_6$ in $D_{12}$ is $C_2$, the cyclic group of  order two generated by $(\dotbvs_1\dotbvs_2)^3$, which acts by interchanging $g$ and $h$ and fixing the underscored entries. Note that $D_{12}$ is the internal direct product of $D_6$ and $C_2$.

We now consider the following sets 
$$
\dotcalY=\{({\bf v},\tilde{v},\partial,\bar{\bf v})~|~ {\bf v}\in \Z^3, \tilde{v},\partial\in \Z, \bar{\bf v}\in \Z^{2} \}
$$
$$
\dotcalZ=\{({\bf v},\tilde{v},\partial,\bar{\bf v})~|~ {\bf v}\in \Z_{\leq 0}^3, \tilde{v},\partial\in \Z_{\leq 0}, \bar{\bf v}\in \Z_{\leq 0}^{2} \}
$$
To an element $\boldaleph=({\bf v},\tilde{v},\partial,\bar{\bf v})$ we can associate the symbol $$\pcoef{\fcoef{v(1)}{v(2)}{\bullet}}{\scoef{\partial}{\tilde{v}}{\bar{v}(1)}{\bullet}{\bullet}}$$
This associated symbol contains all the information about $\boldaleph$ except for the coordinates $v(3)$ and $\bar{v}(2)$. If we let 
 the dihedral group  $D_{12}$ act on $\boldaleph\in\dotcalY$ by acting trivially on $v(3)$ and $\bar{v}(2)$ and on the other entries by the action induced from the action on symbols via the above correspondence, we obtain a $D_{12}$ action on $\dotcalY$. We say that 
$\boldaleph=({\bf v},\tilde{v},\partial,\bar{\bf v})\in \dotcalY$ satisfies the  the local constraint for the generator $\dotbvs_1$ (and  for $\dotbvs_2$ also) if
\begin{equation} \label{II-dvsconditions}
 {v}(3)=\bar{v}(2),\quad v(1)+v(2)=\tilde{v}+\partial+\bar{v}(1)
 \end{equation}

 The reason for which such a local constraint is natural is the following. Let us suspend  for a moment from the convention \eqref{II-convention}  the part that replaces exponents by zero if they are strictly positive and for $\boldaleph$ an element of $\dotcalY$ as above let us use the notation \eqref{II-aleph4}. If $\boldaleph$ satisfies the constraint \eqref{II-dvsconditions} then $(1-t^{\boldaleph})=(1-t^{\dotbvs_1\cdot\boldaleph})$, and similarly for $\dotbvs_2$, as it can be easily checked. Therefore, $(1-t^{\boldaleph})$ does not change if we replace $\boldaleph$ by an element of its constrained orbit $D_{12}\diamond\boldaleph$. If we restore the convention \eqref{II-convention} and we keep the notation \eqref{II-aleph4} then $(1-t^{\boldaleph})$ does not change if we replace an element $\boldaleph\in\dotcalZ$ by an element of $(D_{12}\diamond\boldaleph)\cap\dotcalZ$.


\subsubsection{A $S_4\times S_4$ action}\label{II-BCD4}
We use the symbol $\ppcoef{\fcoef{a}{b}{c}}{\scoef{d}{e}{f}{g}{h}}{\scoef{i}{j}{k}{l}{m}}$ to refer to the fact that all the entries  are integers and that 
\begin{equation}\label{II-eq58}
a+b=c,\quad d+e=g,\quad e+f=h,\quad i+j=l,\quad j+k=m
\end{equation}
Denote by $\ddotY$ the set of all such symbols. 

Remark that in any equilateral triangle having entries as vertices, it is sufficient to specify only two entries. We take advantage of this fact by marking with $\bullet$  any entry that can be deduced using \eqref{II-eq58}. 

We define five involutions in $S(\ddotY)$ as follows 
$$
\ddotbvs_1\cdot \ppcoef{\fcoef{a}{b}{\bullet}}{\scoef{\bullet}{e}{\bullet}{g}{h}}{\scoef{\bullet}{j}{\bullet}{l}{m}}=\ppcoef{\fcoef{e}{b}{\bullet}}{\scoef{\bullet}{a}{\bullet}{h}{g}}{\scoef{\bullet}{j}{\bullet}{l}{m}}
$$
$$
\ddotbvs_2\cdot \ppcoef{\fcoef{a}{b}{\bullet}}{\scoef{\bullet}{e}{\bullet}{g}{h}}{\scoef{\bullet}{j}{\bullet}{l}{m}}=\ppcoef{\fcoef{a}{e}{\bullet}}{\scoef{\bullet}{b}{\bullet}{g}{h}}{\scoef{\bullet}{j}{\bullet}{l}{m}}
$$
$$
\ddotbvs_3\cdot \ppcoef{\fcoef{a}{b}{\bullet}}{\scoef{\bullet}{e}{\bullet}{g}{h}}{\scoef{\bullet}{j}{\bullet}{l}{m}}=\ppcoef{\fcoef{a}{b}{\bullet}}{\scoef{\bullet}{j}{\bullet}{l}{m}}{\scoef{\bullet}{e}{\bullet}{g}{h}}
$$
$$
\ddotbvs_4\cdot \ppcoef{\fcoef{a}{b}{\bullet}}{\scoef{\bullet}{e}{\bullet}{g}{h}}{\scoef{\bullet}{j}{\bullet}{l}{m}}=\ppcoef{\fcoef{a}{b}{\bullet}}{\scoef{\bullet}{e}{\bullet}{g}{m}}{\scoef{\bullet}{j}{\bullet}{l}{h}}
$$
$$
\ddotbvs\cdot \ppcoef{\fcoef{a}{\bullet}{c}}{\scoef{\bullet}{e}{f}{g}{\bullet}}{\scoef{\bullet}{j}{\bullet}{l}{m}}=\ppcoef{\fcoef{e}{\bullet}{c}}{\scoef{\bullet}{a}{f}{g}{\bullet}}{\scoef{\bullet}{j}{\bullet}{l}{h}}
$$
It is clear that the action of $\ddotbvs_1$ and $\ddotbvs_2$ coincides with the action of $\dotbvs_1$ and $\dotbvs_2$ on the first two segments of a symbol. The subgroup of $S(\ddotY)$ generated by $(\ddotbvs_1\ddotbvs_2)^2$, $\ddotbvs_2$, and $(\ddotbvs_2\ddotbvs_3)^2$ is isomorphic to $S_4$ and is acting by permuting the underscored entries and is keeping $g,h,l,m$ fixed. The subgroup of $S(\ddotY)$ generated by $(\ddotbvs_1\ddotbvs_2)^3$ and $(\ddotbvs_2\ddotbvs_3)^3$ is isomorphic to $D_8$ the dihedral group of order eight and is acting by permuting $g,h,l,m$ and keeping the underscored entries fixed. The subgroup of $S(\ddotY)$ generated by this copy of $D_8$ and by $\ddotbvs_4$ is isomorphic to $S_4$ and is also acting by permuting $g,h,l,m$ and keeping the underscored entries fixed. The group generated by $\{\ddotbvs_i\}_{i=1}^3$ is isomorphic to $S_4\times D_8$  and the group 
by $\{\ddotbvs_i\}_{i=1}^4$ is isomorphic to $S_4\times S_4$.

Therefore, we have a group action $$
(S_4\times S_{4})\times \ddotY\to \ddotY,\quad (w,y)\mapsto w\cdot y
$$
that extends the  action of $\{\ddotbvs_i\}_{i=1}^4$.

We now consider the following sets 
$$
\ddotcalY=\{({\bf v},\tilde{v}^\prime,\partial^\prime,\bar{\bf v}^\prime,\tilde{v}^{\dprime},\partial^{\dprime},\bar{\bf v}^{\dprime})~|~ {\bf v}\in \Z^3, \tilde{v}^\prime,\partial^\prime, \tilde{v}^{\dprime},\partial^{\dprime}\in \Z, \bar{\bf v}^\prime,\bar{\bf v}^{\dprime}\in \Z^{2} \}
$$
$$
\ddotcalZ=\{({\bf v},\tilde{v}^\prime,\partial^\prime,\bar{\bf v}^\prime,\tilde{v}^{\dprime},\partial^{\dprime},\bar{\bf v}^{\dprime})~|~ {\bf v}\in \Z_{\leq 0}^3, \tilde{v}^\prime,\partial^\prime, \tilde{v}^{\dprime},\partial^{\dprime}\in \Z_{\leq 0}, \bar{\bf v}^\prime,\bar{\bf v}^{\dprime}\in \Z_{\leq 0}^{2} \}
$$
For an element $$\boldaleph=({\bf v},\tilde{v}^\prime,\partial^\prime,\bar{\bf v}^\prime,\tilde{v}^{\dprime},\partial^{\dprime},\bar{\bf v}^{\dprime})$$ 
we will use the notation 
\begin{align*}
&\boldaleph^\prime:=({\bf v},\tilde{v}^\prime,\partial^\prime,\bar{\bf v}^\prime)& &{}^\prime\boldaleph:=(\tilde{v}^\dprime,\partial^\dprime,\bar{\bf v}^\dprime)\\
&\boldaleph^\dprime:=({\bf v},\tilde{v}^\dprime,\partial^\dprime,\bar{\bf v}^\dprime)& &{}^\dprime\boldaleph:=(\tilde{v}^\prime,\partial^\prime,\bar{\bf v}^\prime)
\end{align*}
To an element $\boldaleph$ as above  we associate the symbol 
$$\ppcoef{\fcoef{v(1)}{v(2)}{\bullet}}{\scoef{\partial^\prime}{\tilde{v}^\prime}{\bar{v}^\prime(1)}{\bullet}{\bullet}}{\scoef{\partial^\dprime}{\tilde{v}^\dprime}{\bar{v}^\dprime(1)}{\bullet}{\bullet}}$$
This associated symbol contains all the information about $\boldaleph$ except for the coordinates $v(3)$, $\bar{v}^\prime(2)$, and $\bar{v}^\dprime(2)$. If we let 
$S_{4}\times S_4$ and respectively  $\ddotbvs$ act on $\boldaleph\in\ddotcalY$ by acting trivially on  $v(3)$, $\bar{v}^\prime(2)$, and $\bar{v}^\dprime(2)$ and on the other entries by the action induced from the action on symbols via the above correspondence, we obtain a $S_4\times S_4$ action and respectively a $\ddotbvs$ action on $\ddotcalY$. We say that 
$\boldaleph\in \ddotcalY$ satisfies the  the local constraint for the generator $\ddotbvs_1$ (and  for $\ddotbvs_2$ also) if
\begin{equation} \label{II-dvsconditionsbis}
 {v}(3)=\bar{v}^\prime(2),\quad v(1)+v(2)=\tilde{v}^\prime+\partial^\prime+\bar{v}^\prime(1)
 \end{equation}
 There is no local constraint for the action of $\ddotbvs_3$. We say that 
$\boldaleph\in \ddotcalY$ satisfies the  the local constraint for the generator $\ddotbvs_4$ if
\begin{equation} \label{II-dvsconditions4}
\partial^\prime=\partial^\dprime,\quad    \bar{v}^\prime(2)=\bar{v}^\dprime(2)
 \end{equation}
Finally, we say that $\boldaleph\in \ddotcalY$ satisfies the  the local constraint for the generator $\ddotbvs$ if
\begin{equation} \label{II-ddvsconditions}
 {v}(3)=\bar{v}^{\prime}(2),\quad v(1)=\tilde{v}^\prime+\bar{v}^\prime(1)
 \end{equation}
 
As before, let us see what these local constraints allow us to do.  Denote
\begin{equation}\label{II-eq43}
(1-t^{\boldaleph}):=(1-t^{\bf v})-t^{\tilde{v}^\prime}(1-t^{\partial^\prime})(1-t^{\bar{\bf v}^\prime})-t^{\tilde{v}^\dprime}(1-t^{\partial^\dprime})(1-t^{\bar{\bf v}^\dprime})
\end{equation}
Suspend from the convention \eqref{II-convention}  the part that replaces exponents by zero if they are strictly positive. If $\boldaleph$ satisfies the constraint for one $\ddotbvs_i$,  or for $\ddotbvs$, then $(1-t^{\boldaleph})$ equals $(1-t^{\ddotbvs_i\cdot\boldaleph})$ or, respectively  $(1-t^{\ddotbvs\cdot\boldaleph})$, as it can be directly verified. 

Therefore, $(1-t^{\boldaleph})$ does not change if we replace $\boldaleph$ by an element of its $\ddotbvs$-extended constrained orbit $(S_{4}\times S_4)\langle\ddotbvs\rangle(S_{4}\times S_4)\diamond\boldaleph$. If we restore the convention \eqref{II-convention} and we keep the notation \eqref{II-aleph4} then $(1-t^{\boldaleph})$ does not change if we replace an element $\boldaleph\in\ddotcalZ$ by an element of $((S_{4}\times S_4)\langle\ddotbvs\rangle(S_{4}\times S_4)\diamond\boldaleph)\cap\ddotcalZ$.


\subsection{Normal weights}\label{II-normal}

We are now ready to indicate the general shape of  the formula for Fourier coefficients $c_\l(t)$ for a normal weight $\l$. For any normal weight $\l$ assume there is vector $$\ag_\l=(\ag_1(\l),\dots,\ag_{\ell^*(\l)}(\l))$$ in 
$\Z^{\ell^*(\l)}$ such that the following properties are satisfied
\begin{subequations}\label{II-normaleq}
\begin{alignat}{2}
&\ag_{s_\th(\l)}\not<{\bf 0},&\quad& \text{if $\l$ is dominant}\label{II-eq18}\\
&\ag_{\l-\a_i}\not< {\bf 0} \text{ and } \ag_{s_i(\l)}\not< {\bf 0},&\quad& \text{if $\ag_\l\not< {\bf 0}$ and $(\l,\a_i^\vee)>0$} \label{II-eq19}\\
\intertext{The conditions \eqref{II-eq20}--\eqref{II-eq22} hold if $\ag_\l<{\bf 0}$.}
&\ag_\l=\ag_{s_i(\l)}=\ag_{\l-\a_i},&\quad& \text{if $\a_i\in\A_1(\l)\cup\A_2^{nr}(\l)$} \label{II-eq20}\\
\begin{split}
&\hat{\ag}^j_\l = \hat{\ag}^j_{s_i(\l)} \ =\ \ag_{\l-\a_i} \\
&\ag_{s_i(\l)}(j) = \ag_\l(j)+1 \label{II-eq21} 
\end{split}&\quad& \text{for some $j$, if $\a_i\in\A_2^r(\l)$}\\
\begin{split}
&\hat{\ag}^k_\l = \hat{\ag}^j_{s_i(\l)} \ =\ \ag_{\l-\a_i} \\
&\ag_{s_i(\l)}(j) = \ag_\l(k)+2 \label{II-eq22} 
\end{split}&\quad& \text{for some } j,k, \text{ if } \a_i\in\A_3(\l)
\end{alignat}
\end{subequations}
\begin{Thm}\label{II-thm6} Let $\l$ be a normal weight and $\ag_\l$ a vector with the properties described above. Then, 
$$
c_\l(t)=t^{\hght(\l)}(1-t^{\ag_\l})
$$
\end{Thm}
\begin{proof} From Theorem \ref{II-thm3} we know that if  $\Gamma$ is the union of the convex hulls of $W$-orbits of all normal weights for $R$ then $\Gamma\cap Q$ consists only of normal weights. Therefore, from \cite[Corollary 1.3]{ion-geI}, it is enough to check that the proposed formula for $c_\l(t)$ satisfies the system ${\rm Sys}(\Gamma)$. 

For the equation \cite[(1.8b)]{ion-geI} we have to check that if $\l$ is dominant then
$$
c_{s_\th(\l)}(t) = 0
$$
which is equivalent to \eqref{II-eq18} above. We will now verify \cite[(1.8a)]{ion-geI}
$$
c_{s_i(\l)}(t)-t^{-1}c_\l(t)=- c_{\l-\a_i}(t)+
t^{-1}c_{s_i(\l)+\a_i}(t) \ \ \   \text{if}\  \  (\l,\a_i^\vee)>0 
$$
If $\ag_{\l-\a_i}\not<{\bf 0}$ then \eqref{II-eq19}--\eqref{II-eq22}  imply that all the terms that appear are zero so the desired equality is trivially satisfied. 

If  $\ag_{\l-\a_i}<{\bf 0}$ we need to verify that 
$$
(1-t^{\ag_{s_i(\l)}})-(1-t^{\ag_{s_i(\l-\a_i)}})=t^{(\l,\a_i^\vee)-1}\left((1-t^{\ag_{\l}})-(1-t^{\ag_{\l-\a_i}})\right)
$$
When $\a_i\in\A_1(\l)\cup\A_2^{nr}(\l)$  then \eqref{II-eq20} assures that both the left-hand side and the right-hand side are zero. If $\a_i\in\A_2^r(\l)\cup\A_3(\l)$ then \eqref{II-eq21} and\eqref{II-eq22} imply the desired equality. By Proposition \ref{II-prop8} here are no other possibilities. 
\end{proof}
A  vector $\ag_\l$ with the required properties can be constructed from information extracted from the combinatorics of minimal expressions of $\l$.
This vector, which we refer to as the aggregate vector of $\l$, is a difference of two vectors
\begin{equation}\label{II-eq37}
\ag_\l:=\de_\l-\boldeta_{[\l]}
\end{equation}
The first vector
$$\de_\l=(\de_1(\l),\dots,\de_{\ell^*(\l)}(\l))\in\Z^{\ell^*(\l)}_{\geq 0}$$  will be referred to as the defect vector of $\l$. The terminology is motivated by the fact that  the  sum of coordinates equals the total defect 
\begin{equation}\label{II-defectsum}
\sum_{i=1}^{\ell^*(\l)}\de_i(\l)=\Do^-_\l
\end{equation}
The second vector, which will be referred to as the cut-off vector, is defined by $$\boldeta_{[1,1]}:=(\eta_{[1,1]},\eta_{[1]})$$ and 
$$\boldeta_{[\l]}:=(\eta_{[1^{\ell^*(\l)}]},\dots,\eta_{[1]})$$  for all the other normal weights. Except for the normal weights of co-length two in type $D_n$, the cut-off vector depends only on the co-length of $\l$. 

Theorem \ref{II-thm6} should be compared with Theorem 3.5 in \cite{ion-geI} which is the corresponding result in type $A_n$. Upon inspection, it can be seen that the aggregate vector in \cite{ion-geI} also satisfies \eqref{II-eq37} and \eqref{II-defectsum}.

 As confirmed by Theorem \ref{II-thm6} the aggregate vector and, in consequence, the defect vector (or at least their coordinates) are canonical entities. In principle, the coordinates of the defect  vector encode a canonical partition of the set  of \emph{negative} $\l$-relevant roots. We have not been able to give a conceptual, concise description of how the total defect is distributed among the components of the defect vector. One description had to do with the graph constructed in the following fashion: the vertices are $\l$-relevant roots $\a$ (counted with multiplicity $(\l,\a^\vee)-1$) and two roots are connected by and edge if and only if they are orthogonal. For $\l$ normal of co-length $N$ this graph is strongly regular  of type $$(\Do_{[1^N]}, \Do_{[1^{N-1}]}, \Do_{[1^{N-2}]}, \Do_{[1^{N-1}]}\frac{N-2}{N-1})$$
These graphs seem to admit  a partition of the set of vertices into $\eta_{[1^N]}$ subsets with $N$ elements such that sub-graphs induced by each subset are all isomorphic and as disconnected as possible. The coordinates of the defect vector count the number of vertices represented by negative roots in each sub-graph.  Such considerations are of independent interest and not very illuminating for the present discussion. We hope to return to this point of view in a future publication. 

For now, we postpone  the definition of $\de_\l$ until Appendix \ref{defects} where we will explicitly specify the defect vector for all small weights.

For normal weights of co-length one or two the defect vector can be easily defined. If $\l$ is of co-length one then $\de_\l=\Do^-_{\l}$. If $\l$ is of co-length two then
$\de_\l=(\Do_\l^-, 0)$ if $\Do_\l^-\leq \eta_{[\l]}$ and $\de_\l=(\eta_{[\l]},\Do_\l^--\eta_{[\l]})$ if $\Do_\l^-\geq \eta_{[\l]}$. It is straightforward  to verify that the properties \eqref{II-normaleq} are satisfied.

\subsection{Non-normal weights}\label{II-nonnormal}
The Fourier coefficients of non-normal weights have very similar formulas, generically the only difference being a correction term which appears because $\delta>0$. There are however exceptions due to peculiar features of some weights: the top weights in exceptional root systems ($[2^2]$ in type $E_6$ and $[2,1]$ in type $E_7,~ E_8$) need two corrections term. For all the other non-normal weights the formulas are uniform but the relationship between two coefficients as in \eqref{II-normaleq}  will require some special attention for weights with symbol $[2,1]$ in type $E_6$. Hence, we will need treat   some weights separately. Below, the word generic  refers specifically to weights with symbol $[2,1^k]$, $k\geq 0$ in type $D$, or weights with symbol $[2]$ in type  $E$.

\subsubsection{The generic case}\label{II-generic} 
 
 Let us first fix some notation. For  an integer $N\geq 2$  the data describing the Fourier coefficient for a generic non-normal weight of co-length $N$ will consist of an element of the set $\mathcal{Y}(N)$ defined as in Section \ref{II-Sn}.

Let $\boldaleph=({\bf v},\tilde{v},\partial,\bar{\bf v})\in\mathcal{Y}(N)$.
For any fixed $1\leq j\leq N$ and $1\leq k\leq N-1$ we denote 
$$
\widehat{\boldaleph}^{j|k}=(\widehat{\bf v}^j,\tilde{v},\partial,\widehat{\bar{\bf v}}^k)
$$
We use a similar notation if we need to omit more than one coordinate.
 
Below,  whenever $\l$ is a normal weight, $\ag_\l$ refers to a  vector with properties \eqref{II-normaleq}. For any $\l$ a generic non-normal weight of co-length $N$, assume there exists $$\boldaleph_\l=(\ag_\l,\tilde{a}_\l,\delta_\l,\bar{\ag}_\l)\in\mathcal{Y}(N)$$  such that, eventually after replacing it with an element of $(S_{N+1}\diamond \mathcal{Y}(N))\cap\mathcal{Z}(N)$ if $\boldaleph_\l\in\mathcal{Z}(N)$,  the conditions  \eqref{II-nonnormaleq} are satisfied
\begin{subequations} \label{II-nonnormaleq}
 \begin{alignat}{2} \label{II-eq23} 
 & \text{either } \bar{\ag}_\l\not<{\bf 0}\text{ or } \delta_\l=0\text{ (or both) }&\quad& \text{if } \ag_\l\not<{\bf 0} \\ \label{II-eq24} 
 & \ag_{s_\th(\l)}\not<{\bf 0},&\quad& \text{if } \l \text{ is dominant} \\ \label{II-eq25} 
 & \ag_{\l-\a_i}\not< {\bf 0} \text{ and } \ag_{s_i(\l)}\not< {\bf 0},&\quad& \text{if  } \ag_\l\not< {\bf 0} \text{ and } (\l,\a_i^\vee)>0 \\ \label{II-eq25bar} 
 & \bar{\ag}_{\l-\a_i}\not< {\bf 0} \text{ and } \bar{\ag}_{s_i(\l)}\not< {\bf 0},&\quad& \text{if  } \bar{\ag}_\l\not< {\bf 0} \text{ and } \a_i\in\A_2^r(\l)
 \intertext{The remaining conditions   hold if $\ag_\l<{\bf 0}$.} \label{II-eq26} 
&\boldaleph_\l=\boldaleph_{s_i(\l)}&\quad& \text{if } \a_i\in\A_1(\l)\\  \label{II-eq27} 
&\boldaleph_{s_i(\l)}=\boldaleph_\l+({\bf 0},1,0,{\bf 0}),~\ag_\l=\ag_{\l-\a_i}&\quad&\text{if } \a_i\in\A_2^{nr}(\l)\\  \label{II-eq28} 
 \begin{split}
& \widehat{\boldaleph}^{j|k}_{s_i(\l)}=\widehat{\boldaleph}^{j|k}_\l=\boldaleph_{\l-\a_i}\\
&\ag_{s_i(\l)}(j) = \ag_\l(j)+1\\ 
&\bar{\ag}_{s_i(\l)}(k) = \bar{\ag}_\l(k)+1
\end{split}&\quad& \text{for some $j, k$, if $\a_i\in\A_2^r(\l)$}\\  \label{II-eq29} 
\begin{split}
&\hat{\ag}^k_\l = \hat{\ag}^j_{s_i(\l)} \ =\ \ag_{\l-\a_i} \\
&\ag_{s_i(\l)}(j) = \ag_\l(k)+2,~ \bar{\ag}_{s_i(\l)}=\bar{\ag}_\l\\
&\tilde{a}_{s_i(\l)}=\tilde{a}_\l+2,~ \delta_{s_i(\l)}=\delta_\l
\end{split}&\quad& \text{for some } j,k, \text{ if } \a_i\in\A_3(\l)
\end{alignat}
\end{subequations} 
 
 Before we move on to the proof on the next result it is necessary to make a few remarks on the nature of the weights with symbol $[2]$ in type $D_n$. 
For $\l$ a weight with symbol $[2,1^k]$, $k\geq1$ and $\a\in\A_2^r(\l)$ we have $[\l-\a]=[2,1^{k-1}]$; for  $\a\in\A_2^{nr}(\l)$ we have $[\l-\a]=[1^{k+1}]$.
For $\l$ with symbol $[2]$, there are no relevant roots in $\A_2(\l)$ but, in turn, if  $\a\in\A_2^{nr}(\l)$ then $[\l-\a]$ is either $[1^2]$ or $[1,1]$. It turns out that the weights with symbol $[1,1]$ have features similar to those with symbol $[2,1^k]$, $k\geq 0$, and it is in fact useful to think about them as being part of that series. In consequence, we set $[2,1^{-1}]:=[1,1]$
and  for weights with symbol $[2]$ we redefine 
\begin{equation}\label{II-eq31}
\A_2^r(\l):=\{\a\in \A_2(\l)|[\l-\a]=[1,1]\},~ \A_2^{nr}(\l):=\{\a\in \A_2(\l)|[\l-\a]=[1^2]\}
\end{equation}
With this notation, in order to make sense out of \eqref{II-eq28} we need to specify $\boldaleph_\l$ for $[\l]=[1,1]$. For such a weight and $\ag_\l=(\ag_\l(1),\ag_\l(2))$ a vector satisfying 
\eqref{II-normaleq} we set 
\begin{equation}\label{II-eq32}
\boldaleph_\l:=(\ag_\l(2),\ag_\l(1),\ag_\l(2),\emptyset)
\end{equation}
which clearly satisfies 
$$
(1-t^{\boldaleph_\l})=(1-t^{\ag_\l})
$$
 
\begin{Thm}\label{II-thm7} Let $\l$ be a generic non-normal weight and $\boldaleph_\l$ a 4-tuple with the properties described above. Then, 
$$
c_\l(t)=t^{\hght(\l)}(1-t^{\boldaleph_\l})
$$
\end{Thm}
\begin{proof} Let $\Gamma$ be the union of the convex hulls of $W$-orbits of all normal and generic non-normal weights. As before, it is enough to check that the proposed formula for $c_\l(t)$ satisfies the system ${\rm Sys}(\Gamma)$. 

For the equation \cite[(1.8b)]{ion-geI} we have to check that if $\l$ is dominant then
$$
c_{s_\th(\l)}(t) = 0
$$
which immediately follows from \eqref{II-eq24} above. We will now verify \cite[1.8a)]{ion-geI}
$$
c_{s_i(\l)}(t)-t^{-1}c_\l(t)=- c_{\l-\a_i}(t)+
t^{-1}c_{s_i(\l)+\a_i}(t) \ \ \   \text{if}\  \  (\l,\a_i^\vee)>0 
$$
If $\ag_{\l-\a_i}\not<{\bf 0}$ then \eqref{II-eq26}--\eqref{II-eq29} together with \eqref{II-eq23}  imply that all the terms that appear are zero so the desired equality is trivially satisfied. 
Therefore, in what follows we can safely assume that $\ag_{\l-\a_i}\not<{\bf 0}$, which by \eqref{II-eq25} implies also $\ag_\l<{\bf 0}$.

If   $\a_i\in\A_1(\l)$ we need to verify that $$c_{s_i(\l)}(t)=t^{-1}c_\l(t)$$ which is an immediate consequence of \eqref{II-eq26}.  

If  $\a_i\in\A_2^r(\l)$ we need to verify that 
\begin{equation}\label{II-eq30}
c_{s_i(\l)}(t)-t^{-1}c_\l(t)=(t^{-1}-1)c_{\l-\a_i}(t)
\end{equation}
The set $\A_i^r(\l)$ is non-empty only in type $D_n$. Hence, we have $[\l-\a_i]=[2,1^{\ell^*(\l)-3}]$ (keep in mind our convention for $\ell^*(\l)=2$) and 
 \eqref{II-eq28} implies the desired equality.

If  $\a_i\in\A_2^{nr}(\l)$ we need to verify again \eqref{II-eq30}. 
In this case $[\l-\a_i]=[1^{\ell^*(\l)}]$  and \eqref{II-eq27} assures that the equality holds.

If   $\a_i\in\A_3(\l)$ then $[\l-\a_i]=[\l-2\a_i]=[1^{\ell^*(\l)-1}]$ and $\a_i\in\A_1(\l-\a_i)$. Therefore, taking into account \eqref{II-eq20},
we need to verify that
$$c_{s_i(\l)}(t)-t^{-1}c_\l(t)=(t^{-2}-1)c_{\l-\a_i}(t)$$
This equality follows from \eqref{II-eq29}. 
\end{proof}

Again, we are left with the task of specifying $\boldaleph_\l$ (called aggregate data) with the required properties. We will do this in Appendix \ref{defects}; for now we only say a few words regarding the relationship between this data and the combinatorics of minimal expressions. As for the case of normal weights the data takes the form
\begin{equation}\label{II-nonnormalshape}
\begin{split}
\ag_\l&:= ~\de_\l-\boldeta_{[1^{\ell^*(\l)}]}\\
\tilde{a}_{\l}&:=~ \tilde{d}_{\l}-\tilde{\eta}\\
\delta_{\l}&:=~ \partial_{\l}-\delta\\
\bar{\ag}_\l&:=~ \overline{\de}_\l-\boldeta_{[1^{\ell^*(\l)-1}]}
\end{split}
\end{equation}
The first quantities after the equal sign will be referred to as the defect data and the second quantities will be referred to as the cut-off data. The cut-off data is of course the data specified at the end of Section  \ref{II-minimalexp}.

The defect data always consists of non-negative integers; it counts in some fashion the number elements in $\A^-_2(\l)\cup\A^-_3(\l)$. For example,  the sum of coordinates of $\de_\l$  equals  $\Do^-_\l+\partial_\l$;  the sum of coordinates of $\overline{\de}_\l$  equals $|\A_2^{-,r}|$ and  $\tilde{d}_{\l}+\partial_\l$ equals  $|\A_2^{-, nr}(\l)\cup\A^{-}_3(\l)|$.

\subsubsection{Weights with symbol \texorpdfstring{$[2^{E_6},1^{A_2}]$}{[2,1]}}\label{II-E6[2,1]}

We now move on to discuss the non-generic weights. These are the weights with symbols $[2,1]$ in type $E$ and those with symbol $[2^2]$ in type $E_6$. For each, the combinatorics of minimal expressions is special in its own way. In this section we look at weights with symbol $[2,1]$ in type $E_6$. These are very close to being generic, the only difference being that   $\A_2(\l)=\A^{r}_2(\l)$ since there are no small weights with symbol $[1^3]$ in type $E_6$.

For any $\l$ with symbol $[2,1]$ assume there exists $$\boldaleph_\l=(\ag_\l,\tilde{a}_\l,\delta_\l,\bar{\ag}_\l)\in\dotcalY$$  such that, eventually after replacing it with an element of $(D_{12}\diamond \dotcalY)\cap\dotcalZ$ if $\boldaleph_\l\in\dotcalZ$,  the conditions  \eqref{II-nonnormaleq} are satisfied (keep in mind  that $\A_2^{nr}(\l)=\emptyset$). Of course, whenever a normal weight $\mu$  appears in \eqref{II-nonnormaleq}, $\ag_\mu$ refers to a  vector with properties \eqref{II-normaleq} and whenever a weight $\mu$ with symbol $[2]$ appears $\boldaleph_\mu\in\mathcal{Y}(2)$ refers to an element  satisfying properties \eqref{II-nonnormaleq}.

The following result has exactly the same proof as Theorem \ref{II-thm7}.
\begin{Thm}\label{II-thm9} Let $\l$ be a weight with symbol $[2,1]$ and $\boldaleph_\l$ a 4-tuple with the properties described above. Then, 
$$
c_\l(t)=t^{\hght(\l)}(1-t^{\boldaleph_\l})
$$
\end{Thm}

The 4-tuples $\boldaleph_\l$ with the required properties will again take the form \eqref{II-nonnormalshape}  with the note that $\eta_{[1^3]}$, which appears in the definition of $\boldeta_{[1^3]}$, is defined as 
$$
\eta_{[1^3]}:=1+2(\eta_{[1^2]}-1)=7
$$
as it would have been had a weight with symbol $[1^3]$ existed in type $E_6$.

\subsubsection{Weights with symbol \texorpdfstring{$[2^{E_7},1^{A_5}]$}{[2,1]} or  \texorpdfstring{$[2^{E_8},1^{E_6}]$}{[2,1]}}\label{II-E78[2,1]}

In this section we look at the weights with symbol $[2,1]$ in root systems of type $E_7$ and $E_8$.  Below,  whenever a normal weight $\mu$  appears in \eqref{II-nonnormaleq}, $\ag_\mu$ refers to a  vector with properties \eqref{II-normaleq} and whenever a weight $\mu$ with symbol $[2]$ appears, $\boldaleph_\mu\in\mathcal{Y}(2)$ refers to an element  satisfying properties \eqref{II-nonnormaleq}. For any $\l$ with symbol $[2,1]$ assume there exists $$\boldaleph_\l=(\ag_\l,\tilde{a}^\prime_\l,\delta^\prime_\l,\bar{\ag}^\prime_\l,\tilde{a}^\dprime_\l,\delta^\dprime_\l,\bar{\ag}^\dprime_\l)\in\ddotcalY$$  such that  eventually after replacing it with an element of $(S_4\times S_4)\langle\ddotbvs\rangle(S_4\times S_4)\diamond \ddotcalY)\cap\ddotcalZ$ if $\boldaleph_\l\in\dotcalZ$,  the conditions  \eqref{II-e78eq} are satisfied (note that $\A_2^{nr}(\l)=\emptyset$ in type $E_8$) 
\begin{subequations} \label{II-e78eq}
 \begin{alignat}{2} \label{II-eq49} 
\begin{split}
&  \text{either } \bar{\ag}^\prime_\l\not<{\bf 0}\text{ or } \delta^\prime_\l=0\text{ (or both) }\\
 & \text{either } \bar{\ag}^\dprime_\l\not<{\bf 0}\text{ or } \delta^\dprime_\l=0\text{ (or both) }\\
\end{split}
 &\quad& \text{if } \ag_\l\not<{\bf 0} \\ \label{II-eq50} 
 & \ag_{s_\th(\l)}\not<{\bf 0},&& \text{if } \l \text{ is dominant} \\ \label{II-eq51} 
 & \ag_{\l-\a_i}\not< {\bf 0} \text{ and } \ag_{s_i(\l)}\not< {\bf 0},&\quad& \text{if  } \ag_\l\not< {\bf 0} \text{ and } (\l,\a_i^\vee)>0 \\ \label{II-eq51primebar} 
& \bar{\ag}_{\l-\a_i}\not< {\bf 0} \text{ and } \bar{\ag}_{s_i(\l)}\not< {\bf 0},&\quad& \text{if  } \bar{\ag}^\prime_\l\not< {\bf 0} \text{ and } \a_i\in\A_2^r(\l) 
 \intertext{The conditions  \eqref{II-eq52}--\eqref{II-eq55} hold if $\ag_\l<{\bf 0}$.} \label{II-eq52} 
&\boldaleph_\l=\boldaleph_{s_i(\l)}&\quad& \text{if } \a_i\in\A_1(\l)\\  \label{II-eq53} 
\begin{split}
&\boldaleph_{s_i(\l)}=\boldaleph_\l+({\bf 0},1,0,{\bf 0},1,0,{\bf 0})\\
&\ag_\l=\ag_{\l-\a_i}\\
\end{split}
&\quad&\text{if } \a_i\in\A_2^{nr}(\l)\\  \label{II-eq54} 
\begin{split}
&\hat{\ag}^k_\l = \hat{\ag}^j_{s_i(\l)}  = \ag_{\l-\a_i},~ \ag_{s_i(\l)}(j) = \ag_\l(k)+2 \\
&{}^\prime\boldaleph_{s_i(\l)}={}^\prime\boldaleph_{\l}+(2,0,{\bf 0})\\
&{}^\dprime\boldaleph_{s_i(\l)}={}^\dprime\boldaleph_{\l}+(2,0,{\bf 0})
\end{split}&\quad& \text{for some } j,k, \text{ if } \a_i\in\A_3(\l)
\intertext{If $\a_i\in\A_2^r(\l)$ then for some $j, k$ then one of the following sets of equations hold}
\begin{split}
& \widehat{\boldaleph^\prime}^{j|k}_{s_i(\l)}=\widehat{\boldaleph^\prime}^{j|k}_\l=\boldaleph_{\l-\a_i}\\ \label{II-eq55} 
& {}^\prime\boldaleph_{s_i(\l)}={}^\prime\boldaleph_{\l}+(1,0,{\bf 0})\\
&\ag_{s_i(\l)}(j) = \ag_\l(j)+1\\ 
&\bar{\ag}^\prime_{s_i(\l)}(k) = \bar{\ag}^\prime_\l(k)+1
\end{split}&
\begin{split}
& \widehat{\boldaleph^\dprime}^{j|k}_{s_i(\l)}=\widehat{\boldaleph^\dprime}^{j|k}_\l=\boldaleph_{\l-\a_i}\\
& {}^\dprime\boldaleph_{s_i(\l)}={}^\dprime\boldaleph_{\l}+(1,0,{\bf 0})\\
&\ag_{s_i(\l)}(j) = \ag_\l(j)+1\\ 
&\bar{\ag}^\dprime_{s_i(\l)}(k) = \bar{\ag}^\dprime_\l(k)+1
\end{split}
 \end{alignat}
 \end{subequations} 
\begin{Thm}\label{II-thm10} Let $\l$ be a weight with symbol $[2,1]$ and $\boldaleph_\l$ a 7-tuple with the properties described above. Then, 
$$
c_\l(t)=t^{\hght(\l)}(1-t^{\boldaleph_\l})
$$
\end{Thm}
\begin{proof} Let $\Gamma$ be the union of the convex hulls of $W$-orbits of all small weights. We check that the proposed formula for $c_\l(t)$ satisfies the system ${\rm Sys}(\Gamma)$. For small weights of symbols other than $[2,1]$ the relevant equations have been checked in Theorem \ref{II-thm6} and Theorem \ref{II-thm7}.

For the equation \cite[(1.8b)]{ion-geI} we have to check that if $\l$ is dominant then
$$
c_{s_\th(\l)}(t) = 0
$$
which immediately follows from \eqref{II-eq49} and \eqref{II-eq50}. We will now verify \cite[(1.8a)]{ion-geI}
$$
c_{s_i(\l)}(t)-t^{-1}c_\l(t)=- c_{\l-\a_i}(t)+
t^{-1}c_{s_i(\l)+\a_i}(t) \ \ \   \text{if}\  \  (\l,\a_i^\vee)>0 
$$
If $\ag_{\l-\a_i}\not<{\bf 0}$ then \eqref{II-eq52}--\eqref{II-eq55} together with \eqref{II-eq49}  imply that all the terms that appear are zero so the desired equality is trivially satisfied. 
Therefore, in what follows we can safely assume that $\ag_{\l-\a_i}<{\bf 0}$, which by \eqref{II-eq51} implies also $\ag_\l<{\bf 0}$.

If   $\a_i\in\A_1(\l)$ we need to verify that $$c_{s_i(\l)}(t)=t^{-1}c_\l(t)$$ which is an immediate consequence of \eqref{II-eq52}.  

If  $\a_i\in\A_2^r(\l)$ we need to verify that 
\begin{equation}\label{II-eq56}
c_{s_i(\l)}(t)-t^{-1}c_\l(t)=(t^{-1}-1)c_{\l-\a_i}(t)
\end{equation}
We have $[\l-\a_i]=[2]$ and 
 \eqref{II-eq55} implies the desired equality.

If  $\a_i\in\A_2^{nr}(\l)$ we need to verify again \eqref{II-eq56}. 
In this case $[\l-\a_i]=[1^{3}]$ (so this may occur only in type $E_7$) and \eqref{II-eq53} assures that the equality holds.

If   $\a_i\in\A_3(\l)$ then $[\l-\a_i]=[\l-2\a_i]=[1^{2}]$ and $\a_i\in\A_1(\l-\a_i)$. Therefore, taking into account \eqref{II-eq20},
we need to verify that
$$c_{s_i(\l)}(t)-t^{-1}c_\l(t)=(t^{-2}-1)c_{\l-\a_i}(t)$$
This equality follows from \eqref{II-eq54}. 
\end{proof}

An explicit  7-tuple $\boldaleph_\l$ (called aggregate data) with the required properties will be presented in Appendix \ref{defects}. As before, it encodes the information about the combinatorics of minimal expressions. The data (modulo the constrained action of the corresponding symmetry groups) takes the form 
\begin{alignat}{6}\label{II-e78shape}
\begin{split}
& \ag_\l:= ~\de_\l-\boldeta_{[1^{3}]}-(\delta,0,0) \\ 
& \\
&
\end{split}&\quad\quad&
\begin{split}
& \tilde{a}^\prime_{\l}:=~ \tilde{d}^\prime_{\l}-{\eta}_{[1^3]} \\
&  \delta^\prime_{\l}:=~ \partial^\prime_{\l}-\delta\\
& \bar{\ag}^\prime_\l:= ~\de^\prime_\l-\boldeta_{[1^{2}]} 
\end{split}&\quad\quad&
\begin{split}
& \tilde{a}^\dprime_{\l}:=~ \tilde{d}^\dprime_{\l}-\tilde{\eta} \\
&  \delta^\dprime_{\l}:=~ \partial^\dprime_{\l}-\delta \\
& \bar{\ag}^\dprime_\l:= ~\de^\dprime_\l-\boldeta_{[1^{2}]}
\end{split}
 \end{alignat}
The first quantities after the equal sign will be referred to as the defect data and the second quantities will be referred to as the cut-off data. To keep some formulas in Appendix \ref{defects} more compact we will sometimes need to replace the cut-off data by its image under $\ddotbvs_2$ (note that the cut-off data satisfies the local constraint  for $\ddotbvs_2$). We will specify when this is necessary, the assumption being that unless specified otherwise the cut-off data is the one in \eqref{II-e78shape}.

Note that $\eta_{[1^3]}$ that appears above is defined as 
$$
\eta_{[1^3]}:=1+2(\eta_{[1^2]}-1)
$$
even if a small root with symbol $[1^3]$ does not exist (which is the case in type $E_8$). The defect data always consists of non-negative integers which count in some fashion the number of elements in $\A^-_2(\l)\cup\A^-_3(\l)$. 

\subsubsection{Weights with symbol \texorpdfstring{$[2^{E_6},2^{A_2}]$}{[2,2]}}  \label{II-E6[2,2]}

In this section $\l$ will denote a weight with symbol $[2^2]$. If $\a\in \A_2^{nr}(\l)$ then $\l-\a$ would have to be a small weight with symbol $[2,1^2]$ but there are no such weights in $E_6$. We conclude that $\A_{2}(\l)$ must be empty and the set of  roots that have positive scalar product with $\l$ is $\A_1(\l)\cup\A_3(\l)$.

It is important to note that from \eqref{II-eq16} and \eqref{II-eq17} we know that 
$
\E_{[2^2]}=160
$
so there are 40 minimal expressions for a weight with this symbol. Each relevant root appears $\E_{[2,1]}/2=10$ times in all these expressions. Hence, there are exactly 16 elements in $\A_3(\l)$ for any $\l$ with symbol $[2^2]$. The integer $\eta_{[2^2]}$ is exactly half of this number. If $\a$ is in $\A_3(\l)$ then $\l-\a$ and $\l-2\a$ have symbol $[2,1]$.

The defect data  consists of a single integer
$$
\de_\l:=|\A^-_3(\l)|
$$ 
In order to be able to write the formula more compactly we define
$$
\de_\l^\prime:=\min\{\de_\l,4\}\quad\text{and}\quad \de_\l^\dprime:=\max\{\de_\l-4,0\}
$$
Of course,
$$
\de_\l=\de_\l^\prime+\de_\l^\dprime
$$
The defect data for weights with symbol $[2,1]$ that will be defined in Appendix \ref{defects} has (up to the constrained $D_{12}$ action on $\dotcalY$) the following properties. Let $\a_i\in\A_3(\l)$. Then
\begin{subequations} \label{II-e622eq}
 \begin{alignat}{2} \label{II-eq40} 
&\de_{s_\th(\l)}\geq 8&\quad& \text{if } \l \text{ is dominant} \\ \label{II-eq39} 
& \ag_{\l-\a_i}\not<{\bf 0}  &\quad& \text{if } \de_\l\geq 8 \\ \label{II-eq41} 
\begin{split}
 & \de_{\l-\a_i}=(\de_\l,0,0),~\tilde{d}_{\l-\a_i}=0 \\  
 & \partial_{\l-\a_i}=0,~\bar{\de}_{\l-\a_i}=(\de_\l,0)
 \end{split}&\quad& \text{if } 0\leq \de_\l\leq3\\    
 \intertext{If $4\leq \de_\l\leq7$ then one of the following is true}
  \label{II-eq42}  
\begin{split}
 & \de_{\l-\a_i}=(\de_\l,0,0),~\tilde{d}_{\l-\a_i}=3 \\  
 & \partial_{\l-\a_i}=0,~\bar{\de}_{\l-\a_i}=(\de_\l-3,0)
 \end{split}&\quad\text{or}\quad& \ag_{\l-\a_i}\not<{\bf 0}
 \end{alignat}
 \end{subequations} 
As in the case of the top small weights in types $E_7$ and $E_8$ the Fourier coefficients for top weights in type $E_6$ can be described by a 7-tuple (called aggregate data)
$$\boldaleph_\l=(\ag_\l,\tilde{a}^\prime_\l,\delta^\prime_\l,\bar{\ag}^\prime_\l,\tilde{a}^\dprime_\l,\delta^\dprime_\l,\bar{\ag}^\dprime_\l)$$ 
with $\ag_\l\in \Z^4$, $\tilde{a}^\prime_\l,\tilde{a}^\dprime_\l,\delta^\prime_\l, \delta^\dprime_\l\in\Z$, and $\bar{\ag}^\prime_\l, \bar{\ag}^\dprime_\l\in\Z^{3}$. We also use the same notation as in \eqref{II-eq43}. Define, $\ag_\l, \bar{\ag}^\prime_\l, \bar{\ag}^\dprime_\l\not < {\bf 0}$ if $\ag_{\l-\a}\not <{\bf 0}$ for some $\a\in\A_3^+(\l)$. Otherwise, define
\begin{alignat}{6}\label{II-e6[2,2]shape}
\begin{split}
& \ag_\l:= ~ (\de_\l-8,\de_\l-7,-4,-1)\\ 
& \bar{\ag}^\prime_\l:= ~ (\de^\dprime_\l-4,\de^\dprime_\l-3,-1) \\
& \bar{\ag}^\dprime_\l:= ~(\de^\prime_\l-5,\de_\l^\prime-4,-1)
\end{split}&\quad\quad&
\begin{split}
&\\
& \tilde{a}^\prime_{\l}:=~  2\de^\prime_\l-10\\
&  \delta^\prime_{\l}:=~ -2
\end{split}&\quad\quad&
\begin{split}
&\\
& \tilde{a}^\dprime_{\l}:=~ 2\de^\dprime_\l -5\\
&  \delta^\dprime_{\l}:=~  -2
\end{split}
 \end{alignat}
\begin{Thm}\label{II-thm8} Let $\l$ be a weight with symbol  $[2^2]$. If the defect data for weights with symbol $[2^{E_6},1^{A_2}]$ satisfies \eqref{II-e622eq} then 
$$c_\l(t)=t^{\hght(\l)}(1-t^{\boldaleph_\l})$$
\end{Thm}
\begin{proof} We proceed in the usual fashion. 
If $\l$ is dominant then
$$
c_{s_\th(\l)}(t) = 0
$$
as assured by \eqref{II-eq40} above. 

If $\a_i\in \A_1(\l)$ then $$
\de_\l=\de_{s_i(\l)}
$$ 
so, $c_{s_i(\l)}(t)=t^{-1}c_\l(t)$ which is the desired equality. 

If $\a_i\in \A_3(\l)$ then we need to verify that $$c_{s_i(\l)}(t)-t^{-1}c_\l(t)=(t^{-2}-1)c_{\l-\a_i}(t)$$
Now, $\de_{s_i(\l)}=\de_\l+1$. If $\de_\l\geq 8$ then by taking into account \eqref{II-eq39} all the terms above are zero. If  $\de_\l\leq 7$ then \eqref{II-eq41}, \eqref{II-eq42} and Theorem \ref{II-thm7} imply the desired equality.
\end{proof}


\subsection{Proof of Theorem \ref{II-state4}} Straightforward from Theorem \ref{II-state3} and \cite[(1.3)]{ion-geI}.\qed

The Fourier coefficients of small dominant weights are particularly simple.
\begin{Cor} Let $\l$ be a small dominant weight. Then, $c_\l(t)$ equals
\begin{enumerate}
\item[i)] $(1-t^{-\boldeta_{[1,1]}})$ if $[\l]=[1,1]$
\item[ii)] $(1-t^{-\eta_{[1]}})(1-t^{-\eta_{[1^2]}})(1-t^{-\eta_{[2,1]}})(1-t^{-\eta_{[2^2]}})$ if $[\l]=[2^2]$
\item[iii)] $(1-t^{-\boldeta_{[1^{\ell^*(\l)}]}})$ if $\l$ is normal, $[\l]\neq [1,1]$
\item[iv)] $(1-t^{-\boldeta_{[1^{\ell^*(\l)}]}})-t^{-\tilde{\eta}}(1-t^{-\delta})(1-t^{-\boldeta_{[1^{\ell^*(\l)-1}]}})$ if $\l$ is non-normal, $[\l]\neq [2^2]$
\end{enumerate}
\end{Cor}

\appendix
\section{Small weights of co-length two} \label{appA}

We now embark on a through analysis of the small weights of co-length two and their minimal expressions. Our goal is to provide justification for all the statements in Proposition \ref{II-prop3}. Throughout this section $\l$ will denote a weight of co-length two.

\begin{Lm}\label{II-lemma4} Assume that $R$ is not simply laced and let $\l=\b_1+\b_2$ be a minimal expression such that $\b_1$ and $\b_2$ are non-orthogonal and have different lengths.
Then $\l$ is not small.
\end{Lm}
\begin{proof}
Assume that  $\b_1$ is short,  $\b_2$ is long, and they are not orthogonal. Then, $(\b_1,\b_2^\vee)=1$ and  
$$
\frac{1}{2}\l+\frac{1}{2}s_{\b_1-\b_2}(\l)=2\b_1
$$
Therefore, $2\b_1$ is in the convex hull of the $W$ orbit of $\l$ and in consequence in $\wt(\l)$. Hence $\l$ is not a small weight.
\end{proof}
\begin{Lm}\label{II-lemma5} Assume that $R$ is not simply laced and let $\l=\b_1+\b_2$ be a minimal expression such that $\b_1$ and $\b_2$ are non-orthogonal long roots.
Then $\l$ is not small.
\end{Lm}
\begin{proof} The hypothesis is also satisfied by all the elements in the $W$ orbit of $\l$ so we can safely assume that $\l$ is dominant. 

Remark first that  if $\l$ has another minimal expression $\l=\b_3+\b_4$ then,
$$
3=(\b_1^\vee,\l)=(\b_1^\vee,\b_3)+(\b_1^\vee,\b_4)\leq 2
$$
leading to a contradiction. 

Therefore, $\l=\b_1+\b_2$ is the unique (and hence canonical) minimal expression of $\l$. 
Then, $\b_1$ must be the dominant long root $\th_\ell$.   Since $\l$ is dominant it follows that 
$$
(\b_2,\a_k)\geq 0
$$
for all $\a_k\neq \a_{\th_\ell}$. If in addition $(\b_2, \a_{\th_\ell})\geq 0$ then $\b_2$ would also be dominant and 
since it is long it must be $\th_\ell$. Therefore, $\l=2\th_\ell$ which is not small. On the other hand, if
$$
(\b_2, \a_{\th_\ell})< 0
$$
then, 
$$
2\geq (\th_\ell^\vee,s_{\a_{\th_\ell}}(\b_2))\geq 1-(\b_2,\a_{\th_\ell}^\vee)\geq 2
$$
Consequently, $s_{\a_{\th_\ell}}(\b_2)=\b_2+\a_{\th_\ell}$ is a long root such that  $$(\th_\ell^\vee,\b_2+\a_{\th_\ell})=2
 $$
 which implies that $\th_\ell=\b_2+\a_{\th_\ell}$. 

Now, remark that  $\a_{\th_\ell}=s_{\th_\ell}(\a_{\th_\ell}-\th_\ell)$ is also a long root. Furthermore, $\a_{\th_\ell}$ is orthogonal on $\th_s$ and keeping in mind that $(\th_s,\th_\ell)=1$ it is easy to check that
$$\frac{1}{2}s_{\gamma}(\l)+\frac{1}{2}\l=2\th_s $$
with $\gamma=s_{\th_\ell-\th_s}(\a_{\th_\ell})$.
This means that $2\th_s$ is in the convex hull of the $W$ orbit of $\l$ or, equivalently,  that $\l$ is not small. 
\end{proof}

\begin{Lm}\label{II-lemma6} Assume that $R$ is not simply laced and that $\l$ is a small dominant weight. Let 
 $$
\l=\b_1+\b_2
$$
 be a minimal expression such that $\b_1$ and $\b_2$ are non-orthogonal short roots. Then,
 the above expression must necessarily be  $$\l=\th_s+(\th_s-\a_{\th_s})$$ The canonical expression of $\l$ is $$\l=\th_\ell+(2\th_s-\th_\ell-\a_{\th_s})$$ Furthermore, all minimal expressions of $\l$ consisting of orthogonal roots have one root of each length.
\end{Lm}
\begin{proof} Let us assume that the root $\b_3:=\b_1-\b_2$ is positive. Then 
$$
\l=2\b_1-\b_3
$$
and by keeping in mind that $\l$ is dominant it is easy to see that $\b_1$ must be dominant. In consequence, $\b_1=\th_s$. By repeating the argument from the proof of Lemma \ref{II-lemma5} one can show that $\b_2=\th_s-\a_{\th_s}$ and that $\a_{\th_s}$ is a short root.

Now, let us check that 
$$
\l=\th_\ell+(2\th_s-\th_\ell-\a_{\th_s})
$$
is another minimal expression. Indeed, $s_{\th_s}(\th_\ell)=\th_\ell-2\th_s$ and 
$$2\th_s-\th_\ell-\a_{\th_s}=s_{s_{\th_s}(\th_\ell)}(-\a_{\th_s})$$
 is a short root that is orthogonal on $\th_\ell$. Since this expression contains the highest root it must be the canonical expression for $\l$.

We are left to argue that if 
$$
\l=\b_3+\b_4
$$
is a minimal expression for $\l$ consisting of orthogonal roots then $\b_3$ and $\b_4$ must have different lengths. Indeed, remark that both $\b_3$ and $\b_4$ are different from $\th_s$ and 
$$
3=(\l,\th_s^\vee)=(\b_3,\th_s^\vee)+(\b_4,\th_s^\vee)
$$
so one of the scalar products must equal 2 (and the corresponding root is long) and the other must equal 1 (and the corresponding root is short).
\end{proof}

As an immediate consequence of the above considerations we obtain the following
\begin{Cor}\label{II-cor2}
Assume that $R$ is not simply laced and  that $\l$ is small and has a minimal expression consisting non-orthogonal roots. Then this minimal expression is the unique minimal expression consisting of non-orthogonal roots, other minimal expressions of $\l$ exist and they each consist of orthogonal roots of both lengths.  
\end{Cor}

\begin{Lm}\label{II-lemma3} Assume that $R$ is not simply laced and that $\l$ is a dominant weight. Fix $$\l=\b_1+\b_2$$
a minimal expression for $\l$.  If  $|\b_1|< |\b_2|$ then this expression  is not canonical.
\end{Lm}
\begin{proof}

Assume that $|\b_1|< |\b_2|$ and that the given expression is canonical. In particular, $\b_1$ is dominant so it must be $\th_s$.  

The scalar product between $\th_s$ and $\b_2$ is non-negative. Let us consider first  the case when they are orthogonal. In this case, $\b_2$ must lie in $R_{\th_s}$ and since $\l$ is dominant $\b_2$ must be a dominant root in one of the irreducible components of $R_{\th_s}$. Note that $\b_2$ and $\th_\ell$ are not orthogonal since if that is the case
 there must be an irreducible component of $R_{\th_s}$ which is also an irreducible component of $R_{\{\th_s,\th_\ell\}}$ and $\b_2$ is the long dominant root for this irreducible root system. This situation is only encountered in type $C$ and it is easy to check that $\th_s+\b_2$ is not dominant. Therefore, the scalar product $(\b_2,\th_\ell)=1$ and 
$$\th_s+\b_2-\th_\ell=s_{\b_2}s_{\th_\ell}(\th_s)$$ is a root. In consequence, 
$$\l=\th_\ell+(\th_s+\b_2-\th_\ell) $$
is the canonical expression for $\l$.
\end{proof}
\begin{Lm}\label{II-lemma7}
Assume that $R$ is  simply laced and  that $\l$ is small and has a minimal expression consisting non-orthogonal roots. Then this minimal expression is the unique minimal expression  of $\l$.   If $\l$ is dominant then $$
\l=\th+(\th-\a_{\th})
$$
\end{Lm}
\begin{proof} It is of course enough to assume that $\l$ is in addition dominant. Then, by repeating the argument from the proof of Lemma \ref{II-lemma5} one can show that
$$
\l=\th+(\th-\a_{\th})
$$
Keeping in mind that $\th$ is the highest root it is clear that $\l$ can not be written as a different sum of two roots. 
\end{proof}
\begin{Lm}\label{II-lemma8}
Let $\l=\b_1+\b_2$ be a minimal expression consisting of orthogonal roots. Then, any root $\a$ which has positive scalar product with both $\b_1$ and $\b_2$ is $\l$-relevant. Moreover, if $R$ is simply laced a $\l$-relevant root is either one of $\b_1$, $\b_2$, or has positive scalar product with both.
\end{Lm}
\begin{proof}
 If  $\a$ has positive scalar product with both $\b_1$ and $\b_2$ then $\a-\b_1$ is a root  which has positive scalar product with $\b_2$. Therefore $\b_2+\b_1-\a$  is also a root. Of course,
$$
\l=(\b_1+\b_2-\a)+\a
$$
is a minimal expression for $\l$ so $\a$ is $\l$-relevant.

Assume now that $R$ is simply laced and that $\a$ is a $\l$-relevant root different from $\b_1$ and $\b_2$. Then, there exists a root $\b_3$ such that
$$
\l=\a+\b_3
$$
The roots $\a$ and $\b_3$ are orthogonal otherwise by Lemma \ref{II-lemma7} the above expression would be the unique expression for $\l$, contradicting the hypothesis. Therefore,
$$2=(\l,\a^\vee)=(\b_1,\a^\vee)+(\b_2,\a^\vee)$$
and the conclusion follows.
\end{proof}

\begin{Lm}\label{II-lemma9} Assume that $\l$ is a small dominant weight and that $$\l=\b_1+\b_2$$ is a canonical expression consisting of orthogonal roots.
 Then,  $|\b_1|\geq |\b_2|$, $\b_1$ is dominant in $R$ and $\b_2$ is dominant in $R_{\b_1}$.
\end{Lm}
\begin{proof} The claims follow from Lemma \ref{II-lemma3} and the fact that $\l$ is dominant.
\end{proof}
\section{Defect data}\label{defects}
\subsection{Notation}  We will describe now the defect data for each irreducible root system. Thanks to our convention \eqref{II-convention}, once $\ag_\l\not<{\bf 0}$ it does not matter what is the exact formula for $\ag_\l$ and we will take advantage of this fact.

The defect data we will define in what follows does satisfy the conditions \eqref{II-normaleq},  \eqref{II-nonnormaleq},  \eqref{II-e622eq}, and  \eqref{II-e78eq}. For non-normal weights the aggregate data might be replaced by any element of its constrained orbit under the relevant symmetry group. Although some of the symmetry groups are large, in practice it turns out that the constrained orbits are generally very small. While checking the conditions mentioned above we rarely need to alter the aggregate data and when do need to alter it is always by applying less than four generators in the symmetry group. For example, in type $D_n$ the aggregate data might need to be altered by applying at most one generator of the symmetric group. Using the full $\ddotbvs$-extended orbit is also seldom necessary. In type $E_7$, for example, we only need it once when inspecting the defect data in \eqref{II-E7defect20}.

The verification itself is entirely trivial for all weights in classical root systems. It is also very easy for all the weights of exceptional root systems except for those with symbol $[2,1]$ in types $E_7$ and $E_8$. For this latter case  the verification of  \eqref{II-e78eq} is still straightforward but admittedly  a lot more tedious.

Henceforward, we will refer to  the realization of root systems in \cite[pg. 265-290]{bourbaki}. In all that follows we denote by  $\{\e_i\}_{1\leq i\leq n}$ the standard basis of  $\Re^{n}$ and by $(\cdot,\cdot)$ its canonical scalar product.

Next, we introduce some notation that will be used throughout this section. For a vector $\l=(\l_1,\cdots,\l_{n})\in\Z^n$, we denote by $\parallel \l \parallel$ the sum of the absolute values of its coordinates. 

Assume that $\l=(\l_1,\cdots,\l_{n})\in\Z^n$ has all coordinates of absolute value at most 1. Define 
$$
I(\l):=\{i~|~\l_i=-1\}
$$
As usual we write $|I(\l)|$ for the number of elements in $I(\l)$. The maximum number of elements of $I(\l)$ is $\parallel \l \parallel$.  If $I(\l)$ is not empty we  list its elements $$I(\l)=\{i_1,i_2,\dots\}$$ in increasing order. For $1\leq j\leq \parallel\l\parallel$, define 
\begin{equation}\label{II-basicdefect}
d_\l(j):= \sum_{k>i_j}|\l_k| 
\end{equation}
with the convention that $d_\l(j)$ is understood to be zero if $|I(\l)|<j$. Also, denote $$\parallel d_\l\parallel:=\sum_{j\geq 1} d_\l(j)$$

Sometimes we need similar quantities but for a different ordering of the elements in $I(\l)$. Assume $I(\l)$ is not empty and  list its elements $$I(\l)=\{j_1,j_2,\dots\}$$ in decreasing order. For $1\leq k\leq \parallel\l\parallel$, define 
\begin{equation}\label{II-basicdefectrev}
d^{rev}_\l(k):= \sum_{s<j_k}|\l_s| 
\end{equation}
with the convention that $d^{rev}_\l(k)$ is understood to be zero if $|I(\l)|<j$. Keep the same notation as above for  $\parallel d^{rev}_\l\parallel$.

Also for $\l=(\l_1,\cdots,\l_{n})\in\Z^n$ with coordinates of absolute value at most 1, let 
\begin{equation}\label{II-vecnatural}
\l^\natural:=\begin{cases}   
\l&\text{if } \parallel \l\parallel \text{ even}\\
(\l_1,\dots,\l_n,1)\in\Z^{n+1}&\text{if } \parallel \l\parallel \text{ odd}
\end{cases}
\end{equation}

Assume now that $\l=(\l_1,\cdots,\l_{n})\in\Z^n$ and its coordinates have absolute value at most 1 with the exception of exactly one coordinate whose absolute value is 2. Let us assume that $|\l_i|=2$. We denote by $\l^\sharp\in\Z^{n+1}$ and $\l^\flat\in\Z^n$ the vectors defined by
\begin{equation}\label{II-lambdaa}
\begin{split}
\l^\sharp_k&:=\l_k,\quad  \quad1\leq k<i \\
\l^\sharp_k&:=\l_i/2,\quad k=i,i+1 \\
\l^\sharp_k&:=\l_{k-1},\quad  i+2 \leq k\leq n+1
\end{split}
\end{equation}
and 
\begin{equation}\label{II-lambdab}
\begin{split}
\l^\flat_k&:=\l_k,\quad  k\neq i \hphantom{aaaaaaaaaaaa}\\
\l^\flat_i&:=0
\end{split}
\end{equation}
Note that, in particular, $\l-\l^\flat=\pm 2\e_i$. We define $\sgn(\l)$ as 
\begin{equation}\label{II-sgnD}
\begin{split}
\sgn(\l)&:=0,\quad  \text{if }\l-\l^\flat=2\e_i\\
\sgn(\l)&:=1, \quad  \text{if } \l-\l^\flat=-2\e_i
\end{split}
\end{equation}
and define $\ind(\l)$ as 
\begin{equation}\label{II-indD}
\ind(\l):=i,\quad  \text{if }\l-\l^\flat=\pm2\e_i
\end{equation}

Assume now that $\l=(\l_1,\cdots,\l_{n})\in\frac{1}{2}\Z^n$ and its coordinates have absolute value at most $1/2$ with the exception of exactly one coordinate whose absolute value is $5/2$. In this case define 
\begin{equation}\label{II-indD5/2}
\ind(\l):=i,\quad  \text{if }\l_i=\pm5/2
\end{equation}
and 
\begin{equation}\label{II-sgnD5/2}
\begin{split}
\sgn(\l)&:=0,\quad  \text{if }\l_{\ind(\l)}=5/2\\
\sgn(\l)&:=1, \quad  \text{if } \l_{\ind(\l)}=-5/2
\end{split}
\end{equation}
Furthermore, let
$$
o_\l=
\begin{cases}
|\{j~|~\ind(\l)<j,~\l_j=\frac{1}{2}\}|& \text{ if  } \sgn(\l)=0\\
|\{j~|~j<\ind(\l),~|\l_j|=\frac{1}{2}\}|+| \{j~|~\ind(\l)<j,~\l_j=\frac{1}{2}\}|& \text{ if  } \sgn(\l)=1
\end{cases}
$$
and
$$
p_\l=|\{j~|~\ind(\l)<j,~\l_j=-\frac{1}{2}\}|
$$

Finally, assume that $\l=(\l_1,\cdots,\l_{n})\in\frac{1}{2}\Z^n$ and its coordinates have absolute value $0,1/2$ or $3/2$. Define 
\begin{equation}\label{II-Jdef}
J(\l):=\{i~|~|\l_i|=\frac{3}{2}\}\quad \text{and}\quad J^{\pm}(\l):=\{i~|~\l_i=\pm\frac{3}{2}\}
\end{equation}
For $i\in J(\l)$, define 
\begin{equation}\label{II-basicdefect3/2}
J(\l,i)=
\begin{cases}
\{j~|~i<j,~\l_j=-\frac{1}{2}\}& \text{ if  } i\in J^+(\l)\\
\{j~|~j<i,~|\l_j|=\frac{1}{2}\}\cup \{j~|~i<j,~\l_j=-\frac{1}{2}\}& \text{ if } i\in J^-(\l)
\end{cases}
\end{equation}
Also, let 
$$
h_{\l}^{(N)}:=\sum_{\stackrel{S\subseteq J(\l)}{|S|=N}}|\cap_{i\in S} J(\l,i)|\quad \text{and}\quad h_\l=\sum_{i\in J^-(\l)}|\{j~|~j<i,~j\in J(\l)\}|
$$
In Section \ref{II-F4defects} we also need the following definition
\begin{equation}\label{II-basicdefect3/2rev}
J^{rev}(\l,i)=
\begin{cases}
\{j~|~j<i,~\l_j=-\frac{1}{2}\}& \text{ if  } i\in J^+(\l)\\
\{j~|~i<j,~|\l_j|=\frac{1}{2}\}\cup \{j~|~j<i,~\l_j=-\frac{1}{2}\}& \text{ if } i\in J^-(\l)
\end{cases}
\end{equation}
and 
$$
h_{\l}^{(N),rev}:=\sum_{\stackrel{S\subseteq J^{rev}(\l)}{|S|=N}}|\cap_{i\in S} J^{rev}(\l,i)|
$$

\subsection{Classical root systems} 

In this section we will work with a root system of type $B_n$, $C_n$, or $D_n$.  The Euclidean vector space $\left(\h^*_\Re, (\cdot,\cdot)\right)$ can be identified to $\left(\Re^n, (\cdot,\cdot)\right)$. Under this identification  the root lattice $Q$ is $\Z^{n}$ if $R=B_n$ and consists of the elements in $\Z^n$ whose sum of the coordinates is even if $R=C_n, D_n$. 
 
The simple roots are $$\a_i=\e_i-\e_{i+1}, \quad 1\leq i\leq n-1$$ and $\a_n$ is $\e_n$, $2\e_n$, or $\e_{n-1}+\e_n$ if $R$ is $B_n$, $C_n$, or $D_n$, respectively.
The dominant roots are $\th_\ell=\e_1+\e_2$, $\th_s=\e_1$ in type $B_n$, $\th_\ell=2\e_1$, $\th_s=\e_1+\e_2$ in type $C_n$, and $\th=\e_1+\e_{2}$ in type $D_n$.  For $R=B_n, C_n$ the Weyl group is  $S_{n}\ltimes\Z_2^n$ where the symmetric group acts by permuting the coordinates and $\Z_2^n$ by changing their signs. For $D_n$  the Weyl group is  the subgroup of $S_{n}\ltimes\Z_2^n$ which allows only for an even number of sign changes. 

 With the definitions below, \eqref{II-defectsum} holds  and  the verification of \eqref{II-normaleq} is straightforward.

\subsubsection{Type \texorpdfstring{$B_n$}{Bn}} The dominant non-zero small weights are $\e_1+\cdots+\e_k$, $1\leq k\leq n$. For $\l$ a small dominant weight, the co-length equals $$\ell^*(\l)=\frac{1}{2}\parallel \l^\natural \parallel$$ 

If $|I(\l)|>\ell^*(\l)$ then $\ag_\l\not<{\bf 0}$. Else, the defect vector $\de_\l=(\de_1(\l),\dots,\de_{\ell^*(\l)}(\l))$ is defined as 
\begin{equation}\label{II-Bdefect}
\de_\l(j):= d_{\l^\natural}(j), \quad 1\leq j\leq \ell^*(\l)
\end{equation}
The aggregate vector $\ag_\l$ is the defined by \eqref{II-eq37}. As anticipated, the components of the defect vector can be described as 
$$
\de_\l(j)=|\{-\e_{i_j}\pm\e_k~|~k>i_j\}\cap\A_2(\l)|+|\{-\e_{i_j}\}\cap\A_2(\l)|
$$
The second term is obviously 1. 
\subsubsection{Type \texorpdfstring{$C_n$}{Cn}} The dominant non-zero small weights are $\e_1+\cdots+\e_{2k}$, $1\leq k\leq n/2$ and  $2\e_1+\e_2+\cdots+\e_{2k-1}$, $1\leq k\leq (n+1)/2$.
For $\l$ a small dominant weight, the co-length equals $$\ell^*(\l)=\frac{1}{2}\parallel \l \parallel$$
Let $\l$ be in the Weyl group orbit of a weight of the form  $\e_1+\cdots+\e_{2k}$, $1\leq k\leq n/2$. If $|I(\l)|>\ell^*(\l)$ then $\ag_\l\not<{\bf 0}$. Else, the defect vector $\de_\l=(\de_1(\l),\dots,\de_{\ell^*(\l)}(\l))$ is defined as 
\begin{equation}\label{II-Cdefect1}
\de_\l(j):= d_\l(j), \quad 1\leq j\leq \ell^*(\l)
\end{equation}
Let $\l$ be  in the Weyl group orbit of a weight of the form  $2\e_1+\e_2+\cdots+\e_{2k-1}$, $1\leq k\leq (n+1)/2$. If $|I(\l^\sharp)|>\ell^*(\l)$ then $\ag_\l\not<{\bf 0}$. Otherwise,  the defect vector $\de_\l=(\de_1(\l),\dots,\de_{\ell^*(\l)}(\l))$ is defined as
\begin{equation}\label{II-Cdefect2}
\de_\l(j):= d_{\l^\sharp}(j), \quad 1\leq j\leq \ell^*(\l)
\end{equation}

The aggregate vector $\ag_\l$ is the defined by \eqref{II-eq37}. The components of the defect vector can be described in the same fashion as  for type $B_n$ with the exception that 
for the second class of weights the elements of $\A^-_3(\l)$ are also counted, with multiplicity 2. 

\subsubsection{Type \texorpdfstring{$D_n$}{Dn}}\label{II-Ddefectdata} The dominant non-zero small weights are, up to diagram automorphisms,   $2\e_1$, $2\e_1+\e_2+\cdots+\e_{2k-1}$, $2\leq k\leq (n+1)/2$, and  $\e_1+\cdots+\e_{2k}$, $1\leq k\leq n/2$.

 The weights $\pm 2\e_i$ have symbol $[1,1]$ and co-length two. For small weights of co-length two the  defect vector $\de_\l=(\de_\l(1),\de_\l(2))$ was described in the last paragraph of Section \ref{II-normal}. More explicitly,
\begin{equation}\label{II-Ddefect[1,1]}
\begin{split}
&\de_{2\e_i}(1):= i-1,\;\;\;\quad\de_{2\e_i}(2)=0\\
&\de_{-2\e_i}(1):= n-1,\quad\de_{-2\e_i}(2)=n-i
\end{split}
\end{equation}

For all the other small weights the co-length equals $$\ell^*(\l)=\frac{1}{2}\parallel \l \parallel$$
The weights  $\l$ in the Weyl group orbit of a weight of the form  $\e_1+\cdots+\e_{2k}$, $1\leq k\leq n/2$ are normal. If $|I(\l)|>\ell^*(\l)$ then $\ag_\l\not<{\bf 0}$. Otherwise, the defect vector $\de_\l=(\de_1(\l),\dots,\de_{\ell^*(\l)}(\l))$ is defined as 
\begin{equation}\label{II-Ddefect1}
\de_\l(j):= d_\l(j), \quad 1\leq j\leq \ell^*(\l)
\end{equation}
The aggregate vector $\ag_\l$ is the defined by \eqref{II-eq37}. The verification of \eqref{II-normaleq} is straightforward.

The weights $\l$ in the Weyl group orbit of a weight of the form  $2\e_1+\e_2+\cdots+\e_{2k-1}$, $1\leq k\leq (n+1)/2$ are non-normal. If $|I(\l^\sharp)|>\ell^*(\l)$ then $\ag_\l\not<{\bf 0}$. Otherwise, the defect data  is defined as
\begin{equation}\label{II-Ddefect2}
\begin{split}
\de_\l(j)&:= d_{\l^\sharp}(j), \quad 1\leq j\leq \ell^*(\l)\\
\bar{\de}_\l(j)&:= d_{\l^\flat}(j), \quad 1\leq j\leq \ell^*(\l)-1\\
\partial_\l&:=\sgn(\l)\\
\tilde{d}_\l&:=(2(n-\ind(\l))-1)\sgn(\l)+\sum_{j<\ind(\l)}(1-\l_j)
\end{split}
\end{equation}

For later use let define
\begin{equation}\label{II-Ddefect2rev}
\tilde{d}^{rev}_\l:=(2(\ind(\l)-1)-1)\sgn(\l)+\sum_{\ind(\l)<j\leq 5}(1-\l_j)
\end{equation}
The aggregate vector $\ag_\l$ is the defined by \eqref{II-nonnormalshape}. The definition of $\tilde{d}_\l$ has the following interpretation. When $\l-\l^\flat=2\e_i$, $\tilde{d}_\l+\partial_\l$ equals the number of zero coordinates of $\l$ to the left of $2\e_i$ (that is $|\A_2^{-,nr}(\l)|$) plus twice the number of negative entries to the left of $2\e_i$ (that is $2|\A^-_3(\l)|$). When $\l-\l^\flat=-2\e_i$, $\tilde{d}_\l+\partial_\l$ equals the number of zero coordinates of $\l$ plus the number of zero coordinates to the right of $-2\e_i$ (that is $|\A_2^{-,nr}(\l)|$) plus twice the number of negative entries to the left of $-2\e_i$ plus the number of non-zero coordinates to the right of $-2\e_i$(that is $2|\A^-_3(\l)|$). 

With the definition above, it is  necessary to replace $\boldaleph$ by an element in its constrained orbit only when 
$\sgn(\l)\neq\sgn(s_i(\l))$ in which case we need to apply $\boldsigma_k$ where $k$ is the smallest integer for which $\de_\l(k)=0$.

\subsection{The root system \texorpdfstring{$E_6$}{E6}}

Denote $$\teps_6:=(-\e_6-\e_7+\e_8)\in\Re^8 \quad \text{and} \quad V=\Re\e_1+\Re\e_2+\Re\e_3+\Re\e_4+\Re\e_5+\Re\teps_6\subset\Re^8$$

The Euclidean vector space $\left(\h^*_\Re, (\cdot,\cdot)\right)$ can be identified to $\left(V, (\cdot,\cdot)\right)$. The simple roots are $\a_1=\frac{1}{2}(\e_1-\e_2-\e_3-\e_4-\e_5)+\frac{1}{2}\teps_6$, $\a_2=\e_1+\e_2$, $\a_i=-\e_{i-2}+\e_{i-1}$, $3\leq i\leq 6$. The dominant root is $\th=\frac{1}{2}(\e_1+\e_2+\e_3+\e_4+\e_5)+\frac{1}{2}\teps_6$. The Weyl group orbit of an element of the root lattice cannot be described concisely in general. For this reason we will split an orbit into orbits under $W_{D_5}$, the parabolic Weyl group of type $D_5$ obtained by excluding the simple root $\a_1$. This group acts only on the first five coordinates of a vector in $\Re^8$ by permuting them and possibly changing an even number of signs. 

The positive roots are $\pm\e_i+\e_j$, $1\leq i<j\leq 5$ and the elements of the $W_{D_5}$ orbit of $\th$. The dominant non-zero small weights are $\frac{1}{2}(\e_1+\e_2+\e_3+\e_4+\e_5)+\frac{1}{2}\teps_6$ (with symbol $[1]$),  $\e_5+\teps_6$ (with symbol $[1^2]$), $\e_3+\e_4+\e_5+\teps_6$ (with symbol $[2]$), $\e_4+2\e_5+\teps_6$ and $\frac{1}{2}(-\e_1+\e_2+\e_3+\e_4+\e_5)+\frac{3}{2}\teps_6$ (with symbol $[2,1]$), $3\e_5+\teps_6$ and $2\teps_6$ (with symbol $[2^2]$). We will describe separately the defect data for each of these weights. 
One common feature is that  $\ag_\l\not<{\bf 0}$ if the coefficient of $\teps_6$ in $\l$ is strictly negative.

Some of these weights are in fact in the root lattice of the parabolic root subsystem of type $D_5$ mentioned above. For these, it is convenient to use the detect data defined in Section 
\ref{II-Ddefectdata} and it is this data we refer to when we use  notation such as $\de_\l^{D_5}$. However, remark that the indexing of the nodes of the Dynkin diagram in Section  \ref{II-Ddefectdata} differs from the one induced here on the parabolic roots subsystem of type $D_5$. To obtain the correct formulas one needs to use $d_\l^{rev}$ in the formulas \eqref{II-Ddefect1} and \eqref{II-Ddefect2} and also $\tilde{d}^{rev}_\l$ as defined by \eqref{II-Ddefect2rev}.
\subsubsection{Weights with symbol $[1]$} These are  the roots. The defect vector is defined by $\de_\a=0$ if $\a\in R^+$ and $\de_\a=1$ if $\a\in R^-$.

\subsubsection{Weights with symbol $[1^2]$}  We describe the  Weyl group orbit of $\e_5+\teps_6$  by specifying the $W_{D_5}$ orbits it contains. These are the $W_{D_5}$ orbits of 
$\e_5+\teps_6$, $2\e_5$, $\e_2+\e_3+\e_4+\e_5$, $\frac{1}{2}(-\e_1+\e_2+\e_3+\e_4+3\e_5)+\frac{1}{2}\teps_6$, $\e_5-\teps_6$, and $\frac{1}{2}(\e_1+\e_2+\e_3+\e_4+3\e_5)-\frac{1}{2}\teps_6$. For the elements of the last two $W_{D_5}$ orbits  the coordinate of $\teps_6$ is strictly negative and  we set $\ag_\l\not<{\bf 0}$. We focus our attention on the first four orbits. 

For $\l$ in $W_{D_5}(\e_5+\teps_6)$, the defect vector is defined by $\de_\l=(0,0)$. For $\l$ in $W_{D_5}(2\e_5)$, the defect vector is defined by $\de_\l=\de_\l^{D_5}$. Let $\l$ in $W_{D_5}(\e_2+\e_3+\e_4+\e_5)$. Note that $|\A_2(\e_2+\e_3+\e_4+\e_5)|=1$. If $|I(\l)|>2$ then $\ag_\l\not<{\bf 0}$. Otherwise, define $$\de_\l=\de_\l^{D_5}+(1,0)$$

For  $\l$ in $W_{D_5}(\frac{1}{2}(-\e_1+\e_2+\e_3+\e_4+3\e_5)+\frac{1}{2}\teps_6)$ the defect vector is $\de_\l=(f_{\l-\frac{1}{2}\teps_6},0)$. The defect data presented in the section is nothing else but the one defined in the last paragraph of Section \ref{II-normal}.
\subsubsection{Weights with symbol $[2]$} We describe the  Weyl group orbit of $\e_3+\e_4+\e_5+\teps_6$  by specifying the $W_{D_5}$ orbits it contains. These are the $W_{D_5}$ orbits of 
$\e_3+\e_4+\e_5+\teps_6$, $\e_3+\e_4+2\e_5$, $\frac{1}{2}(\e_1+\e_2+\e_3+3\e_4+3\e_5)+\frac{1}{2}\teps_6$, $\e_3+\e_4+\e_5-\teps_6$, and $\frac{1}{2}(-\e_1+\e_2+\e_3+3\e_4+3\e_5)-\frac{1}{2}\teps_6$. For the elements of the last two $W_{D_5}$ orbits  the coordinate of $\teps_6$ is strictly negative and  we set $\ag_\l\not<{\bf 0}$. We focus our attention on the first three orbits. 

Let $\l$ in $W_{D_5}(\e_3+\e_4+\e_5+\teps_6)$. The defect vector is defined by  
\begin{alignat}{3}\label{II-E6defect1}
\begin{split}
\de_\l&:= (0,0),\\
\tilde{d}_\l&:=\parallel d^{rev}_\l\parallel,
\end{split}
&\quad \quad&
\begin{split}
\bar{\de}_\l&:= 0\\
\partial_\l&:=0
\end{split}
\end{alignat}

Let $\l$ in $W_{D_5}(\e_3+\e_4+2\e_5)$. If $|I(\l^\sharp)|>2$ then set $\ag_\l\not<{\bf 0}$. Otherwise,  the defect data is defined by
\begin{alignat}{3}\label{II-E6defect2}
\begin{split}
\de_\l&:= \de_{\l}^{D_5}+(2,0), \\
 \tilde{d}_\l&:=\tilde{d}_\l^{D_5}+1,
\end{split}
&\quad\quad &
\begin{split}
\bar{\de}_\l&:= \bar{\de}_{\l}^{D_5}\\
\partial_\l&:=\partial^{D_5}_\l+1
\end{split}
\end{alignat}
Note that $|\A^{-,nr}_2(\e_2+\e_3+\e_4+2\e_5)|=2$ and this is reflected in the constants present in the formulas above. 

Let $\l$ in $W_{D_5}(\frac{1}{2}(\e_1+\e_2+\e_3+3\e_4+3\e_5)+\frac{1}{2}\teps_6)$. The defect data is defined by
\begin{alignat}{3}\label{II-E6defect3}
\begin{split}
\de_\l&:= (\min\{h^{(2)}_{\l-\frac{1}{2}\teps_6},\delta\}+2h_{\l},0),\\
\tilde{d}_\l&:=h^{(1)}_{\l-\frac{1}{2}\teps_6}-\min\{h^{(2)}_{\l-\frac{1}{2}\teps_6},\delta\}+2h_{\l},
\end{split}
&\quad\quad&
\begin{split}
\bar{\de}_\l&:=0\\
\partial_\l&:=\min\{h^{(2)}_{\l-\frac{1}{2}\teps_6},\delta\}
\end{split}
\end{alignat}

\subsubsection{Weights with symbol $[2,1]$} There are two dominant weights with this symbol, as described above. The first orbit is the union of the  $W_{D_5}$ orbits of $\e_4+2\e_5+\teps_6$, $\e_1+\e_2+\e_3+\e_4+2\e_5$, $\frac{1}{2}(\e_1+\e_2+\e_3+\e_4+5\e_5)+\frac{1}{2}\teps_6$, $\frac{1}{2}(\e_1+\e_2+3\e_3+3\e_4+3\e_5)-\frac{1}{2}\teps_6$, $\e_1+\e_2+\e_3+\e_4+\e_5-\teps_6$, and $\frac{1}{2}(\e_1+\e_2+\e_3+\e_4+\e_5)-\frac{3}{2}\teps_6$. The second orbit is the union of the  $W_{D_5}$ orbits of $\frac{1}{2}(-\e_1+\e_2+\e_3+\e_4+\e_5)+\frac{3}{2}\teps_6$, $-\e_1+\e_2+\e_3+\e_4+\e_5+\teps_6$, $\frac{1}{2}(-\e_1+\e_2+3\e_3+3\e_4+3\e_5)+\frac{1}{2}\teps_6$, $-\e_1+\e_2+\e_3+\e_4+2\e_5$, $\frac{1}{2}(-\e_1+\e_2+\e_3+\e_4+5\e_5)-\frac{1}{2}\teps_6$, and $\e_4+2\e_5-\teps_6$. For the elements of the $W_{D_5}$ orbits  for which the coordinate of $\teps_6$ is strictly negative  we set $\ag_\l\not<{\bf 0}$.

Let $\l$ in $W_{D_5}(e_4+2\e_5+\teps_6)$. The defect data  is defined as
\begin{equation}\label{II-E6defect5}
\begin{split}
\de_\l(1)&:= (2(\ind(\l)-1)-1)\sgn(\l)+\sum_{\ind(\l)<j\leq 5}(1-\l_j)\\
\bar{\de}_\l(1)&:=(-2 +2\sum_{j<\ind(\l)}(1-|\l_j|))\sgn(\l)+\sum_{\ind(\l)<j\leq 5}(1-|\l_j|)\\
\de_\l(2)&:=0,\quad\de_\l(3):=0, \quad \bar{\de}_\l(2):=0\\
\tilde{d}_\l&:=(1+2\sum_{j<\ind(\l)}|\l_j|)\sgn(\l)+\sum_{\ind(\l)<j\leq 5}(|\l_j|-\l_j),\quad\partial_\l:=0
\end{split}
\end{equation}

The definitions of $\de_\l(1),\bar{\de}_\l(1)$, and $ \tilde{d}_\l$ have the following interpretation. When $\sgn(\l)=0$, $\de_\l(1)+\sgn(\l)$ equals the number of zero coordinates of $\l$ to the right  of $2\e_i$ (that is $|\A_2^{-,nr}(\l)|$) plus twice the number of negative entries to the right of $2\e_i$ (that is $2|\A^-_3(\l)|$). When $\sgn(\l)=1$, $\de_\l(1)+\sgn(\l)$ equals the number of zero coordinates of $\l$ plus the number of zero coordinates to the right of $-2\e_i$ (that is $|\A_2^{-,nr}(\l)|$) plus twice the number of negative entries to the right of $-2\e_i$ plus the number of non-zero coordinates to the left of $-2\e_i$(that is $2|\A^-_3(\l)|$). Similarly, it can be seen that $\bar{\de}_\l(1)+2\sgn(\l)$ equals $|\A_2^{-,nr}(\l)|$ and $\tilde{d}_\l-\sgn(\l)$ equals $2|\A^-_3(\l)|$.

Let $\l$ in $W_{D_5}(\frac{1}{2}(\e_1+\e_2+\e_3+\e_4+5\e_5)+\frac{1}{2}\teps_6)$. The defect data  is defined as
\begin{equation}\label{II-E6defect6}
\begin{split}
\de_\l(1)&:= 2(\ind(\l)-1)\sgn(\l)+   2p_{\l-\frac{1}{2}\teps_6}+\min\{p_{\l-\frac{1}{2}\teps_6},\delta\}+o_{\l-\frac{1}{2}\teps_6}-\de_\l(2)\\  
\de_\l(2)&:=  2\left\lceil\frac{1}{2}(\ind(\l)-1)\right\rceil\sgn(\l)+   2\left\lfloor\frac{1}{2}\max\{p_{\l-\frac{1}{2}\teps_6}-\delta,0\}\right\rfloor+\min\{p_{\l-\frac{1}{2}\teps_6},\delta\}\\
\de_\l(3)&:=0\\
\bar{\de}_\l(1)&:=\min\{p_{\l-\frac{1}{2}\teps_6},\delta\}+ o_{\l-\frac{1}{2}\teps_6}, \quad \bar{\de}_\l(2):=0\\
\partial_\l&:=\min\{p_{\l-\frac{1}{2}\teps_6},\delta\}\\
\tilde{d}_\l&:=2(\ind(\l)-1)\sgn(\l)+2p_{\l-\frac{1}{2}\teps_6}-\min\{p_{\l-\frac{1}{2}\teps_6},\delta\}
\end{split}
\end{equation}
The quantities that appear above have the following meaning: $2(\ind(\l)-1)\sgn(\l)+   2p_{\l-\frac{1}{2}\teps_6}$ equals $2|\A^-_3(\l)|$ and $o_{\l-\frac{1}{2}\teps_6}$ equals $|\A^-_2(\l)|$.

Let $\l$ in $W_{D_5}(\pm\e_1+\e_2+\e_3+\e_4+2\e_5)$. If $|I(\l^\sharp)|>3$ then set $\ag_\l\not<{\bf 0}$. Otherwise,  the defect data is defined by
\begin{alignat}{3}\label{II-E6defect4}
\begin{split}
\de_\l&:= \de_{\l}^{D_5}+(2,1,0),\\
\tilde{d}_\l&:=\tilde{d}_\l^{D_5}+1,
\end{split}
&\quad\quad&
\begin{split}
\bar{\de}_\l&:= \bar{\de}_{\l}^{D_5}+(1,0)\\
\partial_\l&:=\partial^{D_5}_\l+1
\end{split}
\end{alignat}
Note that $|\A^{-,r}_2(\pm\e_1+\e_2+\e_3+\e_4+2\e_5)|=4$ and this is reflected in the constants present in the formulas above. 

For  $\l$ in $W_{D_5}(\frac{1}{2}(-\e_1+\e_2+\e_3+\e_4+\e_5)+\frac{3}{2}\teps_6)$ define the defect data to be 
\begin{alignat}{3}\label{II-E6defect11}
\begin{split}
\de_\l&:=(0,0,0),\\
\tilde{d}_\l&:=0,
\end{split}
&\quad\quad&
\begin{split}
\bar{\de}_\l&:=(0,0)\\
\partial_\l&:=0
\end{split}
\end{alignat}

For  $\l$ in $W_{D_5}(-\e_1+\e_2+\e_3+\e_4+\e_5+\teps_6)$ define the defect data to be 
\begin{equation}\label{II-E6defect7}
\begin{split}
\de_\l&:= (\parallel d^{rev}_{\l-\teps_6}\parallel-\parallel d^{rev}_{\l-\l_1\e_1-\l_2\e_2-\teps_6}\parallel,0,0)\\
\bar{\de}_\l&:= (\parallel d^{rev}_{\l-\teps_6}\parallel-2\parallel d^{rev}_{\l-\l_1\e_1-\l_2\e_2-\teps_6}\parallel,0)\\
\tilde{d}_\l&:=\parallel d^{rev}_{\l-\l_1\e_1-\l_2\e_2-\teps_6}\parallel,\quad \partial_\l:=0\\
\end{split}
\end{equation}

Let $\l$ in $W_{D_5}(\frac{1}{2}(-\e_1+\e_2+3\e_3+3\e_4+3\e_5)+\frac{1}{2}\teps_6)$. If $h_{\l}=3$, or if $h_{\l}=2$ and $h^{(1)}_{\l-\frac{1}{2}\teps_6}\geq 4$, then we set $\ag_\l\not<{\bf 0}$. Otherwise,  the defect data is defined as follows. If  $h_{\l}\leq 1$ the defect data is defined by 
\begin{subequations}\label{II-E6defect10}
\begin{alignat}{2}
\begin{split}\label{II-E6defect8}
\de_\l&:= (h^{(1)}_{\l-\frac{1}{2}\teps_6}+2h_{\l},h^{(3)}_{\l-\frac{1}{2}\teps_6},0)+(1,0,0)\\
\bar{\de}_\l&:= (h^{(1)}_{\l-\frac{1}{2}\teps_6}-h^{(2)}_{\l-\frac{1}{2}\teps_6}+h^{(3)}_{\l-\frac{1}{2}\teps_6}-\partial_\l,0)+(1,0)\\
\partial_\l&:=\max\{h^{(3)}_{\l-\frac{1}{2}\teps_6}-\left\lfloor\frac{1}{2}h^{(2)}_{\l-\frac{1}{2}\teps_6}\right\rfloor,0\}\\
\tilde{d}_\l&:=h^{(2)}_{\l-\frac{1}{2}\teps_6}+2h_{\l}
\end{split}
\intertext{If  $h_{\l}= 2$, $h^{(1)}_{\l-\frac{1}{2}\teps_6}\leq 2$, and $h^{(2)}_{\l-\frac{1}{2}\teps_6}\leq 1$ the defect data is defined by }
\begin{split}\label{II-E6defect9}
\de_\l&:= (h^{(1)}_{\l-\frac{1}{2}\teps_6}, h^{(2)}_{\l-\frac{1}{2}\teps_6},0)+(5,1,0)\\
\bar{\de}_\l&:= (h^{(1)}_{\l-\frac{1}{2}\teps_6}-h^{(3)}_{\l-\frac{1}{2}\teps_6},0)+(2,0)\\
\partial_\l&:=h^{(2)}_{\l-\frac{1}{2}\teps_6}+1\\
\tilde{d}_\l&:=h^{(3)}_{\l-\frac{1}{2}\teps_6}+3
\end{split}
\end{alignat}
\end{subequations}
Note that $|\A^{-,r}_2(\frac{1}{2}(-\e_1+\e_2+3\e_3+3\e_4+3\e_5)+\frac{1}{2}\teps_6)|=1$ and this is reflected in the constants present in the formulas above.
\subsubsection{Weights with symbol $[2^2]$} There are two dominant weights with this symbol, as described above. The first orbit is the union of the  $W_{D_5}$ orbits of $3\e_5+\teps_6$, $\frac{1}{2}(3\e_1+3\e_2+3\e_3+3\e_4+3\e_5)-\frac{1}{2}\teps_6$,  and $-{2}\teps_6$. The second orbit is the union of the  $W_{D_5}$ orbits of $-{2}\teps_6$,  $\frac{1}{2}(-3\e_1+3\e_2+3\e_3+3\e_4+3\e_5)+\frac{1}{2}\teps_6$,  and $3\e_5-\teps_6$. For the elements of the $W_{D_5}$ orbits  for which the coordinate of $\teps_6$ is strictly negative or the coordinate of $\e_5$  is $-\frac{3}{2}$, we set $\ag_\l\not<{\bf 0}$. For all the other weights we need to specify $\de_\l$. 

If $\l=3\e_i+\teps_6$ then $\de_\l=5-i$. If $\l=-3\e_i+\teps_6$ then $\de_\l=3+i$. For $\l=2\teps_6$ we have $\de_\l=0$ and  for $\l\in W_{D_5}(\frac{1}{2}(-3\e_1+3\e_2+3\e_3+3\e_4+3\e_5)+\frac{1}{2}\teps_6)$ we have  $\de_\l=h_{\l-\frac{1}{2}\teps_6}+1$.


\subsection{The root system \texorpdfstring{$E_7$}{E7}}
Denote $$\teps_7:=(-\e_7+\e_8)\in\Re^8 \quad \text{and} \quad V_7=\Re\e_1+\Re\e_2+\Re\e_3+\Re\e_4+\Re\e_5+\Re\e_6+\Re\teps_7\subset\Re^8$$

The Euclidean vector space $\left(\h^*_\Re, (\cdot,\cdot)\right)$ can be identified to $\left(V_7, (\cdot,\cdot)\right)$. The simple roots are $\a_1=\frac{1}{2}(\e_1-\e_2-\e_3-\e_4-\e_5-\e_6)+\frac{1}{2}\teps_7$, $\a_2=\e_1+\e_2$, $\a_i=-\e_{i-2}+\e_{i-1}$, $3\leq i\leq 7$. The dominant root is $\th=\teps_7$. The Weyl group orbit of an element of the root lattice cannot be described concisely in general. For this reason we will split an orbit into orbits under $W_{D_6}$, the parabolic Weyl group of type $D_6$ obtained by excluding the simple root $\a_1$. This group acts only on the first six coordinates of a vector in $\Re^8$ by permuting them and possibly changing an even number of signs. 

The positive roots are $\pm\e_i+\e_j$, $1\leq i<j\leq 6$, $\teps_7$, and the elements of the $W_{D_6}$ orbit of $\a_1$. The dominant non-zero small weights are $\teps_7$ (with symbol $[1]$),  $\e_5+\e_6+\teps_7$ (with symbol $[1^2]$), $2\e_6+\teps_7$ (with symbol $[1^3]$), $\frac{1}{2}(-\e_1+\e_2+\e_3+\e_4+\e_5+\e_6)+\frac{3}{2}\teps_7$ (with symbol $[2]$),  and $\frac{1}{2}(\e_1+\e_2+\e_3+\e_4+\e_5+3\e_6)+\frac{3}{2}\teps_7$ (with symbol $[2,1]$). We will describe separately the defect data for each of these weights. 
One common feature is that  $\ag_\l\not<{\bf 0}$ if the coefficient of $\teps_7$ in $\l$ is strictly negative.

Some of these weights are in fact in the root lattice of the parabolic root subsystem of type $D_6$ mentioned above. For these, it is convenient to use the detect data defined in Section 
\ref{II-Ddefectdata} and it is this data we refer to when we use  notation such as $\de_\l^{D_6}$. Note that the indexing of the nodes of the Dynkin diagram in Section  \ref{II-Ddefectdata} differs from the one induced here on the parabolic roots subsystem of type $D_6$. To obtain the correct formulas one needs to use $d_\l^{rev}$ in the formulas \eqref{II-Ddefect1} and \eqref{II-Ddefect2} and also $\tilde{d}^{rev}_\l$ as defined by \eqref{II-Ddefect2rev}. We will also make reference to the parabolic root system of type $E_6$ obtained by excluding $\a_7$.

\subsubsection{Weights with symbol $[1]$} These are  the roots. The defect vector is defined by $\de_\a=0$ if $\a\in R^+$ and $\de_\a=1$ if $\a\in R^-$.

\subsubsection{Weights with symbol $[1^2]$} We describe the  Weyl group orbit of $\e_5+\e_6+\teps_7$  by specifying the $W_{D_6}$ orbits it contains. These are the $W_{D_6}$ orbits of 
$\e_5+\e_6+\teps_7$, $2\e_6$, $\e_3+\e_4+\e_5+\e_6$, $\frac{1}{2}(\e_1+\e_2+\e_3+\e_4+\e_5+3\e_6)+\frac{1}{2}\teps_7$, $\e_5+\e_6-\teps_7$,  and $\frac{1}{2}(\e_1+\e_2+\e_3+\e_4+\e_5+3\e_6)-\frac{1}{2}\teps_6$. For the elements of the last two $W_{D_6}$ orbits  the coordinate of $\teps_6$ is strictly negative and  we set $\ag_\l\not<{\bf 0}$. We focus our attention on the first four orbits. 

For $\l$ in $W_{D_6}(\e_5+\e_6+\teps_7)$, the defect vector is defined by $\de_\l=(\de_{\l-\teps_7}^{D_6},0)$. For $\l$ in $W_{D_6}(2\e_6)$, the defect vector is defined by $\de_\l=\de_\l^{D_6}$. Let $\l$ in $W_{D_6}(\e_3+\e_4+\e_5+\e_6)$. Note that $|\A^-_2(\e_3+\e_4+\e_5+\e_6)|=2$. If $|I(\l)|>2$ then $\ag_\l\not<{\bf 0}$. Otherwise, define $$\de_\l=\de_\l^{D_6}+(2,0)$$

For  $\l$ in $W_{D_6}(\frac{1}{2}(\e_1+\e_2+\e_3+\e_4+\e_5+3\e_6)+\frac{1}{2}\teps_7)$ the defect vector is $\de_\l=(h^{(1)}_{\l-\frac{1}{2}\teps_7},0)$. The defect data presented in the section is nothing else but the one defined in the last paragraph of Section \ref{II-normal}.

\subsubsection{Weights with symbol $[1^3]$} The Weyl group orbit of $2\e_6+\teps_7$ consists of the  $W_{D_6}$ orbits of $2\e_6+\teps_7$, $\e_1+\e_2+\e_3+\e_4+\e_5+\e_6$, and 
$2\e_6-\teps_7$. The elements of the last  $W_{D_6}$ orbit have  the coordinate of $\teps_6$ strictly negative and  we set $\ag_\l\not<{\bf 0}$. 

For $\l$ in $W_{D_6}(2\e_6+\teps_7)$, the defect vector is defined by $\de_\l=(\de_{\l-\teps_7}^{D_6}(1),\de_{\l-\teps_7}^{D_6}(2),0)$. Let $\l$ in $W_{D_6}(\e_1+\e_2+\e_3+\e_4+\e_5+\e_6)$. Note that $|\A^-_2(\e_3+\e_4+\e_5+\e_6)|=6$. If $|I(\l)|>3$ then $\ag_\l\not<{\bf 0}$. Otherwise, define $$\de_\l=(\de_\l^{D_6}(2),\de_\l^{D_6}(1),\de_\l^{D_6}(3))+(6,0,0)$$
\subsubsection{Weights with symbol $[2]$} We describe the  Weyl group orbit of $\frac{1}{2}(-\e_1+\e_2+\e_3+\e_4+\e_5+\e_6)+\frac{3}{2}\teps_7$  by specifying the $W_{D_6}$ orbits it contains. These are the $W_{D_6}$ orbits of $\frac{1}{2}(-\e_1+\e_2+\e_3+\e_4+\e_5+\e_6)+\frac{3}{2}\teps_7$, $\e_3+\e_4+\e_5+\e_6+\teps_7$, $\e_4+\e_5+2\e_6$, $\frac{1}{2}(-\e_1+\e_2+\e_3+\e_4+3\e_5+3\e_6)+\frac{1}{2}\teps_7$, $-\e_1+\e_2+\e_3+\e_4+\e_5+\e_6$, $\frac{1}{2}(-\e_1+\e_2+\e_3+\e_4+\e_5+\e_6)-\frac{3}{2}\teps_7$, $\e_3+\e_4+\e_5+\e_6-\teps_7$, and $\frac{1}{2}(-\e_1+\e_2+\e_3+\e_4+3\e_5+3\e_6)-\frac{1}{2}\teps_7$.
For the elements of the last three $W_{6}$ orbits  the coordinate of $\teps_6$ is strictly negative and  we set $\ag_\l\not<{\bf 0}$. We focus our attention on the first five orbits. 

For  $\l$ in $W_{D_6}(\frac{1}{2}(-\e_1+\e_2+\e_3+\e_4+\e_5+\e_6)+\frac{3}{2}\teps_7)$ the defect data is defined by 
\begin{alignat}{3}\label{II-E7defect1}
\begin{split}
\de_\l&:= (0,0),\\
\tilde{d}_\l&:=0,
\end{split}
&\quad\quad&
\begin{split}
 \bar{\de}_\l&:= 0\\
\partial_\l&:=0
\end{split}
\end{alignat}

Let $\l$ in $W_{D_6}(\e_3+\e_4+\e_5+\e_6+\teps_7)$. Then $\l=\l_a\e_a+\l_b\e_b+\l_c\e_c+\l_d\e_d+\teps_7$ for some $1\leq a<b<c<d\leq 6$.
The defect vector is defined by  
\begin{alignat}{3}\label{II-E7defect2}
\begin{split}
\de_\l&:= (\parallel d^{rev}_{\l_c\e_c+\l_d\e_d}\parallel,0), \\
\tilde{d}_\l&:=\parallel d^{rev}_{\l-\teps_7}\parallel -\parallel d^{rev}_{\l_c\e_c+\l_d\e_d}\parallel,
\end{split}
&\quad\quad&
\begin{split}
\bar{\de}_\l&:= 0\\
\partial_\l&:=\parallel d^{rev}_{\l_c\e_c+\l_d\e_d}\parallel
\end{split}
\end{alignat}

Let $\l$ in $W_{D_6}(\e_4+\e_5+2\e_6)$. If $|I(\l^\sharp)|>2$ then set $\ag_\l\not<{\bf 0}$. Otherwise,  the defect data is defined by
\begin{alignat}{3}\label{II-E7defect3}
\begin{split}
\de_\l&:= \de_{\l}^{D_6}+(2,0),\\
\tilde{d}_\l&:=\tilde{d}_\l^{D_6}+2,\\
\end{split}
&\quad\quad&
\begin{split}
\bar{\de}_\l&:= \bar{\de}_{\l}^{D_6}\\
 \partial_\l&:=\partial^{D_6}_\l+2\\
\end{split}
\end{alignat}
Note that $|\A^{-,nr}_2(\e_4+\e_5+2\e_6)|=3$ and $|\A^{-,r}_2(\e_4+\e_5+2\e_6)|=1$ and this is reflected in the constants present in the formulas above. 

Let $\l$ in $W_{D_6}(\frac{1}{2}(-\e_1+\e_2+\e_3+\e_4+3\e_5+3\e_6)+\frac{1}{2}\teps_7)$. The defect data is defined by
\begin{equation}\label{II-E7defect4}
\begin{split}
\de_\l&:= (\min\{h^{(2)}_{\l-\frac{1}{2}\teps_7},\delta\}+2h_{\l},0),\quad\quad \bar{\de}_\l:=0\\
\tilde{d}_\l&:=h^{(1)}_{\l-\frac{1}{2}\teps_7}-\min\{h^{(2)}_{\l-\frac{1}{2}\teps_7},\delta\}+2h_{\l},\\
\partial_\l&:=\min\{h^{(2)}_{\l-\frac{1}{2}\teps_7},\delta\}
\end{split}
\end{equation}

Let $\l$ in $W_{D_6}(-\e_1+\e_2+\e_3+\e_4+\e_5+\e_6)$. If $|I(\l)|>3$ then set $\ag_\l\not<{\bf 0}$. Otherwise,  the defect data is defined by
\begin{alignat}{3}\label{II-E7defect5}
\begin{split}
\de_\l&:= (\de_\l^{D_6}(2), \de_\l^{D_6}(3))+(2,0),\\
\tilde{d}_\l&:=\de_\l^{D_6}(1)+2,\\
\end{split}
&\quad\quad&
\begin{split}
\bar{\de}_\l&:=\de_\l^{D_6}(3)\\
\partial_\l&:=\de_\l^{D_6}(2)\\
\end{split}
\end{alignat}
Note that  $|\A^-_3(-\e_1+\e_2+\e_3+\e_4+\e_5+\e_6)|=1$ and this is reflected in the constants present in the formulas above. 
\subsubsection{Weights with symbol $[2,1]$}  We describe the  Weyl group orbit of $\frac{1}{2}(\e_1+\e_2+\e_3+\e_4+\e_5+3\e_6)+\frac{3}{2}\teps_7$  by specifying the $W_{D_6}$ orbits it contains. These are  the  $W_{D_6}$ orbits of $\frac{1}{2}(\e_1+\e_2+\e_3+\e_4+\e_5+3\e_6)+\frac{3}{2}\teps_7$,  $\e_4+\e_5+2\e_6+\teps_7$, $\e_1+\e_2+\e_3+\e_4+\e_5+\e_6+\teps_7$, $\frac{1}{2}(\e_1+\e_2+\e_3+3\e_4+3\e_5+3\e_6)+\frac{1}{2}\teps_7$, $\frac{1}{2}(-\e_1+\e_2+\e_3+\e_4+\e_5+5\e_6)+\frac{1}{2}\teps_7$, $\e_2+\e_3+\e_4+\e_5+2\e_6$, $\frac{1}{2}(\e_1+\e_2+\e_3+\e_4+\e_5+3\e_6)-\frac{3}{2}\teps_7$,  $\e_4+\e_5+2\e_6-\teps_7$, $\e_1+\e_2+\e_3+\e_4+\e_5+\e_6-\teps_7$, $\frac{1}{2}(\e_1+\e_2+\e_3+3\e_4+3\e_5+3\e_6)-\frac{1}{2}\teps_7$, and $\frac{1}{2}(-\e_1+\e_2+\e_3+\e_4+\e_5+5\e_6)-\frac{1}{2}\teps_7$. For the elements of the last five $W_{D_6}$ orbits  the coordinate of $\teps_7$ is strictly negative  and we set $\ag_\l\not<{\bf 0}$.

To keep some of the formulas below more compact is convenient to sometime replace  the cut-off data in \eqref{II-e78shape} by its image under $\ddotbvs_2$. This is the case for \eqref{II-E7defect9}, \eqref{II-E7defect11}, and \eqref{II-E7defect13}. For all the other formulas the cut-off data is the one in \eqref{II-e78shape}.

Let  $\l$ in $W_{D_6}(\e_1+\e_2+\e_3+\e_4+\e_5+\e_6+\teps_7)$. The defect data is defined as follows. 
If $d^{rev}_{\l-\teps_7}(1)\leq 4$ then 
\begin{subequations}\label{II-E7defect21}
\begin{alignat}{3}\label{II-E7defect7}
\begin{split}
\de_\l&:= (\parallel d^{rev}_{\l-\teps_6}\parallel,0,0)  \\
\bar{\de}^\prime_\l&:= (\de_\l(1) -\tilde{d}_\l^\prime,0),\\
\partial_\l^\prime&:=0,\\
\tilde{d}_\l^\prime&:=d^{rev}_{\l-\teps_7}(1)+d^{rev}_{\l-\l_1\e_1+\l_2\e_2-\teps_6}(2),
\end{split}
&\quad&
\begin{split}
& \\
\bar{\de}^\dprime_\l&:= (\tilde{d}^\prime_\l,0)\\
\partial_\l^\dprime&:=0\\
\tilde{d}_\l^\dprime&:=\de_\l^\prime(1)
\end{split}
\end{alignat}
\text{If $d^{rev}_{\l-\teps_7}(1)=5$ then $\l$ lies in the parabolic subsystem of type $E_6$ and }
\begin{alignat}{3}\label{II-E7defect8}
\begin{split}
\de_\l&:= \de_\l^{E_6}+(5,1,0)\\
\bar{\de}^\prime_\l&:= \bar{\de}^{E_6}_\l+(1,0),\\
\partial^\prime_\l&:=\partial_\l^{E_6}+1,\\
\tilde{d}^\prime_\l&:=\tilde{d}_\l^{E_6}+4,
\end{split}
&\quad\quad&
\begin{split}
&\\
\bar{\de}^\dprime_\l&:= (5,0)\\
\partial^\dprime_\l&:=0\\
\tilde{d}^\dprime_\l&:=\parallel d^{rev}_{\l-\teps_7}\parallel-5
\end{split}
\end{alignat}
\end{subequations}

Let $\l$ in $W_{D_6}(\frac{1}{2}(-\e_1+\e_2+\e_3+\e_4+\e_5+5\e_6)+\frac{1}{2}\teps_7)$. The defect data  is defined as
\begin{equation}\label{II-E7defect9}
\begin{split}
\de_\l(1)&:= 2(\ind(\l)-1)\sgn(\l)+   2p_{\l-\frac{1}{2}\teps_7}+\min\{p_{\l-\frac{1}{2}\teps_7},\delta\}+o_{\l-\frac{1}{2}\teps_7}-\de_\l(2)\\  
\de_\l(2)&:=  2\left\lceil\frac{1}{2}(\ind(\l)-1)\right\rceil\sgn(\l)+   2\left\lfloor\frac{1}{2}\max\{p_{\l-\frac{1}{2}\teps_7}-\delta,0\}\right\rfloor+\min\{p_{\l-\frac{1}{2}\teps_7},\delta\}\\
\de_\l(3)&:=0\\
\bar{\de}^\prime_\l&:=(0,0),\quad \bar{\de}^\dprime_\l:=(\min\{p_{\l-\frac{1}{2}\teps_7},\delta\}+ o_{\l-\frac{1}{2}\teps_7},0)\\
\partial_\l^\prime&:=0,\quad \partial_\l^\dprime:=\min\{p_{\l-\frac{1}{2}\teps_7},\delta\}\\
 \tilde{d}^\prime_\l&:=1+\de_\l(1)+\de_\l(2)\\
\tilde{d}^\dprime_\l&:=1+2(\ind(\l)-1)\sgn(\l)+2p_{\l-\frac{1}{2}\teps_7}-\min\{p_{\l-\frac{1}{2}\teps_7},\delta\}
\end{split}
\end{equation}
The quantities that appear above have the following meaning: $2(\ind(\l)-1)\sgn(\l)+   2p_{\l-\frac{1}{2}\teps_6}$ equals $2|\A^-_3(\l)|$ and $o_{\l-\frac{1}{2}\teps_6}$ equals $|\A^{-,r}_2(\l)|$.
Note that  $|\A_2^{-,nr}(\frac{1}{2}(-\e_1+\e_2+\e_3+\e_4+\e_5+5\e_6)+\frac{1}{2}\teps_7)|=1$ and this is reflected in the constants present in the formulas above. 

Let $\l$ in $W_{D_6}(\frac{1}{2}(\e_1+\e_2+\e_3+\e_4+\e_5+3\e_6)+\frac{3}{2}\teps_7)$. The defect data is defined by
\begin{alignat}{3}\label{II-E7defect6}
\begin{split}
\de_\l&:= (h^{(1)}_{\l-\frac{3}{2}\teps_7},0,0)\\
\bar{\de}^\prime_\l&:=(0,0),\\
\partial_\l^\prime&:=0,\\
\tilde{d}^\prime_\l&:=0,
\end{split}
&\quad\quad&
\begin{split}
&\\
\bar{\de}^\dprime_\l&:=(h^{(1)}_{\l-\frac{3}{2}\teps_7},0)\\
\partial_\l^\dprime&:=0\\
\tilde{d}^\dprime_\l&:=0
\end{split}
\end{alignat}

Let $\l$ in $W_{D_6}(\e_2+\e_3+\e_4+\e_5+2\e_6)$. If $|I(\l^\sharp)|>3$ then set $\ag_\l\not<{\bf 0}$. Otherwise,  the defect data is defined by 
\begin{alignat}{3}\label{II-E7defect10}
\begin{split}
\de_\l&:= \de_{\l}^{D_6}+(7,2,0)\\
\bar{\de}^\prime_\l&:= \bar{\de}_{\l}^{D_6}+(2,0),\\
\partial^\prime_\l&:=\partial^{D_6}_\l+2,\\
\tilde{d}^\prime_\l&:=\tilde{d}_\l^{D_6}+4,
\end{split}
&\quad\quad&
\begin{split}
&\\
\bar{\de}^\dprime_\l&:= (5,0)\\
\partial^\dprime_\l&:=1\\
\tilde{d}^\dprime_\l&:=\tilde{d}_\l^{D_6}+\partial^{D_6}_\l+\parallel\bar{\de}_{\l}^{D_6}\parallel+1
\end{split}
\end{alignat}
Note that $|\A^{-,r}_2(\e_2+\e_3+\e_4+\e_5+2\e_6)|=4$ and $|\A^-_3(\e_2+\e_3+\e_4+\e_5+2\e_6)|=1$ and this is reflected in the constants present in the formulas above. 

Let $\l$ in $W_{D_6}(e_4+\e_5+2\e_6+\teps_7)$. The defect data  is defined as follows. If $\sgn(\l)=0$ then
\begin{subequations}\label{II-E7defect14}
\begin{alignat}{2}
\begin{split}\label{II-E7defect11}
\de_\l&:=(\sum_{\ind(\l)<j\leq 6}(1-\l_j), \max\{d^{rev}_{(\l-\teps_7)^\sharp}(1)-2, 0\}, 0)\\
\bar{\de}^\prime_\l&:=(0,0), \quad \bar{\de}^\dprime_\l:=(\de_\l(2)+\sum_{\ind(\l)<j\leq 6}(1-|\l_j|),0)\\
\partial_\l^\prime&:=0, \quad \partial^\dprime_\l:=\de_\l(2)\\
\tilde{d}^\prime_\l&:=\de_\l(1)+\de_\l(2), \quad \tilde{d}^\dprime_\l:=\tilde{d}^\prime_\l-\bar{\de}^\dprime_\l(1)
\end{split}
\intertext{If $\sgn(\l)=1$ and $\sum_{\ind(\l)<j\leq 6}\l_j=2$ then}
\begin{split}\label{II-E7defect12}
\de_\l&:=(\ind(\l)+2,0,0)\\
\bar{\de}^\prime_\l&:=(\ind(\l)-1,0), \quad \bar{\de}^\dprime_\l:=(3,0)\\
\partial_\l^\prime&:=0, \quad \partial^\dprime_\l:=0\\
\tilde{d}^\prime_\l&:=3, \quad \tilde{d}^\dprime_\l:=\ind(\l)-1
\end{split}
\intertext{Note that in this case $|\A^r_2(\l)|=3+(\ind(\l)-1)$.}
\intertext{If $\sgn(\l)=1$ and $\sum_{\ind(\l)<j\leq 6}\l_j\neq 2$ then }
\begin{split}\label{II-E7defect13}
\de_\l&:=(2\ind(\l)-3+\sum_{\ind(\l)<j\leq 6}(1-\l_j),1,0)\\
\bar{\de}^\prime_\l&:=(0,0), \quad \bar{\de}^\dprime_\l:=(2+\sum_{1\leq j<\ind(\l)}(1-|\l_j|),0)\\
\partial_\l^\prime&:=0, \quad \partial^\dprime_\l:=1\\
\tilde{d}^\prime_\l&:=\parallel\de_\l\parallel+\parallel d^{rev}_{\l^\flat-\teps_7}\parallel, \quad \tilde{d}^\dprime_\l:=\parallel d^{rev}_{(\l-\teps_7)^\sharp} \parallel-1
\end{split}
\end{alignat}
\end{subequations}

Let $\l$ in $W_{D_6}(\frac{1}{2}(\e_1+\e_2+\e_3+3\e_4+3\e_5+3\e_6)+\frac{1}{2}\teps_7)$. Note that $|\A^{-,r}_2(\frac{1}{2}(\e_1+\e_2+\e_3+3\e_4+3\e_5+3\e_6)+\frac{1}{2}\teps_7)|=3$ and this constant will be reflected in the formulas below. We set $\ag_\l\not<{\bf 0}$ unless $\l$ satisfies one of the conditions below. Otherwise,  the defect data is defined as follows.

If  $h_{\l}\leq 1$ the defect data is defined by 
\begin{subequations}\label{II-E7defect20}
\begin{alignat}{2}
\begin{split}\label{II-E7defect15}
\de_\l&:= (h^{(1)}_{\l-\frac{1}{2}\teps_7}+h^{(2)}_{\l-\frac{1}{2}\teps_7}+2h_{\l-\frac{1}{2}\teps_7},h^{(3)}_{\l-\frac{1}{2}\teps_7},0)+(3,0,0)\\
\bar{\de}^\prime_\l&:= (h^{(1)}_{\l-\frac{1}{2}\teps_7},0)+(2,0),\quad \bar{\de}^\dprime_\l:= (h^{(2)}_{\l-\frac{1}{2}\teps_7},0)+(1,0)\\
\partial^\prime_\l&:=h^{(3)}_{\l-\frac{1}{2}\teps_7},\quad \partial^\dprime_\l:=h^{(3)}_{\l-\frac{1}{2}\teps_7}\\
\tilde{d}^\prime_\l&:=3+h^{(2)}_{\l-\frac{1}{2}\teps_7}+2h_{\l-\frac{1}{2}\teps_7}, \quad \tilde{d}^\dprime_\l:=h^{(1)}_{\l-\frac{1}{2}\teps_7}+2h_{\l-\frac{1}{2}\teps_7}
\end{split}
\intertext{ If $h_{\l}=2$ and $h^{(3)}_{\l-\frac{1}{2}\teps_7}\leq 1$ the defect data is defined by }
\begin{split}\label{II-E7defect17}
\de_\l&:= (h^{(1)}_{\l-\frac{1}{2}\teps_7}+h^{(2)}_{\l-\frac{1}{2}\teps_7},h^{(3)}_{\l-\frac{1}{2}\teps_7},0)+(7,2,0)\\
\bar{\de}^\prime_\l&:= (h^{(1)}_{\l-\frac{1}{2}\teps_7},0)+(2,0),\quad \bar{\de}^\dprime_\l:= (h^{(2)}_{\l-\frac{1}{2}\teps_7},0)+(3,0)\\
\partial^\prime_\l&:=2+h^{(3)}_{\l-\frac{1}{2}\teps_7},\quad \partial^\dprime_\l:=h^{(3)}_{\l-\frac{1}{2}\teps_7}\\
\tilde{d}^\prime_\l&:=5+h^{(2)}_{\l-\frac{1}{2}\teps_7}, \quad \tilde{d}^\dprime_\l:=4+h^{(1)}_{\l-\frac{1}{2}\teps_7}
\end{split}
\intertext{If  $h_{\l}= 2$, and  $h^{(i)}_{\l-\frac{1}{2}\teps_7}=2$, $i\leq 3$  the defect data is defined by }
\begin{split}\label{II-E7defect18}
\de_\l&:= (11,3,0)\\
\bar{\de}^\prime_\l&:= (3,0),\quad \bar{\de}^\dprime_\l:=(3,0) \\
\partial^\prime_\l&:=3,\quad \partial^\dprime_\l:=3\\
\tilde{d}_\l^\prime&:=8, \quad \tilde{d}_\l^\dprime:=7
\end{split}
\intertext{If  $h_{\l}= 3$,   and $h^{(i)}_{\l-\frac{1}{2}\teps_7}\leq 3-i$, $i\leq 3$  the defect data is defined by }
\begin{split}\label{II-E7defect19}
\de_\l&:= (h^{(1)}_{\l-\frac{1}{2}\teps_7}+h^{(2)}_{\l-\frac{1}{2}\teps_7},0,0)+(6,4,0)\\
\bar{\de}^\prime_\l&:=(h^{(1)}_{\l-\frac{1}{2}\teps_7}+h^{(2)}_{\l-\frac{1}{2}\teps_7},0),\quad \bar{\de}^\dprime_\l:=(4, 0) \\
\partial^\prime_\l&:=2,\quad \partial^\dprime_\l:=2\\
\tilde{d}_\l^\prime&:=h^{(2)}_{\l-\frac{1}{2}\teps_7}+6, \quad \tilde{d}_\l^\dprime:=h^{(1)}_{\l-\frac{1}{2}\teps_7}+3
\end{split}
\end{alignat}
\end{subequations}


\subsection{The root system \texorpdfstring{$E_8$}{E8}}
The Euclidean vector space $\left(\h^*_\Re, (\cdot,\cdot)\right)$ can be identified to $\left(\Re^8, (\cdot,\cdot)\right)$. The simple roots are $\a_1=\frac{1}{2}(\e_1-\e_2-\e_3-\e_4-\e_5-\e_6-\e_7)+\frac{1}{2}\e_8$, $\a_2=\e_1+\e_2$, $\a_i=-\e_{i-2}+\e_{i-1}$, $3\leq i\leq 8$. The dominant root is $\th=\e_7+\e_8$. The Weyl group orbit of an element of the root lattice cannot be described concisely in general. For this reason we will split an orbit into orbits under $W_{D_7}$, the parabolic Weyl group of type $D_7$ obtained by excluding the simple root $\a_1$. This group acts only on the first seven coordinates of a vector in $\Re^8$ by permuting them and possibly changing an even number of signs. 

The positive roots are $\pm\e_i+\e_j$, $1\leq i<j\leq 8$, and the elements of the $W_{D_7}$ orbit of $\a_1$. The dominant non-zero small weights are $\e_7+\e_8$ (with symbol $[1]$),  $2\e_8$ (with symbol $[1^2]$), $\e_6+\e_7+2\e_8$ (with symbol $[2]$), and $\frac{1}{2}(\e_1+\e_2+\e_3+\e_4+\e_5+\e_6+\e_7)+\frac{5}{2}\e_8$  (with symbol $[2,1]$). We will describe separately the defect data for each of these weights. One common feature is that  $\ag_\l\not<{\bf 0}$ if the coefficient of $\e_8$ in $\l$ is strictly negative.

Some of these weights are in fact in the root lattice of the parabolic root subsystem of type $D_7$ mentioned above. For these, it is convenient to use the detect data defined in Section 
\ref{II-Ddefectdata} and it is this data we refer to when we use  notation such as $\de_\l^{D_7}$. Note that the indexing of the nodes of the Dynkin diagram in Section  \ref{II-Ddefectdata} differs from the one induced here on the parabolic roots subsystem of type $D_7$. To obtain the correct formulas one needs to use $d_\l^{rev}$ in the formulas \eqref{II-Ddefect1} and \eqref{II-Ddefect2} and also $\tilde{d}^{rev}_\l$ as defined by \eqref{II-Ddefect2rev}. We will also make reference to the parabolic root system of type $E_7$ obtained by excluding $\a_8$.

\subsubsection{Weights with symbol $[1]$} These are  the roots. The defect vector is defined by $\de_\a=0$ if $\a\in R^+$ and $\de_\a=1$ if $\a\in R^-$.

\subsubsection{Weights with symbol $[1^2]$} We describe the  Weyl group orbit of $2\e_8$  by specifying the $W_{D_7}$ orbits it contains. These are the $W_{D_7}$ orbits of 
$2\e_8$, $2\e_7$,  $\e_4+\e_5+\e_6+\e_7$,  $\e_5+\e_6+\e_7+\e_8$, $\frac{1}{2}(-\e_1+\e_2+\e_3+\e_4+\e_5+\e_6+\e_7)+\frac{3}{2}\e_8$, $\frac{1}{2}(-\e_1+\e_2+\e_3+\e_4+\e_5+\e_6+3\e_7)+\frac{1}{2}\e_8$, $-2\e_8$, $\e_5+\e_6+\e_7-\e_8$, $\frac{1}{2}(\e_1+\e_2+\e_3+\e_4+\e_5+\e_6+\e_7)-\frac{3}{2}\e_8$, and $\frac{1}{2}(\e_1+\e_2+\e_3+\e_4+\e_5+\e_6+3\e_7)-\frac{1}{2}\e_8$. For the elements of the last four $W_{D_7}$ orbits  the coordinate of $\e_8$ is strictly negative and  we set $\ag_\l\not<{\bf 0}$. We focus our attention on the first six orbits. 

For $\l=2\e_8$ and $\l\in W_{D_7}(\frac{1}{2}(-\e_1+\e_2+\e_3+\e_4+\e_5+\e_6+\e_7)+\frac{3}{2}\e_8)$ the defect vector is defined by $\de_\l=(0,0)$. For $\l$ in $W_{D_7}(2\e_7)$, the defect vector is defined by $\de_\l=\de_{\l}^{D_7}+(1,0)$.  Let $\l$ in $W_{D_7}(\e_4+\e_5+\e_6+\e_7)$. Note that $|\A^-_2(\e_3+\e_4+\e_5+\e_6)|=4$. If $|I(\l)|>2$ then $\ag_\l\not<{\bf 0}$. Otherwise, define $$\de_\l=\de_\l^{D_7}+(4,0)$$

For  $\l$ in $W_{D_7}(\frac{1}{2}(-\e_1+\e_2+\e_3+\e_4+\e_5+\e_6+3\e_7)+\frac{1}{2}\e_8)$ the defect vector is $$\de_\l=(h^{(1)}_{\l-\frac{1}{2}\e_8}+1,0)$$

Let $\l$ in $W_{D_7}(\e_5+\e_6+\e_7+\e_8)$. The defect vector is defined by  
$$
\de_\l:= (\parallel d^{rev}_\l\parallel,0)
$$

The defect data presented in the section is nothing else but the one defined in the last paragraph of Section \ref{II-normal}.

\subsubsection{Weights with symbol $[2]$} We describe the  Weyl group orbit of $\e_6+\e_7+2\e_8$  by specifying the $W_{D_7}$ orbits it contains. These are the $W_{D_7}$ orbits of 
$\e_6+\e_7+2\e_8$, $\e_5+\e_6+2\e_7$, $\e_6+2\e_7+\e_8$, $\e_3+\e_4+\e_5+\e_6+\e_7+\e_8$, $\e_2+\e_3+\e_4+\e_5+\e_6+\e_7$, $\frac{1}{2}(\e_1+\e_2+\e_3+\e_4+\e_5+\e_6+3\e_7)+\frac{3}{2}\e_8$, $\frac{1}{2}(\e_1+\e_2+\e_3+\e_4+\e_5+3\e_6+3\e_7)+\frac{1}{2}\e_8$, $\e_6+\e_7-2\e_8$, $\e_6+2\e_7-\e_8$, $\e_3+\e_4+\e_5+\e_6+\e_7-\e_8$, $\frac{1}{2}(-\e_1+\e_2+\e_3+\e_4+\e_5+\e_6+3\e_7)-\frac{3}{2}\e_8$, and $\frac{1}{2}(-\e_1+\e_2+\e_3+\e_4+\e_5+3\e_6+3\e_7)-\frac{1}{2}\e_8$.

For the elements of the last five $W_{D_7}$ orbits  the coordinate of $\e_8$ is strictly negative and  we set $\ag_\l\not<{\bf 0}$. We focus our attention on the first seven orbits. 

Let $\l$ in $W_{D_7}(\e_6+\e_7+2\e_8)$. The defect vector is defined by  
\begin{alignat}{3}\label{II-E8defect1}
\begin{split}
\de_\l&:= (0,0),\\
\tilde{d}_\l&:=\parallel d^{rev}_{\l^\flat}\parallel,
\end{split}
&\quad\quad&
\begin{split}
\bar{\de}_\l&:= 0\\
\partial_\l&:=0
\end{split}
\end{alignat}

Let $\l$ in $W_{D_7}(\e_5+\e_6+2\e_7)$. If $|I(\l^\sharp)|>2$ then set $\ag_\l\not<{\bf 0}$. Otherwise,  the defect data is defined by
\begin{alignat}{3}\label{II-E8defect2}
\begin{split}
\de_\l&:= \de_{\l}^{D_7}+(4,0),\\
\tilde{d}_\l&:=\tilde{d}_\l^{D_7}+5,
\end{split}
&\quad\quad&
\begin{split}
\bar{\de}_\l&:= \bar{\de}_{\l}^{D_7}\\
\partial_\l&:=\partial^{D_7}_\l+4
\end{split}
\end{alignat}
Note that $|\A^{-,nr}_2(\e_4+\e_5+2\e_6)|=9$ and this is reflected in the constants present in the formulas above. 

Let $\l$ in $W_{D_7}(\e_6+2\e_7+\e_8)$. The defect data is defined by 
\begin{equation}\label{II-E8defect3}
\begin{split}
\de_\l&:= (\parallel d^{rev}_{\l^\sharp} \parallel, 0),\quad \bar{\de}_\l:= 0\\
\partial_\l&:=\sgn(\l)\\
\tilde{d}_\l&:=\sgn(\l)(2(\ind(\l)-1)-1)+\sum_{\ind(\l)<j\leq 7} (1-\l_j)
\end{split}
\end{equation}

Let $\l$ in $W_{D_7}(\e_3+\e_4+\e_5+\e_6+\e_7+\e_8)$. Then $\l=\l_a\e_a+\l_b\e_b+\l_c\e_c+\l_d\e_d+\l_e\e_e+\e_8$ for some $1\leq a<b<c<d<e\leq 7$.
The defect vector is defined by  
\begin{alignat}{3}\label{II-E8defect4}
\begin{split}
\de_\l&:= (\parallel d^{rev}_{\l_c\e_c+\l_d\e_d+\l_e\e_e}\parallel,0),\\
\tilde{d}_\l&:=\parallel d^{rev}_{\l}\parallel-\de_\l(1)+2,
\end{split}
&\quad\quad&
\begin{split}
\bar{\de}_\l&:= 0\\
\partial_\l&:=\de_\l(1)
\end{split}
\end{alignat}
Note that $|\A^{-,nr}_2(\e_3+\e_4+\e_5+\e_6+\e_7+\e_8)|=2$ and this is reflected in the constants present in the formulas above. 

Let $\l$ in $W_{D_7}(\frac{1}{2}(\e_1+\e_2+\e_3+\e_4+\e_5+\e_6+3\e_7)+\frac{3}{2}\e_8)$. The defect data is defined by
\begin{alignat}{3}\label{II-E8defect5}
\begin{split}
\de_\l&:= (0,0),\\
\tilde{d}_\l&:=h^{(1)}_{\l},
\end{split}
&\quad\quad&
\begin{split}
\bar{\de}_\l&:=0\\
\partial_\l&:=0\\
\end{split}
\end{alignat}

Let $\l$ in $W_{D_7}(\e_2+\e_3+\e_4+\e_5+\e_6+\e_7)$. If $|I(\l)|>3$ then set $\ag_\l\not<{\bf 0}$. Otherwise,  the defect data is defined by
\begin{alignat}{3}\label{II-E8defect6}
\begin{split}
\de_\l&:= (\de_\l^{D_7}(1), \de_\l^{D_7}(3))+(2,0),\\
\tilde{d}_\l&:=\de_\l^{D_7}(2)+8,
\end{split}
&\quad\quad&
\begin{split}
\bar{\de}_\l&:=\de_\l^{D_7}(3)\\
\partial_\l&:=\de_\l^{D_7}(1)
\end{split}
\end{alignat}
Note that  $|\A^-_3(\e_2+\e_3+\e_4+\e_5+\e_6+\e_7+\e_8)|=1$ and $|\A^{-,nr}_2(\e_2+\e_3+\e_4+\e_5+\e_6+\e_7+\e_8)|=6$ and this is reflected in the constants present in the formulas above. 

Let $\l$ in $W_{D_7}(\frac{1}{2}(\e_1+\e_2+\e_3+\e_4+\e_5+3\e_6+3\e_7)+\frac{1}{2}\e_8)$. The defect data is defined by
\begin{alignat}{3}\label{II-E8defect7}
\begin{split}
\de_\l&:= (\min\{h^{(2)}_{\l},\delta\}+2h_{\l},0),\\
\tilde{d}_\l&:=h^{(1)}_{\l}-\min\{h^{(2)}_{\l},\delta\}+2h_{\l},
\end{split}
&\quad\quad&
\begin{split}
\bar{\de}_\l&:=0\\
\partial_\l&:=\de_\l(1)
\end{split}
\end{alignat}
Note that   $|\A^{-,nr}_2(\frac{1}{2}(\e_1+\e_2+\e_3+\e_4+\e_5+3\e_6+3\e_7)+\frac{1}{2}\e_8)|=5$ and this is reflected in the constants present in the formulas above. 

\subsubsection{Weights with symbol $[2,1]$}   We describe the  Weyl group orbit of $\frac{1}{2}(\e_1+\e_2+\e_3+\e_4+\e_5+\e_6+\e_7)+\frac{5}{2}\e_8$  by specifying the $W_{D_7}$ orbits it contains. These are  the  $W_{D_7}$ orbits of $\frac{1}{2}(\e_1+\e_2+\e_3+\e_4+\e_5+\e_6+\e_7)+\frac{5}{2}\e_8$,  $\frac{1}{2}(-\e_1+\e_2+\e_3+\e_4+\e_5+3\e_6+3\e_7)+\frac{3}{2}\e_8$, $\frac{1}{2}(-\e_1+\e_2+\e_3+\e_4+3\e_5+3\e_6+3\e_7)+\frac{1}{2}\e_8$, $\frac{1}{2}(\e_1+\e_2+\e_3+\e_4+\e_5+\e_6+5\e_7)+\frac{1}{2}\e_8$, $\e_4+\e_5+\e_6+\e_7+2\e_8$, $\e_4+\e_5+\e_6+2\e_7+\e_8$, $-\e_1+\e_2+\e_3+\e_4+\e_5+\e_6+\e_7+\e_8$,  $\e_3+\e_4+\e_5+\e_6+2\e_7$, $\frac{1}{2}(-\e_1+\e_2+\e_3+\e_4+\e_5+\e_6+\e_7)-\frac{5}{2}\e_8$,  $\frac{1}{2}(\e_1+\e_2+\e_3+\e_4+\e_5+3\e_6+3\e_7)-\frac{3}{2}\e_8$, $\frac{1}{2}(\e_1+\e_2+\e_3+\e_4+3\e_5+3\e_6+3\e_7)-\frac{1}{2}\e_8$, $\frac{1}{2}(-\e_1+\e_2+\e_3+\e_4+\e_5+\e_6+5\e_7)-\frac{1}{2}\e_8$, $\e_4+\e_5+\e_6+\e_7-2\e_8$, $\e_4+\e_5+\e_6+2\e_7-\e_8$, and $\e_1+\e_2+\e_3+\e_4+\e_5+\e_6+\e_7-\e_8$.
 For the elements of the last seven $W_{D_7}$ orbits  the coordinate of $\e_8$ is strictly negative  and we set $\ag_\l\not<{\bf 0}$. 

Let $\l$ in $W_{D_7}(\frac{1}{2}(\e_1+\e_2+\e_3+\e_4+\e_5+\e_6+\e_7)+\frac{5}{2}\e_8)$. The defect data is defined by
\begin{alignat}{3}\label{II-E8defect8}
\begin{split}
\de_\l&:= (0,0,0)\\
\bar{\de}^\prime_\l&:=(0,0),\\
\partial_\l^\prime&:=0,\\
\tilde{d}^\prime_\l&:=0,
\end{split}
&\quad\quad&
\begin{split}
&\\
\bar{\de}^\dprime_\l&:=(0,0)\\
\partial_\l^\dprime&:=0\\
\tilde{d}^\dprime_\l&:=0
\end{split}
\end{alignat}

Let $\l$ in $W_{D_7}(\frac{1}{2}(\e_1+\e_2+\e_3+\e_4+\e_5+\e_6+5\e_7)+\frac{1}{2}\e_8)$. The defect data  is defined by
\begin{equation}\label{II-E8defect9}
\begin{split}
\de_\l(1)&:= 2(\ind(\l)-1)\sgn(\l)+   2p_{\l}+\min\{p_{\l},\delta\}+o_{\l-\frac{1}{2}\e_8}-\de_\l(2)+8\\  
\de_\l(2)&:=  2\left\lceil\frac{1}{2}(\ind(\l)-1)\right\rceil\sgn(\l)+   2\left\lfloor\frac{1}{2}\max\{p_{\l}-\delta,0\}\right\rfloor+\min\{p_{\l},\delta\}+1\\
\de_\l(3)&:=0\\
\bar{\de}^\prime_\l(1)&:=\min\{p_{\l},\delta\}+ o_{\l-\frac{1}{2}\e_8}+1, \quad \bar{\de}^\prime_\l(2):=0,\quad \bar{\de}^\dprime_\l:=(7,0)\\
\partial_\l^\prime&:=\min\{p_{\l},\delta\}+1, \quad \partial_\l^\dprime:=0\\
\tilde{d}^\prime_\l&:=6+2(\ind(\l)-1)\sgn(\l)+2p_{\l}-\min\{p_{\l},\delta\},\quad \tilde{d}^\dprime_\l:=\de_\l(1)+\de_\l(2)
\end{split}
\end{equation}
The quantities that appear above have the following meaning: $2(\ind(\l)-1)\sgn(\l)+   2p_{\l}$ equals $2|\A^-_3(\l)|$ and $o_{\l}+7$ equals $|\A^-_2(\l)|$.
Note that  $|\A_2^{-,r}(\frac{1}{2}(-\e_1+\e_2+\e_3+\e_4+\e_5+\e_6+5\e_7)+\frac{1}{2}\e_8)|=7$ and this is reflected in the constants present in the formulas above. 

Let $\l$ in $W_{D_7}(\e_4+\e_5+\e_6+\e_7+2\e_8)$. Then $\l=\l_a\e_a+\l_b\e_b+\l_c\e_c+\l_d\e_d+2\e_8$ for some $1\leq a<b<c<d\leq 7$. The defect data  is defined as 
\begin{equation}\label{II-E8defect10}
\begin{split}
\de_\l&:=(\parallel d^{rev}_{\l-2\e_8}\parallel,0 , 0) \\  
\bar{\de}^\prime_\l&:=(\parallel d^{rev}_{\l-\l_a\e_a-\l_b\e_b-2\e_8}\parallel, 0),\quad \bar{\de}^\dprime_\l:=(\de_\l(1)-\bar{\de}^\prime_\l(1),0)\\
\partial_\l^\prime&:=0, \quad \partial_\l^\dprime:=0\\
\tilde{d}^\prime_\l&:=\bar{\de}^\dprime_\l,\quad \tilde{d}^\dprime_\l:=\bar{\de}^\prime_\l
\end{split}
\end{equation}

Let $\l$ in $W_{D_7}(\frac{1}{2}(-\e_1+\e_2+\e_3+\e_4+\e_5+3\e_6+3\e_7)+\frac{3}{2}\e_8)$. The defect data is defined by
\begin{equation}\label{II-E8defect11}
\begin{split}
\de_\l&:= (h^{(1)}_{\l}+h_\l^{(2)}+2h_\l+1,0,0)\\
\bar{\de}^\prime_\l&:=(h_\l^{(2)},0),\quad \bar{\de}^\dprime_\l:=(h_\l^{(1)}+1,0)\\
\partial_\l^\prime&:=0,\quad  \partial_\l^\dprime:=0\\
\tilde{d}^\prime_\l&:=\bar{\de}^\dprime_\l(1)+2h_\l,\quad \tilde{d}^\dprime_\l:=\bar{\de}^\prime_\l(1)+2h_\l
\end{split}
\end{equation}
Note that  $|\A_2^-(\frac{1}{2}(-\e_1+\e_2+\e_3+\e_4+\e_5+3\e_6+3\e_7)+\frac{3}{2}\e_8)|=1$ and this is reflected in the constants present in the formulas above. 

For  $\l$ in $W_{D_7}(-\e_1+\e_2+\e_3+\e_4+\e_5+\e_6+\e_7+\e_8)$ define the defect data as follows (need to specify the cut off vector too).
If $d^{rev}_{\l}(1)\leq 5$ then 
\begin{subequations}\label{II-E8defect25}
\begin{equation}\label{II-E8defect12}
\begin{split}
\de_\l&:= (\parallel d^{rev}_{\l}\parallel-\parallel d^{rev}_{\l_5\e_5+\l_6\e_6}\parallel,\parallel d^{rev}_{\l_5\e_5+\l_6\e_6}\parallel,0)+(2,0,0)\\
\bar{\de}^\prime_\l&:= ( \tilde{d}_\l^\dprime-2,0)  \\
\bar{\de}^\dprime_\l&:= (d^{rev}_{\l}(1)-\de_\l(2)+d^{rev}_{\l-\l_1\e_1-\l_2\e_2}(2)+d^{rev}_{\l-\l_1\e_1-\l_2\e_2}(3),0)\\
\partial_\l^\prime&:=\parallel d^{rev}_{\l_5\e_5+\l_6\e_6}\parallel,\quad  \partial_\l^\dprime:=\parallel d^{rev}_{\l_5\e_5+\l_6\e_6}\parallel\\
\tilde{d}_\l^\prime&:=\bar{\de}^\dprime_\l(1)+2,\quad  \tilde{d}_\l^\dprime:=\de_\l(1)-\bar{\de}_\l^\dprime(1)
\end{split}
\end{equation}
\text{If $d^{rev}_{\l}=6$ then $\l$ lies in the parabolic subsystem of type $E_7$ and }
\begin{alignat}{3}\label{II-E8defect13}
\begin{split}
\de_\l&:= \de_\l^{E_7}+(6,2,0)\\
\bar{\de}^\prime_\l&:= \bar{\de}^{\prime E_7}_\l+(2,0),\\
\partial^\prime_\l&:=\partial_\l^{\prime E_7}+2,\\
\tilde{d}^\prime_\l&:=\tilde{d}_\l^{\prime E_7}+4,
\end{split}
&\quad\quad&
\begin{split}
&\\
\bar{\de}^\dprime_\l&:=\bar{\de}^{\dprime E_7}_\l+(2,0)\\
\partial^\dprime_\l&:=\partial_\l^{\dprime E_7}+2\\
\tilde{d}^\dprime_\l&:=\tilde{d}_\l^{\dprime E_7}+4
\end{split}
\end{alignat}
\end{subequations}
Note that  $|\A_3^-(-\e_1+\e_2+\e_3+\e_4+\e_5+\e_6+\e_7+\e_8)|=1$ and this is reflected in the constants present in the formulas above.

Let $\l$ in $W_{D_7}(\e_3+\e_4+\e_5+\e_6+2\e_7)$. If $|I(\l^\sharp)|>3$ then set $\ag_\l\not<{\bf 0}$. Otherwise,  the defect data is defined by
\begin{alignat}{3}\label{II-E8defect14}
\begin{split}
\de_\l&:= \de_{\l}^{D_7}+(13,4,0)\\
\bar{\de}^\prime_\l&:= \bar{\de}_{\l}^{D_7}+(4,0),\\
\partial^\prime_\l&:=\partial^{D_7}_\l+4,\\
\tilde{d}^\prime_\l&:=\tilde{d}_\l^{D_7}+7,
\end{split}
&\quad\quad&
\begin{split}
&\\
\bar{\de}^\dprime_\l&:= (7,0)\\
\partial^\dprime_\l&:=2\\
\tilde{d}^\dprime_\l&:=\parallel \bar{\de}_{\l}^{D_7}\parallel+\tilde{d}_\l^{D_7}+\partial^{D_7}_\l+6
\end{split}
\end{alignat}
Note that $|\A^{-r}_2(\e_3+\e_4+\e_5+\e_6+2\e_7)|=9$ and $|\A^-_3(\e_3+\e_4+\e_5+\e_6+2\e_7)|=2$ and this is reflected in the constants present in the formulas above. 

Let $\l$ in $W_{D_7}(\e_4+\e_5+\e_6+2\e_7+\e_8)$. Note that $|\A_2^-(\e_4+\e_5+\e_6+2\e_7+\e_8)|=4$ and this is reflected in the constants present in the formulas below.
The defect data  is defined as follows. If $\sgn(\l)=0$ then
\begin{subequations}\label{II-E8defect18}
\begin{alignat}{2}
\begin{split}\label{II-E8defect15}
\de_\l(1)&:=\parallel d^{rev}_{\l^\sharp}\parallel + \sum_{\ind(\l)<j\leq 7}(1-|\l_j|)+4,\quad \de_\l(3):=0 \\
\de_\l(2)&:=\max\{d^{rev}_{\l^\sharp}(1)-2, 0\}+\max\{d^{rev}_{\l^\sharp}(2)-2, 0\}\\
\bar{\de}^\prime_\l&:=(\de_\l(2)+\sum_{\ind(\l)<j\leq 7}(1-|\l_j|),0), \quad \bar{\de}^\dprime_\l:=(\parallel d^{rev}_{\l^\flat}\parallel+4,0)\\
\partial_\l^\prime&:=\de_\l(2), \quad \partial^\dprime_\l:=0\\
\tilde{d}^\prime_\l&:=\de_\l(1)-\bar{\de}^\prime_\l(1), \quad \tilde{d}^\dprime_\l:=\de_\l(1)-\bar{\de}^\dprime_\l(1)
\end{split}
\intertext{If $\sgn(\l)=1$ and $\sum_{\ind(\l)<j\leq 7}\l_j=3$ then }
\begin{split}\label{II-E8defect16}
\de_\l&:=(\ind(\l)+6,0,0)\\
 \bar{\de}^\prime_\l&:=(3,0), \quad \bar{\de}^\dprime_\l:=(\ind(\l)+3,0)\\
\partial_\l^\prime&:=0, \quad \partial^\dprime_\l:=0\\
\tilde{d}^\prime_\l&:=\ind(\l)+3, \quad\tilde{d}^\dprime_\l:=3
\end{split}
\intertext{Note that in this case $|\A^{-,r}_2(\l)|=7+(\ind(\l)-1)$. }
\intertext{If $\sgn(\l)=1$ and $\sum_{\ind(\l)<j\leq 7}\l_j\neq 3$ then $$\l^\sharp=\l_a\e_a+\l_b\e_b+\l_c\e_c+\l_d\e_d+\l_e\e_e+\e_8$$ for some $1\leq a<b<c<d<e\leq 7$. The defect data in this case is defined as}
\begin{split}\label{II-E8defect17}
\de_\l&:=(\tilde{d}^\dprime_\l+\bar{\de}^\dprime_\l(1)+2, \partial^\dprime_\l, 0)\\
\bar{\de}^\prime_\l&:=(\parallel d^{rev}_{\l^\flat}\parallel+4,0), \quad \bar{\de}^\dprime_\l:=(1+\partial^\dprime_\l+\sum_{ j<\ind(\l)}(1-|\l_j|),0)\\
\partial_\l^\prime&:=1, \quad \partial^\dprime_\l:=1+\max\{\parallel d^{rev}_{\l_c\e_c+\l_d\e_d+\l_e\e_e}\parallel-1,0\}\\
\tilde{d}^\prime_\l&:=\de_\l(1)-\bar{\de}^\prime_\l(1)+2, \quad \tilde{d}^\dprime_\l:=\parallel d^{rev}_{\l^\sharp} \parallel-\partial^\dprime_\l+2
\end{split}
\end{alignat}
\end{subequations}

Let $\l$ in $W_{D_7}(\frac{1}{2}(-\e_1+\e_2+\e_3+\e_4+3\e_5+3\e_6+3\e_7)+\frac{1}{2}\e_8)$. Note that $|\A^{-}_2(\frac{1}{2}(-\e_1+\e_2+\e_3+\e_4+3\e_5+3\e_6+3\e_7)+\frac{1}{2}\e_8)|=6$ and $|\A^{-}_3(\frac{1}{2}(-\e_1+\e_2+\e_3+\e_4+3\e_5+3\e_6+3\e_7)+\frac{1}{2}\e_8)|=1$ and these constants will be reflected in the formulas below. We set $\ag_\l\not<{\bf 0}$ unless $\l$ satisfies one of the conditions below. Otherwise,  the defect data is defined as follows.

If  $h_{\l}\leq 1$ the defect data is defined by 
\begin{subequations}\label{II-E8defect24}
\begin{alignat}{2}
\begin{split}\label{II-E8defect19}
\de_\l&:= (h^{(1)}_{\l}+h^{(2)}_{\l}+2h_{\l},h^{(3)}_{\l},0)+(7,1,0)\\
\bar{\de}^\prime_\l&:= (h^{(1)}_{\l},0)+(3,0),\quad \bar{\de}^\dprime_\l:= (h^{(2)}_{\l},0)+(2,0)\\
\partial^\prime_\l&:=h^{(3)}_{\l}+1,\quad \partial^\dprime_\l:=h^{(3)}_{\l}+1\\
\tilde{d}^\prime_\l&:=h^{(2)}_{\l}+2h_{\l}+6, \quad \tilde{d}^\dprime_\l:=h^{(1)}_{\l}+2h_{\l}+3
\end{split}
\intertext{ If $h^{(3)}_{\l}\leq 2$ and $h_{\l}=2$ the defect data is defined by }
\begin{split}\label{II-E8defect20}
\de_\l&:= (h^{(1)}_{\l}+h^{(2)}_{\l},h^{(3)}_{\l},0)+(11,3,0)\\
\bar{\de}^\prime_\l&:= (h^{(1)}_{\l},0)+(3,0),\quad \bar{\de}^\dprime_\l:= (h^{(2)}_{\l},0)+(4,0)\\
\partial^\prime_\l&:=h^{(3)}_{\l}+3,\quad \partial^\dprime_\l:=h^{(3)}_{\l}+1\\
\tilde{d}^\prime_\l&:=h^{(2)}_{\l}+8, \quad \tilde{d}^\dprime_\l:=h^{(1)}_{\l}+7
\end{split}
\intertext{If  $h_{\l}= 2$,  $h^{(3)}_{\l}=3$ and $h^{(1)}_{\l}=3$, $h^{(2)}_{\l}=3$ the defect data is defined by }
\begin{split}\label{II-E8defect21}
\de_\l&:= (17,5,0)\\
\bar{\de}^\prime_\l&:= (5,0),\quad \bar{\de}^\dprime_\l:=(5,0) \\
\partial^\prime_\l&:=5,\quad \partial^\dprime_\l:=5\\
\tilde{d}_\l^\prime&:=12, \quad \tilde{d}_\l^\dprime:=11
\end{split}
\intertext{If  $h_{\l}= 3$,  $h^{(3)}_{\l}=0$ and $h^{(2)}_{\l}\leq 2$ the defect data is defined by }
\begin{split}\label{II-E8defect22}
\de_\l&:= (\bar{\de}^\prime_\l(1)+\tilde{d}_\l^\prime, \partial^\prime_\l+2,0)\\
\bar{\de}^\prime_\l&:=( h^{(1)}_{\l}+3,0),\quad \bar{\de}^\dprime_\l:=(\min(h^{(2)}_\l,1)+6, 0) \\
\partial^\prime_\l&:=3+h^{(2)}_{\l}-\min(h^{(2)}_\l,1),\quad \partial^\dprime_\l:=2\\
\tilde{d}_\l^\prime&:=\min(h^{(2)}_\l,1)+10, \quad \tilde{d}_\l^\dprime:=\bar{\de}^\prime_\l(1)+\partial^\prime_\l+1
\end{split}
\intertext{If  $h_{\l}= 3$,  $h^{(3)}_{\l}=1$ and $h^{(2)}_{\l}\leq 2$,  $h^{(1)}_{\l}\leq 3$  the defect data is defined by }
\begin{split}\label{II-E8defect23}
\de_\l&:= (h^{(1)}_{\l}+h^{(2)}_{\l}+12,6,0)\\
\bar{\de}^\prime_\l&:=(h^{(1)}_{\l}+h^{(2)}_{\l}+2,0),\quad \bar{\de}^\dprime_\l:=(6, 0) \\
\partial^\prime_\l&:=4,\quad \partial^\dprime_\l:=4\\
\tilde{d}_\l^\prime&:=h^{(2)}_{\l}+10, \quad \tilde{d}_\l^\dprime:=h^{(1)}_{\l}+7
\end{split}
\end{alignat}
\end{subequations}


\subsection{The root system \texorpdfstring{$F_4$}{F4}}\label{II-F4defects}

The Euclidean vector space $\left(\h^*_\Re, (\cdot,\cdot)\right)$ can be identified to $\left(\Re^4, (\cdot,\cdot)\right)$. The simple roots are $\a_i=\e_{i+1}-\e_{i+2}$, $1\leq i\leq 2$, $\a_3=\e_4$, and $\a_4=\frac{1}{2}\e_1-\frac{1}{2}(\e_2+\e_3+\e_4)$. The dominant root is $\th=\e_1+\e_2$. The Weyl group orbit of an element of the root lattice cannot be described concisely in general. For this reason we will split an orbit into orbits under $W_{B_3}$, the parabolic Weyl group of type $B_3$ obtained by excluding the simple root $\a_4$. This group acts only on the last three coordinates of a vector in $\Re^4$ by permuting them and possibly changing their signs. 

The positive roots are $\e_i$, $1\leq i\leq 4$,  $\e_i\pm\e_j$, $1\leq i<j\leq 4$, and the elements of the $W_{B_3}$ orbit of $\a_4$. The dominant non-zero small weights are $\e_1$ (with symbol $[1_s]$),  $\e_1+\e_2$ (with symbol $[1_\ell]$),  and $\frac{3}{2}\e_1+\frac{1}{2}(\e_2+\e_3+\e_4)$  (with symbol $[1_\ell,1_s]$). We will describe separately the defect data for each of these weights. One common feature is that  $\ag_\l\not<{\bf 0}$ if the coefficient of $\e_1$ in $\l$ is strictly negative.

Some of these weights are in fact in the root lattice of the parabolic root subsystem of type $B_3$ mentioned above. For these, it is convenient to use the detect data defined in Section 
\ref{II-Ddefectdata} and it is this data we refer to when we use  notation such as $\de_\l^{B_3}$. 

\subsubsection{Weights with symbol $[1]$} These are  the roots. The defect vector is defined by $\de_\a=0$ if $\a\in R^+$ and $\de_\a=1$ if $\a\in R^-$.

\subsubsection{Weights with symbol $[1^2]$} We describe the  Weyl group orbit of $\frac{3}{2}\e_1+\frac{1}{2}(\e_2+\e_3+\e_4)$  by specifying the $W_{B_3}$ orbits it contains. These are the $W_{B_3}$ orbits of $\frac{3}{2}\e_1+\frac{1}{2}(\e_2+\e_3+\e_4)$, $\frac{1}{2}\e_1+\frac{1}{2}(3\e_2+\e_3+\e_4)$, $\e_1+\e_2+\e_3$, $\e_2+\e_3+\e_4$, $-\frac{3}{2}\e_1+\frac{1}{2}(\e_2+\e_3+\e_4)$, $-\frac{1}{2}\e_1+\frac{1}{2}(3\e_2+\e_3+\e_4)$, and $-\e_1+\e_2+\e_3$. For the elements of the last four $W_{B_3}$ orbits  the coordinate of $\e_1$ is strictly negative and  we set $\ag_\l\not<{\bf 0}$. We focus our attention on the first four orbits. 

For  $\l\in W_{B_3}(\frac{3}{2}\e_1+\frac{1}{2}(\e_2+\e_3+\e_4))$ the defect vector is defined by $\de_\l=(0,0)$. 

For $\l$ in $W_{B_3}(\frac{1}{2}\e_1+\frac{1}{2}(3\e_2+\e_3+\e_4))$, the defect vector is defined by $$\de_\l=(h^{rev, (1)}_{\l}+2h_{\l+\e_1},0)$$ 

For $\l$ in $W_{B_3}(\e_1+\e_2+\e_3)$, the defect vector is defined by $$\de_\l:= (\parallel d_{\l^\natural}\parallel,0)$$

For $\l$ in $W_{B_3}(\e_2+\e_3+\e_4)$, the defect vector is defined by $$\de_\l:= \de_\l^{B_3}+(2,0)$$
Note that $|\A^{-}_3(\e_2+\e_3+\e_4)|=1$ and this is reflected in the formula above.

The defect data presented in the section is nothing else but the one defined in the last paragraph of Section \ref{II-normal}.


\end{document}